\documentclass[11pt]{article}
\usepackage{amsmath, amssymb, epsfig, eepic}

\numberwithin{equation}{section}

\newcommand{\eps}{\varepsilon}

\newcommand{\qed}{\hfill $\Box$}
\newcommand{\proof}{\noindent {\bf Proof}. \hspace{2mm}}
\newcommand{\prob}{\mathbb P}
\newcommand{\expec}{\mathbb E}
\newcommand{\smallsup}[1] {{\scriptscriptstyle{({#1}})}}
\newtheorem{theorem}{Theorem}[section]
\newtheorem{lemma}[theorem]{Lemma}
\newtheorem{prop}[theorem]{Proposition}
\newtheorem{corr}[theorem]{Corollary}
\newtheorem{ass}[theorem]{Assumption}
\newtheorem{remark}[theorem]{Remark}

\newcommand*{\picdirrectory}{}
\newcommand*{\fig}[3]{
        \begin{figure}[!h!t]
    \begin{center}
        \includegraphics{\picdirrectory#1}
    \end{center}
        \caption{\footnotesize #2}
        \label{#3}
        \end{figure}}

\newcommand{\vep}{\varepsilon}
\newcommand{\vepp}{\varepsilon^4}
\newcommand{\eq}{\begin{equation}}
\newcommand{\en}{\end{equation}}
\def\eqalign#1\enalign{
    \begin{align}#1\end{align}
    }
\newcommand{\sss}   { \scriptscriptstyle }
\newcommand{\nn}   { \nonumber}
\newcommand{\Q}{\mathbb{Q}}
\newcommand{\Z}{\mathbb{Z}}
\newcommand{\cZ}{{\cal Z}}
\newcommand{\sN}{{\sss N}}
\newcommand{\sM}{\!{\sss M}}
\newcommand{\sZ}{{\sss Z}}
\newcommand{\shift}   {\!\!\!\!}
\newcommand{\DDcup}{~\Dot{\cup}}

\def\NN{{\mathbb N}}

\newcommand{\cB}{{\cal B}}
\newcommand{\argmin}{{\rm argmin}}
\newcommand{\Ymplus}[1]{Y_{m,+}^{\smallsup{#1,N}}}
\newcommand{\Ymmin}[1]{Y_{m,-}^{\smallsup{#1,N}}}
\newcommand{\Zfactors}{P_m(k,k_1)}
\newcommand{\ktauN}{k_{\tau,\sN}}
\newcommand{\FN}{F_{\vep}}

\newcommand*{\sumu}{\displaystyle\sum}

\def\AeN{ {\cal A}_{\varepsilon\kern-.05em,\kern .1em N} }
\def\BeN{ B_{\varepsilon\kern-.05em,\kern .1em N}}
\def\CeN{ {\cal C}_{\varepsilon\kern-.05em,\kern .1em N} }
\def\DeN{ {\cal D}_{\varepsilon\kern-.05em,\kern .1em N} }
\def\EeN{ {\cal E}_{\varepsilon\kern-.05em,\kern .1em N} }

\def\prob{{\mathbb P}}
\def\expec{{\mathbb E}}
\newcommand{\whpl}{{\bf whp~}}
\newcommand{\whps}{{\bf whp}}

\newcommand{\smfrac}[2]{\textstyle{#1\over #2}}

\title{Distances in random graphs with finite mean and infinite variance degrees}
\author{Remco van der Hofstad\footnote{Department of Mathematics and
Computer Science, Eindhoven University of Technology, P.O.\ Box
513, 5600 MB Eindhoven, The Netherlands. E-mail: {\tt
rhofstad@win.tue.nl}}\\
Gerard Hooghiemstra\footnote{Delft University of Technology,
Electrical Engineering, Mathematics and Computer Science, P.O. Box
5031, 2600 GA Delft, The Netherlands. E-mail: {\tt
G.Hooghiemstra@ewi.tudelft.nl}} ~and Dmitri
Znamenski\footnote{EURANDOM, P.O.\ Box 513, 5600 MB Eindhoven, The
Netherlands. E-mail: {\tt znamenski@eurandom.nl}} }

\begin{document}
\maketitle

\begin{abstract}
In this paper we study random graphs with independent and
identically distributed degrees of which the tail of the
distribution function is regularly varying with exponent $\tau\in
(2,3)$.

The number of edges between two arbitrary nodes, also called the
graph distance or hopcount, in a graph with $N$ nodes is
investigated when $N\rightarrow \infty$. When $\tau\in (2,3)$,
this graph distance grows like $2\frac{\log\log
N}{|\log(\tau-2)|}$. In different papers, the cases $\tau>3$ and
$\tau\in (1,2)$ have been studied. We also study the fluctuations
around these asymptotic means, and describe their distributions.
The results presented here improve upon results of Reittu and
Norros, who prove an upper bound only.
\end{abstract}

\smallskip

AMS 1991 {\it subject classifications.} Primary 05C80; secondary
60J80.



\section{Introduction}
The study of complex networks plays an increasingly important role
in science. Examples of complex networks are electrical power
grids and telephony networks, social relations, the World-Wide Web
and Internet, co-authorship and citation networks of scientists,
etc. The structure of networks affects their performance and
function. For instance, the topology of social networks affects
the spread of information and infections. Measurements on complex
networks have shown that many networks have similar properties. A
first key example of such a fundamental network property is the
fact that typical distances between nodes are small, which is
called the {\it `small world' phenomenon}. A second key example
shared by many networks is that the number of nodes with degree
$k$ falls off as an inverse power of $k$, which is called a {\it
power law degree sequence}. See \cite{AB02, Newm03, Watts} and the
references therein for an introduction to complex networks and
many examples where the above two properties hold.

The current paper presents a rigorous derivation for the random
fluctuations of the graph distance between two arbitrary nodes
(also called the geodesic, and in Internet called the hopcount) in
a graph with infinite variance degrees. The model studied here is
a variant of the {\it configuration model}. The infinite variance
degrees include power laws with exponent $\tau \in (2,3)$. In
practice, power exponents are observed ranging between $\tau=1.5$
and $\tau=3.2$ (see \cite{Newm03}).

In a previous paper of the first two authors with Van Mieghem
\cite{HHV03}, we investigated the finite variance case $\tau>3$.
In \cite{HHZ04b}, we study the case where $\tau \in (1,2)$. Apart
from the critical cases $\tau=2$ and $\tau=3$, we have thus
investigated all possible values of $\tau$. The paper
\cite{HHZ04c} serves as a survey to the results and, in
particular, describes how our results can be applied to Internet
data, describes related work on random graphs that are similar,
though not identical to ours, and gives further open problems.
Finally, in \cite{HHZ04c}, we also investigate the structure of
the connected components in the random graphs under consideration.
See \cite{AS00,Boll01,JLR00} for an introduction to classical
random graphs.

This section is organised as follows. In Section \ref{sec-mod} we
start by introducing the model, in Section \ref{sec-res} we state
our main results. Section \ref{sec-RW} is devoted to related work,
and in Section \ref{sec-sim}, we describe some simulations for a
better understanding of the results.


\subsection{Model definition}
\label{sec-mod} Fix an integer $N$. Consider an i.i.d.\ sequence
$D_1,D_2,\ldots,D_{\sN}$. We will construct an undirected graph
with $N$ nodes where node $j$ has degree $D_j$. We will assume
that $L_{\sN}=\sum_{j=1}^N D_j$ is even. If $L_{\sN}$ is odd, then
we increase $D_{\sN}$ is by 1. This change will make hardly any
difference in what follows, and we will ignore this effect. We
will later specify the distribution of $D_1$.

To construct the graph, we have $N$ separate nodes and incident to
node $j$, we have $D_j$ stubs. All stubs need to be connected to
another stub to build the graph. The stubs are numbered in an
arbitrary order from $1$ to $L_{\sN}$. We start by connecting at
random the first stub with one of the $L_{\sN}-1$ remaining stubs.
Once paired, two stubs form a single edge of the graph. We
continue the procedure of randomly choosing and pairing the stubs
until all stubs are connected. Unfortunately, nodes having
self-loops may occur. However, self-loops are scarce when $N \to
\infty$.

The above model is a variant of the configuration model, which,
given a degree sequence, is the random graph with that given
degree sequence. For a graph, the degree sequence of that graph is
the vectors of which the $k^{\rm th}$ coordinate equals the
frequency of nodes with degree $k$. In our model, the degree
sequence is very close to the distribution of the nodal degree $D$
of which $D_1, \ldots, D_{\sN}$ are i.i.d.\ copies. The
probability mass function and the distribution function of the
nodal degree law are denoted by
    \begin{equation}
    \label{kansen}
    \prob(D_1=j)=f_j,\quad j=1,2,\ldots, \quad \mbox{and} \quad
    F(x)=\sum_{j=1}^{\lfloor x \rfloor} f_j,
    \end{equation}
where $\lfloor x \rfloor$ is the largest integer smaller than or
equal to $x$. Our main assumption is that we take
    \begin{equation}
    \label{distribution}
    1-F(x)=x^{-\tau+1}L(x),
    \end{equation}
where $\tau \in (2,3)$ and $L$ is slowly varying at infinity. This
means that the random variables $D_i$ obey a power law, and the
factor $L$ is meant to generalize the model. We work under a
slightly more restrictive assumption:

\begin{ass}
\label{ass-gamma} There exists $\gamma\in [0,1)$ and $C>0$ such
that
    \eq
    \label{Fcond}
    x^{-\tau+1-C(\log{x})^{\gamma-1}}\leq 1-F(x)\leq
    x^{-\tau+1+C(\log{x})^{\gamma-1}},\qquad \mbox{for large $x$}.
    \en
\end{ass}
Comparing with (\ref{distribution}), we see that the slowly
varying function $L$ in (\ref{distribution}) should satisfy
    \eq
    e^{-C (\log{x})^{\gamma}}\leq L(x)\leq e^{C (\log{x})^{\gamma}}.
    \en


\subsection{Main results}
\label{sec-res}

We define the graph distance $H_{\sN}$ between the nodes $1$ and
$2$ as the minimum number of edges that form a path from $1$ to
$2$. By convention, the distance equals $\infty$ if $1$ and $2$
are not connected. Observe that the distance between two randomly
chosen nodes is equal in distribution to $H_{\sN}$, because the
nodes are exchangeable. We now describe our main result.

\begin{theorem}[Fluctuations of the Graph Distance]
\label{thm-tau(2,3)} Assume that Assumption \ref{ass-gamma} holds
and fix $\tau \in (2,3)$ in (\ref{distribution}). Then there exist
random variables $(R_{a})_{a\in (-1,0]}$ such that, as $N\to
\infty$,
    \begin{eqnarray}
    \label{limit law}
    &&\prob\Big(H_{\sN}=2\Big\lfloor \frac{\log\log N}{|\log (\tau -2)|}
    \Big\rfloor+l~\Big|H_{\sN}<\infty \Big)=\prob(R_{a_{\sN}}=l)+o(1), \qquad l\in {\mathbb Z},
    \end{eqnarray}
where $a_{\sN}=\lfloor \frac{\log\log{N}}{|\log (\tau-2)|}\rfloor
-\frac{\log\log{N}}{|\log (\tau-2)|}\in (-1,0]$.
\end{theorem}

In words, Theorem \ref{thm-tau(2,3)} states that for
$\tau\in(2,3)$, the graph distance $H_{\sN}$ between two randomly
chosen connected nodes grows proportional to $\log\log$ of the
size of the graph, and that the fluctuations around this mean
remain uniformly bounded in $N$.

We identify the laws of $(R_{a})_{a\in (-1,0]}$ below. Before
doing so, we state two consequences of the above theorem:

\begin{corr}[Convergence in Distribution along Subsequences]
\label{cor-weak} Along the sequence $N_k=\lfloor
N_1^{(\tau-2)^{-(k-1)}}\rfloor,$ where $k=1,2, \ldots,$ and
conditionally on 1 and 2 being connected, the random variables
    \eq
    H_{N_k}- 2\Big\lfloor \frac{\log\log {N_k}}{|\log (\tau
    -2)|}\Big\rfloor,
    \en
converge in distribution to $R_{a_{N_1}}$, as $k\rightarrow
\infty$.
\end{corr}

Simulations illustrating the weak convergence in Corollary
\ref{cor-weak} are discussed in Section \ref{sec-sim}. In the
corollary below, we write that an event $E$ occurs \whpl for the
statement that $\prob(E)=1-o(1)$.

\begin{corr}[Concentration of the Graph Distance]~
\begin{itemize}
\item[{\rm (i)}] Conditionally on 1 and 2 being connected, the
random variable $H_{\sN}$ is, \whps, in between $2\frac{\log\log
{N}}{|\log (\tau -2)|}(1\pm \vep)$, for any $\vep>0$. \item[{\rm
(ii)}] Conditionally on 1 and 2 being connected, the random
variables $H_{\sN}- \frac{\log\log {N}}{|\log (\tau -2)|}$ form a
tight sequence, i.e.,
    \begin{equation}
    \label{hoptight}
    \lim_{K\rightarrow \infty}
    \limsup_{N\rightarrow \infty}
    \prob\Big(\big|H_{\sN}- 2\frac{\log\log {N}}{|\log (\tau -2)|}\big|\leq K~\Big|H_{\sN}<\infty\Big)
    =1.
    \end{equation}
\end{itemize}
\end{corr}

We need a limit result from branching processes theory before we
can identify the limiting random variables $(R_{a})_{a\in
(-1,0]}$. In Section \ref{sec-BP}, we introduce a {\it delayed}
branching process $\{ {\cal Z}_k\}_{k\geq 1}$, where in the first
generation the offspring distribution is chosen according to
(\ref{kansen}) and in the second and further generations the
offspring is chosen in accordance to $g$ given by
    \begin{equation}
    \label{outgoing degree}
        g_j=\frac{(j+1) f_{j+1}}{\mu},\quad j=0,1,\ldots,
    \end{equation}
where $\mu=\sum_{j=1}^{\infty} jf_j$. The branching process
$\{{\cal Z}_k\}$ has infinite expectation. Branching processes
with infinite expectation have been investigated in \cite{davies,
seneta2, barbour}. Assumption \ref{ass-gamma}, using the results
in \cite{davies}, implies that
    \begin{equation}
    \label{loglaw}
    (\tau-2)^n\cdot\log({\cal Z}_n\vee 1)\to Y,\qquad \mbox{a.s.},
    \end{equation}
where $x\vee y=\max\{x,y\}.$ See Section \ref{sec-BP} and the
references there for more details. Then, we can identify the law
of the random variables $(R_{a})_{a\in (-1,0]}$ as follows:

\begin{theorem}[The Limit Laws]
\label{thm-ll} For $a\in (-1,0]$,
    \eq
    \label{Radef}
    \prob(R_{a}>l)=\prob\Big(\min_{s\in {\mathbb Z}}
    \big[(\tau-2)^{-s} Y^{\smallsup{1}}+
    (\tau-2)^{s-c_l} Y^{\smallsup{2}}\big]
    \leq (\tau-2)^{\lceil l/2 \rceil+a}
    \big|Y^\smallsup{1}Y^\smallsup{2}>0\Big),\nonumber
    \en
where $c_l=1$ if $l$ is even, and zero otherwise, and
$Y^{\smallsup{1}}, Y^{\smallsup{2}}$ are two independent copies of
the limit random variable in (\ref{loglaw}).
\end{theorem}

In Remarks \ref{rem-confmod} and \ref{rem-conf} below, we will
explain that our results also apply to the usual configuration
model, where the number of nodes with a given degree is fixed,
when we study the graph distance between two uniformly chosen
nodes, and the degree distribution satisfied certain conditions.
For the precise conditions, see Remark \ref{rem-conf}.

\subsection{Related work}
\label{sec-RW} There is a wealth of related work which we now
summarize. The model investigated here was also studied in
\cite{norros}, with $1-F(x)=x^{-\tau+1}L(x),$ where $\tau\in
(2,3)$ and $L$ denotes a slowly varying function. It was shown in
\cite{norros} that \whpl the graph distance is bounded from above
by $2\frac{\log\log N}{|\log(\tau-2)|}(1+o(1))$. We improve the
results in \cite{norros} by deriving the asymptotic distribution
of the random fluctuations of the graph distance around $2\lfloor
\frac{\log\log N}{|\log(\tau-2)|}\rfloor$. Note that these results
are in contrast to \cite[Section II.F, below (56)]{NSW00}, where
it was suggested that if $\tau<3$, then an exponential cut-off is
necessary to make the graph distance between an arbitrary pair of
nodes well-defined. The problem of the graph distance between an
arbitrary pair of nodes was also studied non-rigorously in
\cite{CH03}, where also the behavior when $\tau=3$ and $x\mapsto
L(x)$ is the constant function, is included. In the latter case,
the graph distance scales like $\frac{\log{N}}{\log\log N}$. A
related model to the one studied here can also be found in
\cite{NR04}, where a graph process is defined by adding and
removing edges. In \cite{NR04}, the authors prove similar results
as in \cite{norros} for this related model.

The graph distance for $\tau>3$, $\tau \in (1,2)$, respectively
was treated in two previous publications \cite{HHV03} and
\cite{HHZ04b}, respectively. We survey these results together with
results on the connected components in \cite{HHZ04c}. In
\cite{HHZ04c}, we also show that when $\tau>2$, the diameter is
bounded from below by a constant times $\log{N}$, which, when
$\tau\in(2,3)$ should be contrasted with the average graph
distance, which is or order $\log\log{N}$. Finally, in
\cite{HHZ04c} also the connected components are studied under the
condition that $\mu=\expec[D_1]>2$, and the results in this paper
are used to show that \whpl there exists a largest connected
component of size $qN[1+o(1)]$, where $q$ is the survival
probability of the delayed branching process, while all other
connected components are of order at most $\log{N}$.

There is substantial work on random graphs that are, although
different from ours, still similar in spirit. In \cite{ACL01},
random graphs were considered with a degree sequence that is {\it
precisely} equal to a power law, meaning that the number of nodes
with degree $k$ is precisely proportional to $k^{-\tau}$. Aiello
{\em et al.} \cite{ACL01} show that the largest connected
component is of the order of the size of the graph when
$\tau<\tau_0=3.47875\ldots$, where $\tau_0$ is the solution of
$\zeta(\tau-2)-2\zeta(\tau-1)=0$, and where $\zeta$ is the Riemann
zeta function. When $\tau>\tau_0$, the largest connected component
is of smaller order than the size of the graph and more precise
bounds are given for the largest connected component. When
$\tau\in (1,2)$, the graph is \whpl connected. The proofs of these
facts use couplings with branching processes and strengthen
previous results due to Molloy and Reed \cite{MR95,MR98}. For this
same model, Dorogovtsev {\em et al.} \cite{DGM01, DGM02}
investigate the leading asymptotics and the fluctuations around
the mean of the graph distance between arbitrary nodes from a
theoretical physics point of view, using mainly generating
functions.

A second related model can be found in \cite{CL02a, CL02b}, where
edges between nodes $i$ and $j$ are present with probability equal
to $w_iw_j/\sum_l w_l$ for some `expected degree vector' $w=(w_1,
\ldots, w_{\sN})$. It is assumed that $\max_i w_i^2<\sum_i w_i$,
so that $w_iw_j/\sum_l w_l$ are probabilities. In \cite{CL02a},
$w_i$ is often taken as $w_i=c i^{-{\frac 1{\tau-1}}}$, where $c$
is a function of $N$ proportional to $N^{\frac{1}{\tau-1}}$. In
this case, the degrees obey a power law with exponent $\tau$.
Chung and Lu \cite{CL02a} show that in this case, the graph
distance between two uniformly chosen nodes is \whpl proportional
to $\log N(1+o(1))$ when $\tau>3$, and $2\frac{\log\log
N}{|\log(\tau-2)|}(1+o(1))$ when $\tau\in (2,3)$. The difference
between this model and ours is that the nodes are not exchangeable
in \cite{CL02a}, but the observed phenomena are similar. This
result can be heuristically understood as follows. Firstly, the
actual degree vector in \cite{CL02a} should be close to the
expected degree vector. Secondly, for the expected degree vector,
we can compute that the number of nodes for which the degree is at
least $k$ equals
    $$
    |\{i: w_i\geq k\}|=|\{i: ci^{-\frac{1}{\tau-1}}\geq k\}|\propto k^{-\tau+1}.
    $$
Thus, one expects that the number of nodes with degree at least
$k$ decreases as $k^{-\tau+1}$, similarly as in our model. In
\cite{CL02b}, Chung and Lu study the sizes of  the connected
components in the above model. The advantage of this model is that
the edges are {\it independently} present, which makes the
resulting graph closer to a traditional random graph.

All the models described above are {\it static}, i.e., the size of
the graph is {\it fixed}, and we have not modeled the {\it growth}
of the graph. There is a large body of work investigating {\it
dynamical} models for complex networks, often in the context of
the World-Wide Web. In various forms, preferential attachment has
been shown to lead to power law degree sequences. Therefore, such
models intend to {\it explain} the occurrence of power law degree
sequences in random graphs. See \cite{ACL01b, AB99, AB02, BBCR03,
BR02, BR03a, BR03b, BRST01, CF03, KRRSTU00} and the references
therein. In the preferential attachment model, nodes with a fixed
degree $m$ are added sequentially. Their stubs are attached to a
receiving node with a probability proportional to the degree of
the receiving node, thus favoring nodes with large degrees. For
this model, it is shown that the number of nodes with degree $k$
decays proportionally to $k^{-3}$ \cite{BRST01}, the diameter is
of order $\frac{\log{N}}{\log\log{N}}$ when $m\geq 2$ \cite{BR02},
and couplings to a classical random graph $G(N,p)$ are given for
an appropriately chosen $p$ in \cite{BR03b}. See also \cite{BR03a}
for a survey.

Possibly, the configuration model is a snapshot of the above
models, i.e., a realization of the graph growth processes at the
time instant that the graph has a certain prescribed size. Thus,
rather than to describe the growth of the model, we investigate
the properties of the model at a given time instant. This is
suggested in \cite[Section VII.D]{AB02}, and it would be very
interesting indeed to investigate this further mathematically,
i.e., to investigate the relation between the configuration and
the preferential attachment models.

\renewcommand{\epsfsize}[2]{0.5#1}

\begin{figure}[t]
\begin{center}
\epsfbox[20 40 576 469]{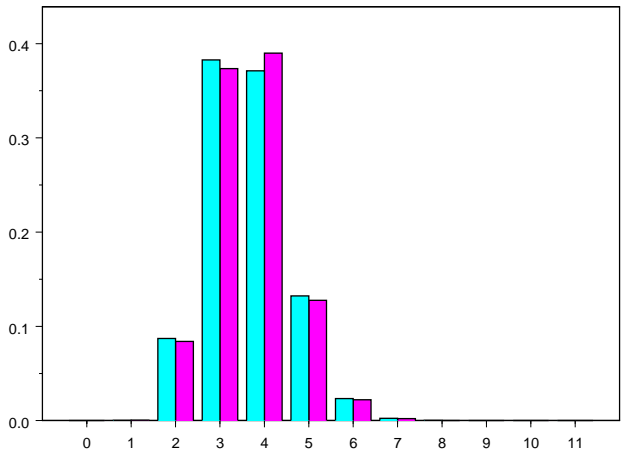} \caption{Histograms of the
AS-count and graph distance in the configuration model with
$N=10,940$, where the degrees have generating function $f_\tau(s)$
in (\ref{genfunctionfirst}), for which the power law exponent
$\tau$ takes the value $\tau=2.25$. The AS-data is lightly shaded,
the simulation is darkly shaded.} \label{fig-AScount}
\end{center}
\end{figure}

We study the above version of the configuration model to describe
the topology of the Internet at a fixed time instant. In a seminal
paper \cite{FFF99}, Faloutsos {\em et al.} have shown that the
degree distribution in Internet follows a power law with exponent
$\tau\approx 2.16-2.25$. Thus, the power law random graph with
this value of $\tau$ can possibly lead to a good Internet model.
In \cite{Tangmunarunkit_sigcom02}, and inspired by the observed
power law degree sequence in \cite{FFF99}, the power law random
graph is proposed as a model for the network of {\it autonomous
systems}. In this graph, the nodes are the autonomous systems in
the Internet, i.e., the parts of the Internet controlled by a
single party (such as a university, company or provider), and the
edges represent the physical connections between the different
autonomous systems. The work of Faloutsos {\em et al.} in
\cite{FFF99} was among others on this graph which at that time had
size approximately 10,000.

In \cite{Tangmunarunkit_sigcom02}, it is argued on a qualitative
basis that the power law random graph serves as a better model for
the Internet topology than the currently used topology generators.
Our results can be seen as a step towards the quantitative
understanding of whether the AS-count in Internet is described
well by the average graph distance in the configuration model. The
AS-count gives the number of physical links connecting the various
autonomous domains between two randomly chosen nodes in the graph.

To validate the model, we compare a simulation of the distribution
of the distance between pairs of nodes in the power law random
graph with the same value of $N$ and $\tau$ to extensive
measurements of the AS-count in Internet. In Figure
\ref{fig-AScount}, we see that AS-count in the model with the
predicted value of $\tau=2.25$ and the value of $N$ from the data
set fits the data remarkably well.

In \cite[Table II]{Newm03}, many other examples are given of real
networks that have power law degree sequences. Interestingly,
there are many examples where the power law exponent is in
$(2,3)$, and it would be of interest to compare the average graph
distance between an arbitrary pair of nodes in such examples.

\subsection{Demonstration of Corollary \ref{cor-weak}}
\label{sec-sim} By a simulation we explain the relevance of
Theorem \ref{thm-tau(2,3)} and especially the relevance of
Corollary \ref{cor-weak}. We have chosen to simulate the
distribution (\ref{outgoing degree}) from the generating function:
    \begin{equation}
    \label{genfunctionfurther}
    g_{\tau}(s)=1-(1-s)^{\tau-2},\qquad \text{ for which }
    \qquad
    g_j=(-1)^{j-1}{{\tau-2} \choose j}\sim \frac{c}{j^{\tau-1}}, \quad j\to \infty.
    \end{equation}
Defining
    \begin{equation}
    \label{genfunctionfirst}
    f_{\tau}(s)=\frac{\tau-1}{\tau-2}s-\frac{1-(1-s)^{\tau-1}}{{\tau-2}},\quad
    \tau \in (2,3),
    \end{equation}
it is immediate that
    $$
    g_{\tau}(s)=\frac{f'_{\tau}(s)}{f'_{\tau}(1)},
    \qquad \mbox{ so that }\qquad g_j=\frac{(j+1)f_{j+1}}{\mu}.
    $$

For fixed $\tau$, we can pick different values of the size of the
simulated graph, so that for each two simulated values $N$ and $M$
we have $a_{\sN}=a_{\sss M}$, i.e., $N=M^{(\tau-2)^{-k}},$ for
some integer $k$. For $\tau=2.8$, we have taken the values
    $$
    N=1,000,\qquad N=5,623,\qquad N=48,697,\qquad N=723,394.
    $$
According to Theorem \ref{thm-tau(2,3)}, the survival functions of
the hopcount $H_{\sN}$, satisfying $N=M^{(\tau-2)^{-k}},$, run
parallel on distance $2$ in the limit for $N\to \infty$. In
Section \ref{sec-BP} below we will show that the distribution with
generating function (\ref{genfunctionfirst}) satisfies Assumption
\ref{ass-gamma}.


\renewcommand{\epsfsize}[2]{0.5#1}
\begin{figure}[t]
\begin{center}
\epsfbox[40 60 556 449]{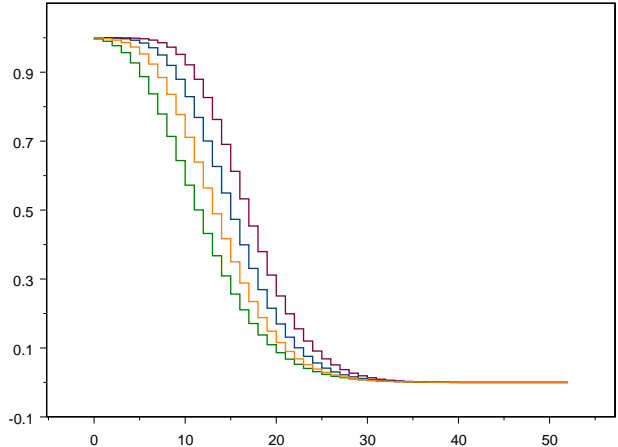} \caption{Empirical survival
functions of the graph distance for $\tau=2.8$ and for the four
values of $N$.}
\label{survival plots 2}
\end{center}
\end{figure}

\subsection{Organization of the paper}
\label{sec-org} The paper is organized as follows. We first review
the relevant literature on branching processes in Section
\ref{sec-BP}. We then describe the growth of shortest path graphs
in Section \ref{subsec-BP1}, and we state coupling results needed
to prove our main results, Theorems
\ref{thm-tau(2,3)}--\ref{thm-ll} in Section \ref{sec-pftau(2,3)}.
In Section \ref{sec-pflemmas}, we prove three technical lemmas
used in Section \ref{sec-pftau(2,3)}. We finally prove the
coupling results in the Appendix.

\section{Review of branching process theory with infinite mean}
\label{sec-BP} Since we heavily rely on the theory of branching
processes (BP's), we now briefly review this theory in the case
where the expected value (mean) of the offspring distribution is
infinite. We follow in particular \cite{davies}, and also refer
the readers to related work in \cite{barbour, seneta2}, and the
references therein.

For the formal definition of the BP we define a double sequence
$\{X_{n,i}\}_{n\geq 0, i\geq 1}$ of i.i.d.\ random variables each
with distribution equal to the offspring distribution $\{g_j\}$
given in (\ref{outgoing degree}) with distribution function
$G(x)=\sum_{j=0}^{\lfloor x\rfloor} g_j$. The BP $\{{\cal Z}_n\}$
is now defined by ${\cal Z}_0=1$ and
    $$
    {\cal Z}_{n+1}=\sum_{i=1}^{{\cal Z}_n} X_{n,i},\quad n\ge 0.
    $$
In case of a delayed BP, we let $X_{0,1}$ have probability mass
function $\{f_j\}$, independently of $\{X_{n,i}\}_{n\geq 1, i\geq
1}$. In this section we restrict to the non-delayed case for
simplicity.

We follow Davies in \cite{davies}, who gives the following
sufficient conditions for convergence of\linebreak $(\tau-2)^n
\log(1+\cZ_n)$. Davies' main theorem states that if for some
non-negative, non-increasing function $\gamma(x)$:
\begin{enumerate}
\item[(i)] $ x^{-\alpha-\gamma(x)}\leq 1-G(x)\leq
x^{-\alpha+\gamma(x)},\qquad \mbox{for large $x$}\quad \mbox{and
\quad} 0<\alpha<1, $ \item[(ii)] $x^{\gamma(x)}$ is
non-decreasing, \item[(iii)] $\int_0^\infty
\gamma(e^{e^x})\,dx<\infty$ or equivalently $\int_e^\infty
\frac{\gamma(y)}{y\log y}\,dy <\infty$,
\end{enumerate}
then $\alpha^n\log(1+\cZ_n)$ converges almost surely to a
non-degenerate finite random variable $Y$ with $\prob(Y=0)$ equal
to the extinction probability of $\{\cZ_n\}$, whereas $Y$ admits a
density on $(0,\infty)$. Therefore, also $\alpha^n \log(\cZ_n\vee
1)$ converges to $Y$ almost surely.

The conditions of Davies quoted as (i-iii) simplify earlier work
by Seneta \cite{seneta2}. For example, for $\{g\}$ in
(\ref{genfunctionfurther}), the above is valid with
$\alpha=\tau-2$ and $\gamma(x) = C(\log{x})^{-1}$, where $C$ is
sufficiently large. We prove in Lemma \ref{lem-G-x} below that for
$F$ as in Assumption \ref{ass-gamma}, and $G$ the distribution
function of $g$ in (\ref{outgoing degree}), the conditions (i-iii)
are satisfied with $\alpha=\tau-2$ and $\gamma(x) =
C(\log{x})^{\gamma-1}$. In particular, for (iii), we need that
$\gamma<1$.

Let $Y^{\smallsup{1}}$ and $Y^{\smallsup{2}}$ be two independent
copies of the limit random variable $Y$. In the course of the
proof, we will encounter the random variable $M=\min_{t\in \Z}
(\kappa^{t} Y^{\smallsup{1}}+\kappa^{c-t} Y^{\smallsup{2}})$, for
some $c\in \{0,1\}$, and where $\kappa=(\tau-2)^{-1}$. The proof
relies on the fact that, conditionally on
$Y^{\smallsup{1}}Y^{\smallsup{2}}>0$, $M$ has a density. The proof
of this fact is as follows. The function $(y_1,y_2)\mapsto
\min_{t\in \Z} (\kappa^{t} y_1 +\kappa^{c-t} y_2)$ is
discontinuous precisely in the points $(y_1,y_2)$ satisfying
$\sqrt{y_2/y_1}=\kappa^{n+\frac12 c},\, n \in \Z$, and,
conditionally on $Y^{\smallsup{1}}Y^{\smallsup{2}}>0$, the random
variables $Y^{\smallsup{1}}$ and $Y^{\smallsup{2}}$ are
independent continuous random variables. Therefore, conditionally
on $Y^{\smallsup{1}}Y^{\smallsup{2}}>0$, the random variable
$M=\min_{t\in \Z} (\kappa^{t} Y^{\smallsup{1}}+\kappa^{c-t}
Y^{\smallsup{2}})$ has a density.

\section{The growth of the shortest path graph}
\label{subsec-BP1}

In this section, we describe the growth of the shortest path graph
(SPG). As a result, we will see that this growth is closely
related to a BP $\{{\hat Z}^{\smallsup{1,N}}_k\}$ with the {\it
random} offspring distribution $\{g_j^{\smallsup{N}}\}$ given by
    \begin{eqnarray}
    g_j^{\smallsup{N}}&=& \sum_{i=1}^N {\bf 1}_{\{D_i=j+1\}} \prob(\mbox{a stub from node
    $i$ is sampled}|D_1, \ldots, D_{\sN}) \nonumber \\
    &=& \sum_{i=1}^N {\bf 1}_{\{D_i=j+1\}}\frac{D_i}{L_{\sN}}=
    \frac{j+1}{L_{\sN}}\sum_{i=1}^N {\bf 1}_{\{D_i=j+1\}} \label{gnN},
    \end{eqnarray}
where, for an event $A$, ${\bf 1}_A$ denotes the indicator
function of the event $A$. By the strong law of large numbers for
$N\to \infty$, almost surely,
    $$
    \frac{L_{\sN}}{N}\to \expec[D],\qquad \mbox{and}\qquad \frac1N \sum_{i=1}^N
    {\bf 1}_{\{D_i=j+1\}}\to \prob(D=j+1),
    $$
so that a.s.,
    \begin{equation}
    \label{convergenceoffspring}
    g_j^{\smallsup{N}} \rightarrow (j+1) \prob(D=j+1)/\expec[D]=g_j, \qquad\qquad N\to \infty.
    \end{equation}
Therefore, the BP $\{{\hat Z}_k^{\smallsup{1,N}}\}$, with
offspring distribution $\{g_j^{\smallsup{N}}\}$, is expected to be
close to a BP with offspring distribution $\{g_j\}$ given in
(\ref{outgoing degree}). Consequently, in Section
\ref{sec-couplingresults}, we state bounds on the coupling of the
BP $\{{\hat Z}^{\smallsup{1,N}}_k\}$ to a BP $\{{\cal
Z}^{\smallsup{1}}_k\}$ with offspring distribution $\{g_j\}$. This
allows us to prove Theorems \ref{thm-tau(2,3)} and \ref{thm-ll} in
Section \ref{sec-pftau(2,3)}.

The shortest path graph (SPG) from node 1 is the power law random
graph as observed from node 1, and consists of the shortest paths
between node 1 and all other nodes $\{2, \ldots, N\}$. As will be
shown below, the SPG is not necessarily a tree because cycles may
occur. Recall that two stubs together form an edge. We define
$Z^{\smallsup{1,N}}_1=D_1$ and, for $k\ge 2$, we denote by
$Z^{\smallsup{1,N}}_k$ the number of stubs attached to nodes at
distance $k-1$ from node 1, but are not part of an edge connected
to a node at distance $k-2$. We refer to such stubs as `free
stubs'.  Thus, $Z^{\smallsup{1,N}}_k$ is the number of outgoing
stubs from nodes at distance $k-1$. By ${\rm SPG}_{k-1}$ we denote
the SPG up to level $k-1$, i.e., up to the moment we have
$Z^{\smallsup{1,N}}_k$ free stubs attached to nodes on distance
$k-1$, and no stubs to nodes on distance $k$. Since we compare
$Z^{\smallsup{1,N}}_k$ to the $k^{\rm th}$ generation of the BP
${\hat Z}^{\smallsup 1}_k$, we call $Z^{\smallsup{1,N}}_k$ the
stubs of level $k$.

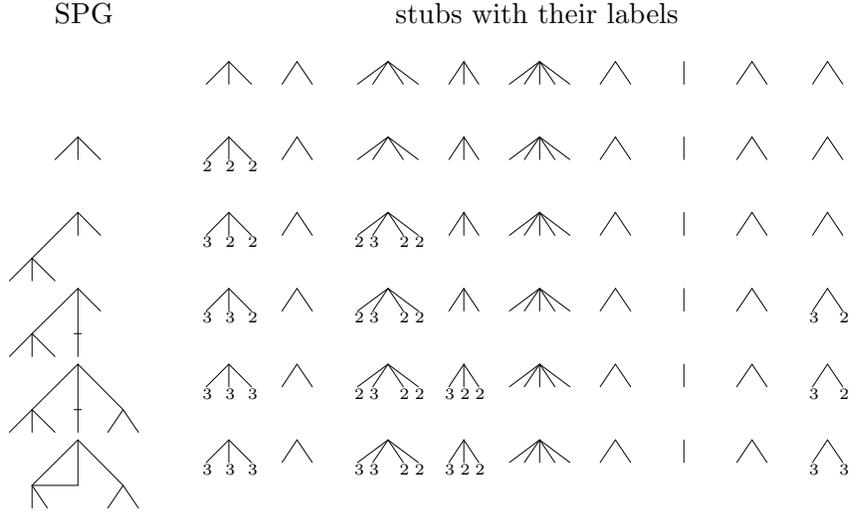
\begin{figure}[t]
\begin{center}
\setlength{\unitlength}{0.0004in}
{
\begin{picture}(6000,9000)(-1000,-6000)

\put(-2000,800) {\makebox(0,0)[lb]{SPG}}

\put(2500,800) {\makebox(0,0)[lb]{stubs with their labels}}

\path(0,0)(300,300)(600,0) \path(300,0)(300,300)

\path(1000,0)(1200,300)(1400,0)

\path(2000,0)(2400,300)(2800,0) \path(2200,0)(2400,300)(2600,0)

\path(3200,0)(3400,300)(3600,0) \path(3400,0)(3400,300)

\path(4000,0)(4400,300)(4800,0) \path(4200,0)(4400,300)(4600,0)
\path(4400,0)(4400,300)

\path(5200,0)(5400,300)(5600,0)

\path(6300,0)(6300,300)

\path(7000,0)(7200,300)(7400,0)

\path(8000,0)(8200,300)(8400,0)


\path(-2000,-1000)(-1700,-700)(-1400,-1000)
\path(-1700,-1000)(-1700,-700)

\path(0,-1000)(300,-700)(600,-1000) \path(300,-1000)(300,-700)
\put(-50,-1150){\makebox(0,0)[lb]{${\scriptscriptstyle 2}$}}
\put(250,-1150){\makebox(0,0)[lb]{${\scriptscriptstyle 2}$}}
\put(550,-1150){\makebox(0,0)[lb]{${\scriptscriptstyle 2}$}}

\path(1000,-1000)(1200,-700)(1400,-1000)

\path(2000,-1000)(2400,-700)(2800,-1000)
\path(2200,-1000)(2400,-700)(2600,-1000)

\path(3200,-1000)(3400,-700)(3600,-1000)
\path(3400,-1000)(3400,-700)

\path(4000,-1000)(4400,-700)(4800,-1000)
\path(4200,-1000)(4400,-700)(4600,-1000)
\path(4400,-1000)(4400,-700)

\path(5200,-1000)(5400,-700)(5600,-1000)

\path(6300,-1000)(6300,-700)

\path(7000,-1000)(7200,-700)(7400,-1000)

\path(8000,-1000)(8200,-700)(8400,-1000)

\path(-2300,-2300)(-1700,-1700)(-1400,-2000)
\path(-1700,-2000)(-1700,-1700) \path(-2300,-2300)(-2300,-2600)
\path(-2600,-2600)(-2300,-2300)(-2000,-2600)

\path(0,-2000)(300,-1700)(600,-2000) \path(300,-2000)(300,-1700)

\put(-50,-2150){\makebox(0,0)[lb]{${\scriptscriptstyle 3}$}}
\put(250,-2150){\makebox(0,0)[lb]{${\scriptscriptstyle 2}$}}
\put(550,-2150){\makebox(0,0)[lb]{${\scriptscriptstyle 2}$}}

\path(1000,-2000)(1200,-1700)(1400,-2000)

\path(2000,-2000)(2400,-1700)(2800,-2000)
\path(2200,-2000)(2400,-1700)(2600,-2000)

\put(1950,-2150){\makebox(0,0)[lb]{${\scriptscriptstyle 2}$}}
\put(2150,-2150){\makebox(0,0)[lb]{${\scriptscriptstyle 3}$}}
\put(2550,-2150){\makebox(0,0)[lb]{${\scriptscriptstyle 2}$}}
\put(2750,-2150){\makebox(0,0)[lb]{${\scriptscriptstyle 2}$}}

\path(3200,-2000)(3400,-1700)(3600,-2000)
\path(3400,-2000)(3400,-1700)

\path(4000,-2000)(4400,-1700)(4800,-2000)
\path(4200,-2000)(4400,-1700)(4600,-2000)
\path(4400,-2000)(4400,-1700)

\path(5200,-2000)(5400,-1700)(5600,-2000)

\path(6300,-2000)(6300,-1700)

\path(7000,-2000)(7200,-1700)(7400,-2000)

\path(8000,-2000)(8200,-1700)(8400,-2000)

\path(-2300,-3300)(-1700,-2700)(-1400,-3000)
\path(-1700,-3600)(-1700,-2700) \path(-1750,-3300)(-1650,-3300)
\path(-2300,-3300)(-2300,-3600)
\path(-2600,-3600)(-2300,-3300)(-2000,-3600)

\path(0,-3000)(300,-2700)(600,-3000) \path(300,-3000)(300,-2700)

\put(-50,-3150){\makebox(0,0)[lb]{${\scriptscriptstyle 3}$}}
\put(250,-3150){\makebox(0,0)[lb]{${\scriptscriptstyle 3}$}}
\put(550,-3150){\makebox(0,0)[lb]{${\scriptscriptstyle 2}$}}

\path(1000,-3000)(1200,-2700)(1400,-3000)

\path(2000,-3000)(2400,-2700)(2800,-3000)
\path(2200,-3000)(2400,-2700)(2600,-3000)

\put(1950,-3150){\makebox(0,0)[lb]{${\scriptscriptstyle 2}$}}
\put(2150,-3150){\makebox(0,0)[lb]{${\scriptscriptstyle 3}$}}
\put(2550,-3150){\makebox(0,0)[lb]{${\scriptscriptstyle 2}$}}
\put(2750,-3150){\makebox(0,0)[lb]{${\scriptscriptstyle 2}$}}

\path(3200,-3000)(3400,-2700)(3600,-3000)
\path(3400,-3000)(3400,-2700)

\path(4000,-3000)(4400,-2700)(4800,-3000)
\path(4200,-3000)(4400,-2700)(4600,-3000)
\path(4400,-3000)(4400,-2700)

\path(5200,-3000)(5400,-2700)(5600,-3000)

\path(6300,-3000)(6300,-2700)

\path(7000,-3000)(7200,-2700)(7400,-3000)

\path(8000,-3000)(8200,-2700)(8400,-3000)
\put(7950,-3150){\makebox(0,0)[lb]{${\scriptscriptstyle 3}$}}
\put(8350,-3150){\makebox(0,0)[lb]{${\scriptscriptstyle 2}$}}

\path(-2300,-4300)(-1700,-3700)(-1100,-4300)
\path(-1700,-4600)(-1700,-3700) \path(-1750,-4300)(-1650,-4300)
\path(-2300,-4300)(-2300,-4600)
\path(-1300,-4600)(-1100,-4300)(-900,-4600)
\path(-2600,-4600)(-2300,-4300)(-2000,-4600)


\path(0,-4000)(300,-3700)(600,-4000) \path(300,-4000)(300,-3700)

\put(-50,-4150){\makebox(0,0)[lb]{${\scriptscriptstyle 3}$}}
\put(250,-4150){\makebox(0,0)[lb]{${\scriptscriptstyle 3}$}}
\put(550,-4150){\makebox(0,0)[lb]{${\scriptscriptstyle 3}$}}

\path(1000,-4000)(1200,-3700)(1400,-4000)

\path(2000,-4000)(2400,-3700)(2800,-4000)
\path(2200,-4000)(2400,-3700)(2600,-4000)

\put(1950,-4150){\makebox(0,0)[lb]{${\scriptscriptstyle 2}$}}
\put(2150,-4150){\makebox(0,0)[lb]{${\scriptscriptstyle 3}$}}
\put(2550,-4150){\makebox(0,0)[lb]{${\scriptscriptstyle 2}$}}
\put(2750,-4150){\makebox(0,0)[lb]{${\scriptscriptstyle 2}$}}

\path(3200,-4000)(3400,-3700)(3600,-4000)
\path(3400,-4000)(3400,-3700)
\put(3150,-4150){\makebox(0,0)[lb]{${\scriptscriptstyle 3}$}}
\put(3350,-4150){\makebox(0,0)[lb]{${\scriptscriptstyle 2}$}}
\put(3550,-4150){\makebox(0,0)[lb]{${\scriptscriptstyle 2}$}}

\path(4000,-4000)(4400,-3700)(4800,-4000)
\path(4200,-4000)(4400,-3700)(4600,-4000)
\path(4400,-4000)(4400,-3700)

\path(5200,-4000)(5400,-3700)(5600,-4000)

\path(6300,-4000)(6300,-3700)

\path(7000,-4000)(7200,-3700)(7400,-4000)

\path(8000,-4000)(8200,-3700)(8400,-4000)
\put(7950,-4150){\makebox(0,0)[lb]{${\scriptscriptstyle 3}$}}
\put(8350,-4150){\makebox(0,0)[lb]{${\scriptscriptstyle 2}$}}

\path(-2300,-5300)(-1700,-4700)(-1100,-5300)
\path(-1700,-5300)(-1700,-4700) \path(-2300,-5300)(-2100,-5600)
\path(-1300,-5600)(-1100,-5300)(-900,-5600)
\path(-2300,-5300)(-2300,-5600) \path(-2300,-5300)(-1700,-5300)

\path(0,-5000)(300,-4700)(600,-5000) \path(300,-5000)(300,-4700)

\put(-50,-5150){\makebox(0,0)[lb]{${\scriptscriptstyle 3}$}}
\put(250,-5150){\makebox(0,0)[lb]{${\scriptscriptstyle 3}$}}
\put(550,-5150){\makebox(0,0)[lb]{${\scriptscriptstyle 3}$}}

\path(1000,-5000)(1200,-4700)(1400,-5000)

\path(2000,-5000)(2400,-4700)(2800,-5000)
\path(2200,-5000)(2400,-4700)(2600,-5000)

\put(1950,-5150){\makebox(0,0)[lb]{${\scriptscriptstyle 3}$}}
\put(2150,-5150){\makebox(0,0)[lb]{${\scriptscriptstyle 3}$}}
\put(2550,-5150){\makebox(0,0)[lb]{${\scriptscriptstyle 2}$}}
\put(2750,-5150){\makebox(0,0)[lb]{${\scriptscriptstyle 2}$}}

\path(3200,-5000)(3400,-4700)(3600,-5000)
\path(3400,-5000)(3400,-4700)
\put(3150,-5150){\makebox(0,0)[lb]{${\scriptscriptstyle 3}$}}
\put(3350,-5150){\makebox(0,0)[lb]{${\scriptscriptstyle 2}$}}
\put(3550,-5150){\makebox(0,0)[lb]{${\scriptscriptstyle 2}$}}

\path(4000,-5000)(4400,-4700)(4800,-5000)
\path(4200,-5000)(4400,-4700)(4600,-5000)
\path(4400,-5000)(4400,-4700)

\path(5200,-5000)(5400,-4700)(5600,-5000)

\path(6300,-5000)(6300,-4700)

\path(7000,-5000)(7200,-4700)(7400,-5000)

\path(8000,-5000)(8200,-4700)(8400,-5000)
\put(7950,-5150){\makebox(0,0)[lb]{${\scriptscriptstyle 3}$}}
\put(8350,-5150){\makebox(0,0)[lb]{${\scriptscriptstyle 3}$}}

\end{picture}
}
\end{center}
\caption{Schematic drawing of the growth of the SPG from the node
1 with $N=9$ and the updating of the labels.} \label{fig-1}
\end{figure}

The first stages of a realization of the generation of the SPG,
with $N=9$ and $L_{\sN}=24$, are drawn in Figure \ref{fig-1}. The
first line shows the $N$ different nodes with their attached
stubs. Initially, all stubs have label 1. The growth process
starts by choosing the first stub of node 1 whose stubs are
labeled by 2 as illustrated in the second line, while all the
other stubs maintain the label 1. Next, we uniformly choose a stub
with label 1 or 2. In the example in line 3, this is the second
stub from node 3, whose stubs are labeled by 2 and the second stub
by label 3. The left hand side column visualizes the growth of the
SPG by the attachment of stub 2 of node 3 to the first stub of
node 1. Once an edge is established the pairing stubs are labeled
3. In the next step, again a stub is chosen uniformly out of those
with label 1 or 2. In the example in line 4, it is the first stub
of the last node that will be attached to the second stub of node
1, the next in sequence to be paired. The last line exhibits the
result of creating a cycle when the second stub of the last node
is chosen to be attached to the first stub of node 3, which is the
next stub in the sequence to be paired. This process is continued
until there are no more stubs with labels 1 or 2. In this example,
we have $Z_1^{\smallsup{1,9}}=3$ and $Z_2^{\smallsup{1,9}}=6$.

We now describe the meaning of the labels. Initially, all stubs
are labeled 1. At each stage of the growth of the SPG, we draw
uniformly at random from all stubs with labels 1 and 2. After each
draw we will update the realization of the SPG, and classify the
stubs according to three categories, which will be labeled 1, 2
and 3. These labels will be updated as the growth of the SPG
proceeds. At any stage of the generation of the SPG, the labels
have the following meaning:

\begin{enumerate}
\item[1.] Stubs with label 1 are stubs belonging to a node that is
not yet attached to the SPG. \item[2.] Stubs with label 2 are
attached to the shortest path graph (because the corresponding
node has been chosen), but not yet paired with another stub. These
are called `free stubs'. \item[3.] Stubs with label 3 in the SPG
are paired with another stub to form an edge in the SPG.
\end{enumerate}

The growth process as depicted in Figure \ref{fig-1} starts by
labelling all stubs by 1. Then, because we construct the SPG
starting from node $1$ we relabel the $D_1$ stubs of node $1$ with
the label $2$. We note that $Z_1^{\smallsup{1,N}}$ is equal to the
number of stubs connected to node 1, and thus
$Z_1^{\smallsup{1,N}}=D_1$. We next identify
$Z_j^{\smallsup{1,N}}$ for $j>1$. $Z_j^{\smallsup{1,N}}$ is
obtained by sequentially growing the SPG from the free stubs in
generation $Z_{j-1}^{\smallsup{1,N}}$. When all free stubs in
generation $j-1$ have chosen their connecting stub,
$Z_j^{\smallsup{1,N}}$ is equal to the number of stubs labeled 2
(i.e., free stubs) attached to the SPG. Note that not necessarily
each stub of $Z_{j-1}^{\smallsup{1,N}}$ contributes to stubs of
$Z_{j}^{\smallsup{1,N}}$, because a cycle may `swallow' two free
stubs. This is the case when a stub with label $2$ is chosen.

When a stub is chosen, we update the labels as follows:

\begin{enumerate}
\item[1.] If the chosen stub has label 1, in the SPG we connect
the present stub to the chosen stub to form an edge and attach the
remaining stubs of the chosen node as children. We update the
labels as follows. The present and chosen stub melt together to
form an edge and both are assigned label 3. All `brother' stubs
(except for the chosen stub) belonging to the same node of the
chosen stub receive label 2.

\item[2.] In this case we choose a stub with label 2, which is
already connected to the SPG.  For the SPG, a self-loop is created
if the chosen stub and present stub are `brother' stubs which
belong to the same node. If they are not `brother' stubs, then a
cycle is formed. Neither a self-loop nor a cycle changes the
distances to the root in the SPG.

The updating of the labels solely consists of changing the label
of the present and the chosen stub from 2 to 3.
\end{enumerate}

The above process stops in the $j^{\rm th}$ generation when there
are no more free stubs in generation $j-1$ for the SPG.

We continue the above process of drawing stubs until there are no
more stubs having label 1 or 2, so that all stubs have label 3.
Then, the SPG from node 1 is finalized, and we have generated the
shortest path graph as seen from node 1. We have thus obtained the
structure of the shortest path graph, and know how many nodes
there are at a given distance from node 1.

The above construction will be performed identically from node
$2$, and we denote the number of free stubs in the SPG of node 2
in generation $k$ by $Z_{k}^{\smallsup{2,N}}$. This construction
is close to being independent. In particular, it is possible to
couple the two SPG growth processes with two independent BP's.
This is described in detail in \cite[Section 3]{HHV03}. We make
essential use of the coupling between the SPG's and the BP's, in
particular, of \cite[Proposition A.3.1]{HHV03} in the appendix.
This completes the construction of the SPG's from both node 1 and
2.

\subsection{Bounds on the coupling}
\label{sec-couplingresults} We now investigate the growth of the
SPG, and its relationship to the BP with law $g$. In its
statement, we write, for $i=1,2$,
    \eq
    Y_n^{\smallsup{i,N}}=(\tau-2)^{n} \log (Z_{n}^{\smallsup{i,N}}\vee 1) \qquad \text{and} \qquad
    Y_n^{\smallsup{i}}=(\tau-2)^{n} \log ({\cal Z}^{\smallsup{i}}_{n}\vee 1),
    \label{Yidef2}
    \en
where $\{{\cal Z}^{\smallsup{1}}_{j}\}_{j\geq 1}$ and $\{{\cal
Z}^{\smallsup{2}}_{j}\}_{j\geq 1}$ are two independent delayed
BP's with offspring distribution $\{g\}$ and where ${\cal
Z}^{\smallsup{i}}_{1}$ has law $\{f\}$. Then the following
proposition shows that the first levels of the SPG are close to
those of the BP's:

\begin{prop}[Coupling at fixed time]
\label{prop-caft} For every $m$ fixed, and for $i=1,2$, there
exist {\rm independent} delayed BP's ${\cal Z}^{\smallsup{1}},
{\cal Z}^{\smallsup{2}}$, such that
    \eq
    \lim_{N\rightarrow \infty}
    \prob(Y_m^{\smallsup{i,N}}=Y_m^{\smallsup{i}})=
    1.
    \en
\end{prop}
In words, Proposition \ref{prop-caft} states that at any {\it
fixed} time, the SPG's from 1 and 2 can be coupled to two
independent BP's with offspring $g$, in such a way that the
probability that the SPG differs from the BP vanishes when
$N\rightarrow \infty$.

In the statement of the next proposition, we write, for $i=1,2$,
    \begin{eqnarray}
    \label{def-TmN}
    {\cal T}^{\smallsup{i,N}}_m={\cal
    T}^{\smallsup{i,N}}_m(\vep)&=&
    \{k> m:\big(Z_m^{\smallsup{i,N}}\big)^{\kappa^{k-m}}
    \leq N^{\frac{1-\vep^2}{\tau-1}}\}\nonumber\\
    &=&
    \{k> m:\kappa^{k}Y_m^{\smallsup{i,N}}
    \leq \frac{1-\vep^2}{\tau-1}\log N\},
    \end{eqnarray}
where we recall that $\kappa=(\tau-2)^{-1}$.

We will see that $Z_k^{\smallsup{i,N}}$ grows super-exponentially
with $k$ as long as $k\in {\cal T}^{\smallsup{i,N}}_m$. More
precisely, $\big(Z_m^{\smallsup{i,N}}\big)^{\kappa^{k-m}}$ is
close to $Z_k^{\smallsup{i,N}}$, and thus,  ${\cal
T}^{\smallsup{i,N}}_m$ can be thought of as the generations for
which the generation size is bounded by
$N^{\frac{1-\vep^2}{\tau-1}}$. The second main result of the
coupling is the following proposition:

\begin{prop}[Super-exponential growth with base $Y_m^{\smallsup{i,N}}$ for large times]
\label{prop-weakconv2bsec} If $F$ satisfies
Assumption~\ref{ass-gamma}, then for $i=1,2$,
    \eq
    \label{coupling}
    (a)\quad\qquad \prob(\vep\leq Y_m^{\smallsup{i,N}}\leq \vep^{-1}, \max_{k\in {\cal T}^{\smallsup{i,N}}_m(\vep)}
    |Y_k^{\smallsup{i,N}}-Y_m^{\smallsup{i,N}}|>\vep^3)=
    o_{\sN,m\vep}(1),
    \en
    \begin{eqnarray}
    \label{bdZT}
    (b)\quad&&\prob(\vep\leq Y_m^{\smallsup{i,N}}\leq \vep^{-1},\,\exists k\in{\cal T}^{\smallsup{i,N}}_{m}(\vep):\,Z_{k-1}^{\smallsup{i,N}}>Z_k^{\smallsup{i,N}})=o_{\sN,m\vep}(1),\\
    &&\prob(\vep\leq Y_m^{\smallsup{i,N}}\leq \vep^{-1},\,\exists k\in{\cal T}^{\smallsup{i,N}}_{m}(\vep):\,Z_k^{\smallsup{i,N}}>N^{\frac{1-\vepp}{\tau-1}})=o_{\sN,m\vep}(1),
    \label{monZT}
    \end{eqnarray}
where $o_{\sN,m\vep}(1)$ denotes a quantity $\gamma_{\sN,m,\vep}$
that converges to zero when first $N\rightarrow \infty$, then
$m\rightarrow \infty$ and finally $\vep\downarrow 0$.
\end{prop}

Proposition \ref{prop-weakconv2bsec} (a), i.e., (\ref{coupling}),
is the main coupling result used in this paper, and says that as
long as $k\in {\cal T}^{\smallsup{i,N}}_m(\vep)$, we have that
$Y_k^{\smallsup{i,N}}$ is close to $Y_m^{\smallsup{i,N}}$, which,
in turn, by Proposition \ref{prop-caft}, is close to
$Y_m^{\smallsup{i}}$. This establishes the coupling between the
SPG and the BP. Part (b) is a technical result used in the proof.
Equation (\ref{bdZT}) is a convenient result, as it shows that,
with high probability, $k\mapsto Z_k^{\smallsup{i,N}}$ is monotone
increasing. Equation (\ref{monZT}) shows that with high
probability $Z_k^{\smallsup{i,N}}\leq N^{\frac{1-\vepp}{\tau-1}}$
for all $k\in {\cal T}^{\smallsup{i,N}}_{m}(\vep)$, which allows
us to bound the number of free stubs in generation sizes that are
in ${\cal T}^{\smallsup{i,N}}_{m}(\vep)$.

We complete this section with a final coupling result, which shows
that for the first $k$ which is {\it not} in ${\cal
T}^{\smallsup{i,N}}_m(\vep)$, the SPG has many free stubs:

\begin{prop}[Lower bound on $Z_{k+1}^{\smallsup{i,N}}$ for
$k+1\not\in {\cal T}^{\smallsup{i,N}}_m(\vep)$] \label{prop-Tub}
Let $F$ satisfy Assumption~\ref{ass-gamma}.  Then,
    \eq
    \label{Zlb}
    \prob(k\in {\cal T}^{\smallsup{i,N}}_m(\vep), k+1\not\in {\cal
    T}^{\smallsup{i,N}}_m(\vep), \vep\leq
    Y_m^{\smallsup{i,N}}\leq \vep^{-1},
    Z_{k+1}^{\smallsup{i,N}}\leq N^{\frac{1-\vep}{\tau-1}})=o_{\sN,m,\vep}(1).
    \en
\end{prop}
 Propositions \ref{prop-caft}, \ref{prop-weakconv2bsec} and \ref{prop-Tub}
 will be proved in the appendix. We now prove the main results in
 Theorems \ref{thm-tau(2,3)} and \ref{thm-ll} subject to
 Propositions \ref{prop-caft}, \ref{prop-weakconv2bsec} and \ref{prop-Tub} in
 Section \ref{sec-pftau(2,3)}.

\section{Proof of Theorems \ref{thm-tau(2,3)} and \ref{thm-ll}}
\label{sec-pftau(2,3)}

In this section we prove Theorem \ref{thm-tau(2,3)} and identify
the limit in Theorem \ref{thm-ll}, using the coupling theory of
the previous section. For $i=1,2$,  we recall that
$Z_j^{\smallsup{i,N}}$ is the number of free stubs connected to
nodes on distance $j-1$ from root $i$. As we show in this section,
the hopcount $H_{\sN}$ is closely related to the SPG's
$\{Z_j^{\smallsup{i,N}}\}_{j\geq 0}, i=1,2$.

\subsection{Outline of the proof}
We start by describing the outline of the proof. The proof is
divided into several key steps proved in 5 subsections.

In the first key step of the proof, in Section \ref{sec-step1}, we
split the probability $\prob(H_{\sN}>k)$ into separate parts
depending on the values of $Y_m^{\smallsup{i,N}}=(\tau-2)^{m}\log
(Z_m^{\smallsup{i,N}}\vee 1)$. We prove that
    \eq
    \prob(H_{\sN}>k, Y_m^{\smallsup{1,N}}Y_m^{\smallsup{2,N}}=0)=1-q_m^2+o(1),
    \qquad N\rightarrow \infty,
    \en
where $1-q_m$ is the probability that the delayed BP $\{{\cal
Z}^{\smallsup{1}}_{j}\}_{j\geq 1}$ dies at or before the $m^{\rm
th}$ generation. When $m$ becomes large, then $q_m\rightarrow q$,
where $q$ equals the survival probability of the BP $\{{\cal
Z}_j^{\smallsup{1}}\}_{j\geq 1}$. This leaves us to determine the
contribution to $\prob(H_{\sN}>k)$ for the cases where
$Y_m^{\smallsup{1,N}}Y_m^{\smallsup{2,N}}>0$. We further show that
for $m$ large enough, and on the event that
$Y_m^{\smallsup{i,N}}>0$, \whps, $Y_m^{\smallsup{i,N}}\in [\vep,
\vep^{-1}],$ for $i=1,2$, where $\vep>0$ is small. This provides
us with a priori bounds on the shortest path graph exploration
processes $\{Z_j^{\smallsup{i,N}}\}$. We denote the event where
$Y_m^{\smallsup{i,N}}\in [\vep, \vep^{-1}],$ for $i=1,2$, by
$E_{m,\sN}(\vep)$.

The second key step in the proof, in Section \ref{sec-step2}, is
to obtain an asymptotic formula for $\prob(\{H_{\sN}>k\}\cap
E_{m,\sN}(\vep))$. Indeed, we prove the existence of
$\lambda=\lambda_{\sN}(k)>0$ such that
    \eq
    \label{firstlimit}
    \prob(\{H_{\sN}>k\}\cap E_{m,\sN}(\vep))= \expec\Big[{\bf 1}_{E_{m,\sN}(\vep)}
    \exp\Big\{-\lambda
    \frac{Z_{k_1+1}^{\smallsup{1,N}} Z_{k-k_1}^{\smallsup{2,N}}}{L_{\sN}}\Big\}\Big],
    \en
where the right-hand side is valid for {\it any} $k_1$ with $0\leq
2k_1\leq k-1$, and where $\lambda=\lambda_{\sN}(k)$ satisfies
$\frac 12 \leq \lambda_{\sN}(k)\leq 4 k$. It is even allowed that
$k_1$ is {\it random}, as long as it is measurable w.r.t.\
$\{Z_j^{\smallsup{i,N}}\}_{j=1}^m$. Even though the estimate on
$\lambda_{\sN}$ is not sharp, it turns out that it gives us enough
information to complete the proof. The bounds $\frac 12 \leq
\lambda_{\sN}(k)\leq 4 k$ play a crucial role in the remainder of
the proof.

In the third key step, in Section \ref{sec-step3}, we show that,
for $k=k_{\sN}\to \infty$, the main contribution of
(\ref{firstlimit}) stems from the term
    \begin{equation}
    \label{mini}
    \expec\Big[{\bf 1}_{E_{m,\sN}(\vep)}\exp\Big\{-\lambda \min_{k_1\in{\cal B}_{\sN}}
    \frac{Z_{k_1+1}^{\smallsup{1,N}}
    Z_{k_{\sN}-k_1}^{\smallsup{2,N}}}{L_{\sN}}
    \Big\}\Big],
    \end{equation}
with ${\cal B}_{\sN}={\cal B}_{\sN}(\vep,k_{\sN})$ defined in
(\ref{BNdef}) and is such that $k_1\in {\cal
B}_{\sN}(\vep,k_{\sN})$ precisely when $k_1+1 \in {\cal
T}_m^{\smallsup{1,N}}(\vep)$ and $k_{\sN}-k_1\in {\cal
T}_m^{\smallsup{2,N}}(\vep)$. Thus, by Proposition
\ref{prop-weakconv2bsec}, it implies that \whps
    $$
    Z_{k_1+1}^{\smallsup{1,N}}\leq N^{\frac{1-\vepp}{\tau-1}}\qquad
    \mbox{and} \qquad
    Z_{k_{\sN}-k_1}^{\smallsup{2,N}} \leq N^{\frac{1-\vepp}{\tau-1}}.
    $$
In turn, these bounds allow us to use Proposition
\ref{prop-weakconv2bsec}(a).

In the fourth key step, in Section \ref{sec-step4}, we proceed by
choosing
    \eq
    k_{\sN}=2 \left\lfloor \frac{\log \log N}{|\log
    (\tau-2)|}\right\rfloor+l,
    \label{kNldef}
    \en
and we show that with probability converging to 1 as
$\vep\downarrow 0$, the results of the coupling in Proposition
\ref{prop-weakconv2bsec} apply, which implies that
$Y_{k_1+1}^{\smallsup{1,N}}\approx Y^{\smallsup{1,N}}_m$ and
$Y_{k_{\sN}-k_1}^{\smallsup{2,N}}\approx Y^{\smallsup{2,N}}_m$.

In the final key step, in Section \ref{sec-step5}, the minimum
occurring in (\ref{mini}), with the approximations
$Y_{k_1+1}^{\smallsup{1,N}}\approx Y^{\smallsup{1,N}}_m$ and
$Y_{k_{\sN}-k_1}^{\smallsup{2,N}}\approx Y^{\smallsup{2,N}}_m$, is
analyzed. The main idea in this analysis is as follows. With the
above approximations, the expression in (\ref{mini}) can be
rewritten as
    \eq
    \label{approx3}
    \expec\Big[{\bf 1}_{E_{m,\sN}(\vep)}\exp\Big\{-\lambda\exp \Big[\min_{k_1\in \cB_{\sN}(\vep, k_{\sN})}
    (\kappa^{k_1+1}Y_{m}^{\smallsup{1,N}}
    +\kappa^{k_{\sN}-k_1}Y_{m}^{\smallsup{2,N}})-\log L_{\sN}\Big]\Big\}\Big]+o_{\sN,m,\vep}(1),
    \en
where $\kappa=(\tau-2)^{-1}>1$. The minimum appearing in the
exponent of (\ref{approx3}) is then rewritten (see (\ref{basis
asymp}) and (\ref{minoverZ})) as
    $$
    \kappa^{\lceil k_{\sN}/2\rceil} \big\{\min_{t \in \Z} (\kappa^{t}
    Y_{m}^{\smallsup{1,N}} +\kappa^{c_{l}-t} Y_{m}^{\smallsup{2,N}})
    -\kappa^{-\lceil k_{\sN}/2 \rceil}\log L_{\sN}\big\}
    .
    $$
Since $\kappa^{\lceil k_{\sN}/2\rceil}\to \infty$, the latter
expression only contributes to (\ref{approx3}) when
    $$
    \min_{t \in \Z} (\kappa^{t}
    Y_{m}^{\smallsup{1,N}} +\kappa^{c_{l}-t} Y_{m}^{\smallsup{2,N}})
    -\kappa^{-\lceil k_{\sN}/2 \rceil}\log L_{\sN}\leq 0.
    $$
Here it will become apparent that the bounds on $\lambda_{\sN}(k)$
are sufficient. The expectation of the indicator of this event
leads to the peculiar limit
    $$
    \prob\left(\min_{t\in \Z} (\kappa^{t} Y^{\smallsup{1}}
    +\kappa^{c_l-t} Y^{\smallsup{2}})\le \kappa^{a_{\sN}-\lceil
    l/2 \rceil}, Y^{\smallsup{1}}Y^{\smallsup{2}}>0\right),
    $$
with $a_{\sN}$ and $c_l$ as defined in Theorem \ref{thm-tau(2,3)}.
We complete the proof by showing that conditioning on the event
that 1 and 2 are connected is asymptotically equivalent to
conditioning on $Y^{\smallsup{1}}Y^{\smallsup{2}}>0$.

\begin{remark}
\label{rem-confmod} In the course of the proof, we will see that
it is not necessary that the degrees of the nodes are i.i.d. In
fact, in the proof below, we need that Propositions
\ref{prop-caft}--\ref{prop-Tub} are valid, as well as that
$L_{\sN}$ is concentrated around its mean $\mu N$. In Remark
\ref{rem-conf} in the appendix, we will investigate what is needed
in the proof of Propositions \ref{prop-caft}-- \ref{prop-Tub}. In
particular, the proof applies also to some instances of the
configuration model where the number of nodes with degree $k$ is
fixed, when we investigate the distance between two {\rm
uniformly} chosen nodes.
\end{remark}

We now go through the details of the proof. \vskip0.5cm

\subsection{A priory bounds on $Y_{m}^{\smallsup{i,N}}$}
\label{sec-step1} We wish to compute the probability
$\prob(H_{\sN}>k)$. To do so, we split $\prob(H_{\sN}>k)$ as
    \eq
    \label{splitprobs}
    \prob(H_{\sN}>k)=
    \prob(H_{\sN}>k, Y_m^{\smallsup{1,N}}Y_m^{\smallsup{2,N}}=0)
    +\prob(H_{\sN}>k, Y_m^{\smallsup{1,N}}Y_m^{\smallsup{2,N}}>0),
    \en
where we take $m$ to be sufficiently large. We will now prove two
lemmas, and use these to compute the first term in the right-hand
side of (\ref{splitprobs}).

\begin{lemma}
\label{lem-Y0}
     \[
     \lim_{N\rightarrow \infty}\prob(Y_m^{\smallsup{1,N}}Y_m^{\smallsup{2,N}}=0)
     =1-q_m^2,
     \]
where
    \[
    q_m=\prob(Y_m^{\smallsup{1}}>0).
    \]
\end{lemma}

\proof By Proposition \ref{prop-caft}, for $N\to\infty$, and
because $Y_m^{\smallsup{1}}$ and $Y_m^{\smallsup{2}}$ are
independent,
    \eqalign
    \prob(Y_m^{\smallsup{1,N}}Y_m^{\smallsup{2,N}}=0)
    &=\prob(Y_m^{\smallsup{1}}Y_m^{\smallsup{2}}=0)+o(1)=1-\prob(Y_m^{\smallsup{1}}Y_m^{\smallsup{2}}>0)+o(1)\\
    &=1-\prob(Y_m^{\smallsup{1}}>0)\prob(Y_m^{\smallsup{2}}>0)+o(1)
    =1-q_m^2 +o(1).\nn
    \enalign
\qed

The following lemma shows that the probability that $H_{\sN}\leq
m$ converges to zero for any {\it fixed} $m$:

\begin{lemma}
\label{lem-Hsmall} For any $m$ fixed,
    \eq
    \nn
    \lim_{N\rightarrow \infty} \prob(H_{\sN}\leq m)
    =0.
    \en
\end{lemma}

\proof As observed above Theorem \ref{thm-tau(2,3)}, by
exchangeability of the nodes $\{1,2, \ldots, N\}$,
    \eq
    \prob(H_{\sN}\leq m)=
    \prob(\widetilde H_{\sN}\leq m),
    \en
where $\widetilde H_{\sN}$ is the hopcount between node 1 and a
uniformly chosen node unequal to 1. We split, for any
$0<\delta<1$,
    \eqalign
    \label{splitdelta}
    &\prob(\widetilde H_{\sN}\leq m)
    =\prob(\widetilde H_{\sN}\leq m,
    \sum_{j\leq m} Z^{\smallsup{1,N}}_j \leq N^{\delta})
    +\prob(\widetilde H_{\sN}\leq m,
    \sum_{j\leq m} Z^{\smallsup{1,N}}_j > N^{\delta}).
    \enalign
The number of nodes at distance at most $m$ from node 1 is bounded
from above by $\sum_{j\leq m} Z^{\smallsup{1,N}}_j$. The event
$\{\widetilde H_{\sN}\leq m\}$ can only occur when the end node,
which is uniformly chosen in $\{2, \ldots, N\}$, is in the SPG of
node 1, so that
    \eq
    \prob(\widetilde H_{\sN}\leq m,
    \sum_{j\leq m} Z^{\smallsup{1,N}}_j \leq N^{\delta})
    \leq \frac{N^{\delta}}{N-1}=o(1),\qquad N\to \infty.
    \en
Therefore, the first term in (\ref{splitdelta}) is $o(1)$, as
required. We will proceed with the second term in
(\ref{splitdelta}). By Proposition \ref{prop-caft}, \whps, we have
that $Y_j^{\smallsup{1,N}}=Y_j^{\smallsup{1}}$ for all $j\leq m$.
Therefore, we obtain, because
$Y_j^{\smallsup{1,N}}=Y_j^{\smallsup{1}}$  implies
$Z_j^{\smallsup{1,N}}={\cal Z}_j^{\smallsup{1}}$,
    \eqalign
    \prob(\widetilde H_{\sN}\leq m,
    \sum_{j\leq m} Z^{\smallsup{1,N}}_j > N^{\delta})
    &\leq \prob(\sum_{j\leq m} Z^{\smallsup{1,N}}_j > N^{\delta})
    = \prob(
    \sum_{j\leq m} {\cal Z}^{\smallsup{1}}_j > N^{\delta})+o(1).\nn
    \enalign
However, when $m$ is fixed, the random variable $\sum_{j\leq m}
{\cal Z}^{\smallsup{1}}_j$ is finite with probability 1, and
therefore,
    \eq
    \lim_{N\rightarrow \infty} \prob(\widetilde H_{\sN}\leq m,
    \sum_{j\leq m} Z^{\smallsup{1,N}}_j > N^{\delta})=0.
    \en
This completes the proof of Lemma \ref{lem-Hsmall}. \qed
\vskip0.5cm

We now use Lemmas \ref{lem-Y0} and \ref{lem-Hsmall} to compute the
first term in (\ref{splitprobs}). We split
    \eq
    \label{Hsmall0}
    \prob(H_{\sN}> k, Y_m^{\smallsup{1,N}}Y_m^{\smallsup{2,N}}=0)
    =\prob(Y_m^{\smallsup{1,N}}Y_m^{\smallsup{2,N}}=0)-\prob(H_{\sN}\leq k,
    Y_m^{\smallsup{1,N}}Y_m^{\smallsup{2,N}}=0).
    \en
By Lemma \ref{lem-Y0}, the first term is equal to $1-q_m^2+o(1)$.
For the second term, we note that when $Y_m^{\smallsup{1,N}}=0$
and $H_{\sN}<\infty$, then $H_{\sN}\leq m-1,$ so that
    \eq
    \label{Hsmall1}
    \prob(H_{\sN}\leq k, Y_m^{\smallsup{1,N}}Y_m^{\smallsup{2,N}}=0)
    \leq \prob(H_{\sN}\leq m-1).
    \en
Using Lemma \ref{lem-Hsmall}, we conclude that
    \begin{corr}
    \label{cor-Y's0}For every $m$ fixed, and each $k\in \mathbb{N}$,
     \[
     \lim_{N\rightarrow \infty}
     \prob(H_{\sN}>k, Y_m^{\smallsup{1,N}}Y_m^{\smallsup{2,N}}=0)
     =1-q_m^2.
     \]
     \end{corr}

By Corollary \ref{cor-Y's0} and (\ref{splitprobs}), we are left to
compute $\prob(H_{\sN}>k,
Y_m^{\smallsup{1,N}}Y_m^{\smallsup{2,N}}>0).$ We first prove a
lemma that shows that if $Y_m^{\smallsup{1,N}}>0,$ then \whpl
$Y_{m}^{\smallsup{1,N}}\in [\vep, \vep^{-1}]$:
\begin{lemma}
\label{lem-mbd} For $i=1,2$,
     \begin{eqnarray}
     &&\limsup_{\vep \downarrow 0} \limsup_{m\rightarrow \infty}
\limsup_{N\rightarrow \infty}
     \prob(0<Y_{m}^{\smallsup{i,N}} < \vep)=
     \limsup_{\epsilon \downarrow 0} \limsup_{m\rightarrow \infty}
\limsup_{N\rightarrow \infty}
     \prob(Y_{m}^{\smallsup{i,N}}> \vep^{-1})=0.\nonumber
     \end{eqnarray}
\end{lemma}

\proof Fix $m$, when $N\rightarrow \infty$ it follows from
Proposition \ref{prop-caft} that
$Y_{m}^{\smallsup{i,N}}=Y_{m}^{\smallsup{i}}$, {\bf whp}. Thus, we
obtain that
     \begin{eqnarray}
     &&\limsup_{\vep \downarrow 0} \limsup_{m\rightarrow \infty}
\limsup_{N\rightarrow \infty}
     \prob(0<Y_{m}^{\smallsup{i,N}} < \vep)=
     \limsup_{\epsilon \downarrow 0} \limsup_{m\rightarrow \infty}
     \prob(0<Y^{\smallsup{i}}_m < \vep),\nonumber
     \end{eqnarray}
and similarly for the second probability. The remainder of the
proof of the lemma follows because $Y_m^{\smallsup{i}}
\stackrel{d}{\to} Y^{\smallsup{i}}$ as $m\rightarrow \infty$ and
is hence a tight sequence. \qed \vskip0.5cm

\noindent Write
    \eqalign
    E_{m,\sN}=E_{m,\sN}(\vep) &=\{Y_{m}^{\smallsup{i,N}}\in [\vep, \vep^{-1}],
    i=1,2\},\label{Edef}\\
    F_{m,\sN}=F_{m,\sN}(\vep)&=\big\{\max_{k\in {\cal T}^{\smallsup{N}}_m(\vep)}
    |Y_k^{\smallsup{i,N}}-Y_m^{\smallsup{i,N}}|\leq \vep^3, i=1,2\big\}.
    \label{Fdefpf}
    \enalign
As a consequence of Lemma \ref{lem-mbd}, we obtain that
     \eq
     \label{Ebd}
     \prob(E_{m,\sN}^c\cap \{Y_m^{\smallsup{1,N}}Y_m^{\smallsup{2,N}}>0\})=o_{\sN,m,\vep}(1).
     \en
In the sequel, we compute
    \eq
    \prob(\{H_{\sN}>k\}\cap E_{m,\sN}),
    \en
and often we make use of the fact that by Proposition
\ref{prop-weakconv2bsec}
     \eq
     \label{Fbdpf}
     \prob(E_{m,\sN}\cap F_{m,\sN}^c)=o_{\sN,m,\vep}(1).
     \en


\subsection{Asymptotics of $\prob(\{H_{\sN}>k\}\cap E_{m,\sN} )$}
\label{sec-step2} We next give a representation of
$\prob(\{H_{\sN}>k\}\cap E_{m,\sN})$. In order to do so, we write
$\Q_{\sZ}^{\smallsup{i,j}}$, where $i,j\geq 0$, for the
conditional probability given $\{Z_s^{\smallsup{1,N}}\}_{s=1}^i$
and $\{Z_s^{\smallsup{2,N}}\}_{s=1}^j$ (where, for $j=0$, we
condition only on $\{Z_s^{\smallsup{1,N}}\}_{s=1}^i$), and
$\expec_{\sZ}^{\smallsup{i,j}}$ for its conditional expectation.
Furthermore, we say that a random variable $k_1$ is {\it
$Z_m$-measurable} if $k_1$ is measurable with respect to the
$\sigma$-algebra generated by $\{Z_s^{\smallsup{1,N}}\}_{s=1}^{m}$
and $\{Z_s^{\smallsup{2,N}}\}_{s=1}^{m}$. The main rewrite is now
in the following lemma:

\begin{lemma}
\label{lem-rewprob} For $k\geq 2m-1$,
    \eq\prob(\{H_{\sN}>k\}\cap E_{m,\sN})=\expec\Big[{\bf 1}_{E_{m,\sN}}
    \Q_{\sZ}^{\smallsup{m,m}}
    (H_{\sN}>2m-1)\Zfactors\Big],
    \label{eqrewprob2}
    \en
where, for any $Z_m$-measurable $k_1$, with $m\leq k_1\leq
(k-1)/2$,
    \eqalign
    \label{eqrewprobQmm}
    \Zfactors &=
    \prod_{i=2m}^{2k_1}\Q_{\sZ}^{\smallsup{\lfloor i/2\rfloor+1,\lceil i/2\rceil}}
    (H_{\sN}>i|H_{\sN}>i-1)\\
    &\qquad \times\prod_{i=1}^{k-2k_1}
    \Q_{\sZ}^{\smallsup{k_1+1,k_1+i}}
    (H_{\sN}>2k_1+i|H_{\sN}>2k_1+i-1).\nn
    \enalign
\end{lemma}

\proof We start by conditioning on
$\{Z_s^{\smallsup{1,N}}\}_{s=1}^{m}$ and
$\{Z_s^{\smallsup{2,N}}\}_{s=1}^{m}$, and note that $E_{m,\sN}$ is
measurable w.r.t.\ $\{Z_s^{\smallsup{1,N}}\}_{s=1}^{m}$ and
$\{Z_s^{\smallsup{2,N}}\}_{s=1}^{m}$, so that we obtain, for
$k\geq 2m-1$,
    \eqalign
    \label{eqrewprob2b}
    \prob(\{H_{\sN}>k\}\cap E_{m,\sN})&=\expec\Big[{\bf 1}_{E_{m,\sN}}\Q_{\sZ}^{\smallsup{m,m}}
    (H_{\sN}>k)\Big]\\
    &=\expec\Big[{\bf 1}_{E_{m,\sN}}
    \Q_{\sZ}^{\smallsup{m,m}}
    (H_{\sN}>2m-1)\Q_{\sZ}^{\smallsup{m,m}}
    (H_{\sN}>k|H_{\sN}>2m-1)\Big].\nn
    \enalign
Moreover, for $i,j$ such that $i+j\leq k$,
    \eqalign
    &\Q_{\sZ}^{\smallsup{i,j}}(H_{\sN}>k|H_{\sN}>i+j-1)\\
    &=\expec_{\sZ}^{\smallsup{i,j}}\big[\Q_{\sZ}^{\smallsup{i,j+1}}(H_{\sN}>k|H_{\sN}>i+j-1)\big]\nn\\
    &=\expec_{\sZ}^{\smallsup{i,j}}\big[\Q_{\sZ}^{\smallsup{i,j+1}}(H_{\sN}>i+j|H_{\sN}>i+j-1)
    \Q_{\sZ}^{\smallsup{i,j+1}}(H_{\sN}>k|H_{\sN}>i+j)\big],\nn
    \enalign
and, similarly,
    \eqalign
    &\Q_{\sZ}^{\smallsup{i,j}}(H_{\sN}>k|H_{\sN}>i+j-1)\\
    &=\expec_{\sZ}^{\smallsup{i,j}}\big[\Q_{\sZ}^{\smallsup{i+1,j}}(H_{\sN}>i+j|H_{\sN}>i+j-1)
    \Q_{\sZ}^{\smallsup{i+1,j}}(H_{\sN}>k|H_{\sN}>i+j)\big].\nn
    \enalign
When we apply the above formulas, we can choose to increase $i$ or
$j$ by one depending on $\{Z_s^{\smallsup{1,N}}\}_{s=1}^i$ and
$\{Z_s^{\smallsup{2,N}}\}_{s=1}^j$.  We iterate the above
recursions until $i+j=k-1$. In particular, we obtain, for
$k>2m-1$,
    \eqalign
    \label{eqrewprob2c}
    \Q_{\sZ}^{\smallsup{m,m}}
    (H_{\sN}>k|H_{\sN}>2m-1)
    &=\expec_{\sZ}^{\smallsup{m,m}}
    \Big[\Q_{\sZ}^{\smallsup{m+1,m}}(H_{\sN}>2m|H_{\sN}>2m-1)\\
    &\qquad\qquad \times
    \Q_{\sZ}^{\smallsup{m+1,m}}(H_{\sN}>k|H_{\sN}>2m)\Big],\nn
    \enalign
so that, using that $E_{m,\sN}$ is
$\Q_{\sZ}^{\smallsup{m,m}}$-measurable and that
$\expec[\expec_{\sZ}^{\smallsup{m,m}}[X]]=\expec[X]$ for any
random  variable $X$,
    \eqalign
    \label{eqrewprob2d}
    &\prob(\{H_{\sN}>k\}\cap E_{m,\sN})\\
    &=\expec\Big[{\bf 1}_{E_{m,\sN}}
    \Q_{\sZ}^{\smallsup{m,m}}
    (H_{\sN}>2m-1)\Q_{\sZ}^{\smallsup{m+1,m}}(H_{\sN}>2m|H_{\sN}>2m-1)
    \Q_{\sZ}^{\smallsup{m+1,m}}(H_{\sN}>k|H_{\sN}>2m)\Big].\nn
    \enalign
We now compute the conditional probability by increasing $i$ or
$j$ as follows. For $i+j\leq 2k_1$, we will increase $i$ and $j$
in turn by 1, and for $i+j>2k_1$, we will only increase the second
component $j$. This leads to
    \eqalign
    \Q_{\sZ}^{\smallsup{m,m}}
    (H_{\sN}>k|H_{\sN}>2m-1)
    &=\expec_{\sZ}^{\smallsup{m,m}}\Big[\prod_{i=2m}^{2k_1}
    \Q_{\sZ}^{\smallsup{\lfloor i/2\rfloor+1,\lceil i/2\rceil}}
    (H_{\sN}>i|H_{\sN}>i-1)
    \label{eqrewprob2c}\\
    &\qquad \times
    \prod_{j=1}^{k-2k_1}
    \Q_{\sZ}^{\smallsup{k_1+1,k_1+j}}
    (H_{\sN}>2k_1+j|H_{\sN}>2k_1+j-1)\Big]\nn\\
    &=\expec_{\sZ}^{\smallsup{m,m}}[\Zfactors].\nn
    \enalign
Here, we use that we can move the expectations
$\expec_{\sZ}^{\smallsup{i,j}}$ outside, as in
(\ref{eqrewprob2d}), so that these do not appear in the final
formula. Therefore, from (\ref{eqrewprob2b}) and
(\ref{eqrewprob2c}),
    \eqalign
    \prob(\{H_{\sN}>k\}\cap E_{m,\sN})
    &=\expec\Big[{\bf 1}_{E_{m,\sN}}
    \Q_{\sZ}^{\smallsup{m,m}}
    (H_{\sN}>2m-1)\expec_{\sZ}^{\smallsup{m,m}}[\Zfactors] \Big]\nonumber\\
    &=\expec\Big[\expec_{\sZ}^{\smallsup{m,m}}[{\bf 1}_{E_{m,\sN}}\Q_{\sZ}^{\smallsup{m,m}}
    (H_{\sN}>2m-1)
    \Zfactors]\Big]\nn\\
    &=\expec\Big[{\bf 1}_{E_{m,\sN}}\Q_{\sZ}^{\smallsup{m,m}}
    (H_{\sN}>2m-1)
    \Zfactors\Big].
    \label{eqrewprob3}
    \enalign
This proves (\ref{eqrewprobQmm}). \qed
\medskip

We note that we can omit the term $\Q_{\sZ}^{\smallsup{m,m}}
(H_{\sN}>2m-1)$ in (\ref{eqrewprob2}) by introducing a small error
term. Indeed, we can write
    \eq
    \label{Qmmsplit}
    \Q_{\sZ}^{\smallsup{m,m}}(H_{\sN}>2m-1)
    =1-\Q_{\sZ}^{\smallsup{m,m}}(H_{\sN}\leq 2m-1).
    \en
For the contribution to (\ref{eqrewprob2}) due to the second term
in (\ref{Qmmsplit}), we bound ${\bf 1}_{E_{m,\sN}}\Zfactors\leq
1$. Therefore, the contribution to (\ref{eqrewprob2}) due to the
second term in (\ref{Qmmsplit}) is bounded by
    \eq
    \label{Qomiss}
    \expec\Big[\Q_{\sZ}^{\smallsup{m,m}}
    (H_{\sN}\leq 2m-1)\Big]=\prob(H_{\sN}\leq 2m-1)=o_{\sN}(1),
    \en
by Lemma \ref{lem-Hsmall}.

We conclude that by (\ref{Qomiss}), (\ref{Fbdpf}) and
(\ref{eqrewprob2}),
    \eqalign
    \label{rewprob3}
    \prob(\{H_{\sN}>k\}\cap E_{m,\sN})&=\expec\Big[{\bf 1}_{E_{m,\sN}\cap F_{m,\sN}}
    \Zfactors\Big]+o_{\sN,m,\vep}(1),
    \enalign
where we recall (\ref{eqrewprobQmm}) for the conditional
probability $\Zfactors$ appearing in (\ref{rewprob3}).

\vskip0.5cm

\noindent We continue with (\ref{rewprob3}) by investigating the
conditional probabilities in $\Zfactors$ defined in
(\ref{eqrewprobQmm}). We have the following bounds for
$\Q_{\sZ}^{\smallsup{i+1,j}}(H_{\sN}>i+j|H_{\sN}>i+j-1)$:
\begin{lemma}
\label{lem-match} For all integers $i,j\geq 0$,
    \[
    \exp\left\{-\frac{4Z_{i+1}^{\smallsup{1,N}}
    Z_{j}^{\smallsup{2,N}}}{L_{\sN}}\right\}\leq \Q_{\sZ}^{\smallsup{i+1,j}}(H_{\sN}>i+j|H_{\sN}>i+j-1)
    \leq\exp\left\{-\frac{Z_{i+1}^{\smallsup{1,N}}
    Z_{j}^{\smallsup{2,N}}}{2L_{\sN}}\right\}.
    \]
The upper bound is always valid, the lower bound is valid whenever
    \eq
    \label{assSPG}
    \sum_{s=1}^{i+1} Z_{s}^{\smallsup{1,N}}+\sum_{s=1}^j Z_{s}^{\smallsup{2,N}}
    \leq \frac{L_{\sN}}{4}.
    \en
\end{lemma}

\proof We start with the upper bound. We fix two sets of $n_1$ and
$n_2$ stubs, and will be interested in the probability that none
of the $n_1$ stubs are connected to the $n_2$ stubs. We order the
$n_1$ stubs in an arbitrary way, and connect the stubs iteratively
to other stubs. Note that we must connect at least $\lceil
n_1/2\rceil$ stubs, since any stub that is being connected removes
at most 2 stubs from the total of $n_1$ stubs. The number $n_1/2$
is reached for $n_1$ even precisely when all the $n_1$ stubs are
connected with each other. Therefore, we obtain that the
probability that the $n_1$ stubs are not connected to the $n_2$
stubs is bounded from above by
    \eq
    \prod_{i=1}^{\lceil n_1/2\rceil} \Big(1-\frac{n_2}{L_{\sN}-2i+1}\Big).
    \en
To complete the upper bound, we note that
    \eq
    1-\frac{n_2}{L_{\sN}-2i+1}\leq 1-\frac{n_2}{L_{\sN}}\leq e^{-\frac{n_2}{L_{\sN}}},
    \en
to obtain that the probability that the $n_1$ stubs are not
connected to the $n_2$ stubs is bounded from above by
    \eq
    e^{-\lceil n_1/2\rceil \frac{n_2}{L_{\sN}}}\leq e^{-\frac{n_1n_2}{2L_{\sN}}}.
    \en
Applying the above bound to $n_1=Z_{i+1}^{\smallsup{1,N}}$ and
$n_2=Z_{j}^{\smallsup{2,N}}$, and noting that the probability that
$H_{\sN}>i+j$ given that $H_{\sN}>i+j-1$ is bounded from above by
the probability that the stubs in $Z_{i+1}^{\smallsup{1,N}}$ are
not connected to the stubs in $Z_{j}^{\smallsup{2,N}}$ leads to
    \eq
    \Q_{\sZ}^{\smallsup{i+1,j}}(H_{\sN}>i+j|H_{\sN}>i+j-1)
    \leq \exp\left\{-\frac{Z_{i+1}^{\smallsup{1,N}}
    Z_{j}^{\smallsup{2,N}}}{2L_{\sN}}\right\},
    \en
which completes the proof of the upper bound.

We again fix two sets of $n_1$ and $n_2$ stubs, and are again
interested in the probability that none of the $n_1$ stubs are
connected to the $n_2$ stubs. However, now we use these bounds
repeatedly, and we assume that in each step there remain to be at
least $M$ stubs available. We order the $n_1$ stubs in an
arbitrary way, and connect the stubs iteratively to other stubs.
We obtain a lower bound by further requiring that the $n_1$ stubs
do not connect to each other. Therefore, the probability that the
$n_1$ stubs are not connected to the $n_2$ stubs is bounded below
by
    \eq
    \prod_{i=1}^{n_1} \Big(1-\frac{n_2}{M-2i+1}\Big).
    \en
When $M-2n_1\geq \frac{L_{\sN}}{2}$, we obtain that
$1-\frac{n_2}{M-2i+1} \geq 1-\frac{2n_2}{L_{\sN}}$. Moreover, when
$x\leq \frac 12$, we have that $1-x \geq e^{-2x}$. Therefore, we
obtain that when $M-2n_1\geq \frac{L_{\sN}}{2}$ and $n_2\leq
\frac{L_{\sN}}{4}$, then the probability that the $n_1$ stubs are
not connected to the $n_2$ stubs when there are still at least $M$
stubs available is bounded below by
    \eq
    \label{boundnmM}
    \prod_{i=1}^{n_1} \Big(1-\frac{n_2}{M-2i+1}\Big)\geq \prod_{i=1}^{n_1} e^{-\frac{4n_2}{L_{\sN}}}
    =e^{-\frac{4n_1n_2}{L_{\sN}}}.
    \en

The event $H_{\sN}>i+j$ conditionally on $H_{\sN}>i+j-1$ precisely
occurs when the stubs $Z_{i+1}^{\smallsup{1,N}}$ are not connected
to the stubs in $Z_{j}^{\smallsup{2,N}}$. We will assume that
(\ref{assSPG}) holds. We have that $M=L_{\sN}-2\sum_{s=1}^{i+1}
Z_{s}^{\smallsup{1,N}}-2\sum_{s=1}^j Z_{s}^{\smallsup{2,N}},$ and
$n_1=Z_{i+1}^{\smallsup{1,N}}$, $n_2=Z_{j}^{\smallsup{2,N}}$.
Thus, $M-2n_1\geq \frac{L_{\sN}}{2}$ happens precisely when
    \eq
    \label{bdMconseq}
    M-2n_1\geq L_{\sN}-2\sum_{s=1}^{i+1} Z_{s}^{\smallsup{1,N}}-2\sum_{s=1}^j Z_{s}^{\smallsup{2,N}}
    \geq \frac{L_{\sN}}{2}.
    \en
This follows from the assumed bound in (\ref{assSPG}). Also, when
$n_2=Z_{j}^{\smallsup{2,N}}$, $n_2\leq \frac{L_{\sN}}{4}$ is
implied by (\ref{assSPG}). Thus, we are allowed to use the bound
in (\ref{boundnmM}). This leads to
    \eqalign
    \Q_{\sZ}^{\smallsup{i+1,j}}(H_{\sN}>i+j|H_{\sN}>i+j-1) &\geq
    \exp\Big\{-\frac{4Z_{i+1}^{\smallsup{1,N}}Z_{j}^{\smallsup{2,N}}}{L_{\sN}}\Big\},
    \enalign
which completes the proof of Lemma \ref{lem-match}. \qed


\subsection{The main contribution to $\prob(\{H_{\sN}>k\}\cap E_{m,\sN})$}
\label{sec-step3}

We rewrite the expression in (\ref{rewprob3}) in a more convenient
form, using Lemma \ref{lem-match}. We derive an upper and a lower
bound. For the upper bound, we bound all terms appearing on the
right hand side of (\ref{eqrewprobQmm}) by 1, except for the term
$\Q_{\sZ}^{\smallsup{k_1+1,k-k_1}} (H_{\sN}>k|H_{\sN}>k-1),$ which
arises when $i=k-2k_1$. Using the upper bound in Lemma
\ref{lem-match}, we thus obtain that
    \eqalign
    \label{rewprobub}
    \Zfactors
    &\leq
    \exp\big\{-\frac{Z_{k_1+1}^{\smallsup{1,N}}
    Z_{k-k_1}^{\smallsup{2,N}}}{2L_{\sN}}\big\}.
    \enalign
The latter inequality is true for any $Z_m$-measurable $k_1$ with
$m\leq k_1\leq (k-1)/2$.

To derive the lower bound, we next assume that
    \eq
    \label{assSPG2}
    \sum_{s=1}^{k_1+1} Z_{s}^{\smallsup{1,N}}+
    \sum_{s=1}^{k-k_1} Z_{s}^{\smallsup{2,N}}
    \leq \frac{L_{\sN}}{4},
    \en
so that (\ref{assSPG}) is satisfied for all $i$ in
(\ref{eqrewprobQmm}). We write, recalling (\ref{def-TmN}),
    \eq
    \label{BN1def}
    {\cal B}_{\sN}^{\smallsup{1}}(\vep,k)=\Big\{m\leq l\leq (k-1)/2:~
    l+1\in {\cal T}^{\smallsup{1,N}}_m(\vep),~k-l\in {\cal T}^{\smallsup{2,N}}_m(\vep)
    \Big\}.
    \en

We restrict ourselves to $k_1\in {\cal
B}_{\sN}^{\smallsup{1}}(\vep,k),$ if ${\cal
B}_{\sN}^{\smallsup{1}}(\vep,k)\neq \varnothing$. When $k_1\in
{\cal B}_{\sN}^{\smallsup{1}}(\vep,k)$, we are allowed to use the
bounds in Proposition \ref{prop-weakconv2bsec}. Note that
$\{k_1\in {\cal B}_{\sN}^{\smallsup{1}}(\vep,k)\}$ is
$Z_m-$measurable. Moreover, it follows from Proposition
\ref{prop-weakconv2bsec} that if $k_1\in {\cal
B}_{\sN}^{\smallsup{1}}(\vep,k),$ that then, with probability
converging to 1 as first $N\rightarrow \infty$ and then
$m\rightarrow \infty,$
    \eq
    \label{implBN1}
    Z_{s}^{\smallsup{1,N}}\leq N^{\frac{1-\vepp}{\tau-1}}, \quad \forall m< s\leq k_1+1,
    \qquad \text{while} \qquad
    Z_{s}^{\smallsup{2,N}}\leq N^{\frac{1-\vepp}{\tau-1}}, \quad \forall m< s\leq k-k_1.
    \en
Therefore, when $k_1\in {\cal B}_{\sN}^{\smallsup{1}}(\vep,k),$
the assumption in (\ref{assSPG2}) is satisfied with probability
$1-o_{\sN,m}(1)$, as long as $k=O(N^{\frac{\tau-2}{\tau-1}})$. The
latter restriction is not serious, as we always have $k$ in mind
for which $k=O(\log\log{N})$ (see e.g.\ Theorem
\ref{thm-tau(2,3)}).

Thus, on the event $E_{m,\sN}\cap  \{k_1\in {\cal
B}^{\smallsup{1}}_{\sN}(\vep,k)\}$, using (\ref{bdZT}) in
Proposition \ref{prop-weakconv2bsec} and the lower bound in Lemma
\ref{lem-match}, with probability $1-o_{\sN,m,\vep}(1)$, and for
all $i\in \{2m, \ldots, 2k-1\}$,
    \eq
    \Q_{\sZ}^{\smallsup{\lfloor i/2\rfloor+1,\lceil i/2\rceil}}
    (H_{\sN}>i|H_{\sN}>i-1)
    \geq \exp\big\{-\frac{4Z_{\lfloor i/2\rfloor+1}^{\smallsup{1,N}}
    Z_{\lceil i/2\rceil}^{\smallsup{2,N}}}{L_{\sN}}\big\}
    \geq \exp\big\{-\frac{4Z_{k_1+1}^{\smallsup{1,N}}
    Z_{k-k_1}^{\smallsup{2,N}}}{L_{\sN}}\big\},
    \en
and, for $1\leq i\leq k-2k_1$,
    \eq
    \Q_{\sZ}^{\smallsup{k_1+1,k_1+i}}
    (H_{\sN}>2k_1+i|H_{\sN}>2k_1+i-1)
    \geq
    \exp\big\{-\frac{4Z_{k_1+1}^{\smallsup{1,N}}
    Z_{k_1+i}^{\smallsup{2,N}}}{L_{\sN}}\big\}\geq
    \exp\big\{-\frac{4Z_{k_1+1}^{\smallsup{1,N}}
    Z_{k-k_1}^{\smallsup{2,N}}}{L_{\sN}}\big\}.
    \en

Therefore, by Lemma \ref{lem-rewprob}, and using the above bounds
for each of the in total $k$ terms, we obtain that when $k_1\in
{\cal B}_{\sN}^{\smallsup{1}}(\vep,k) \neq \varnothing$, and with
probability $1-o_{\sN,m,\vep}(1)$,
    \eqalign
    \label{rewproblb}
    &\Zfactors\geq
    \exp\big\{-4k\frac{Z_{k_1+1}^{\smallsup{1,N}}
    Z_{k-k_1}^{\smallsup{2,N}}}{L_{\sN}}\big\}.
    \enalign

We next use the symmetry for the nodes 1 and 2. Denote
    \eqalign
    \label{BN2def}
    {\cal B}_{\sN}^{\smallsup{2}}(\vep,k)&=\Big\{m\leq l\leq (k-1)/2:
    ~l+1\in {\cal T}^{\smallsup{2,N}}_m(\vep),~k-l\in {\cal T}^{\smallsup{1,N}}_m(\vep)
    \Big\}.
    \enalign
Take $\tilde l=k-l-1$, so that $(k-1)/2\leq \tilde l \leq k-1-m,$
and thus
    \eqalign
    \label{BN2eq}
    {\cal B}_{\sN}^{\smallsup{2}}(\vep,k)=\Big\{(k-1)/2\leq \tilde l \leq k-1-m:
    ~\tilde l+1\in {\cal T}^{\smallsup{1,N}}_m(\vep),~k-\tilde l\in {\cal T}^{\smallsup{2,N}}_m(\vep)
    \Big\}.
    \enalign
Then, since the nodes 1 and 2 are exchangeable, we obtain from
(\ref{rewproblb}), when $k_1\in {\cal
B}_{\sN}^{\smallsup{2}}(\vep,k) \neq \varnothing$, and with
probability $1-o_{\sN,m,\vep}(1)$,
    \eqalign
    \label{rewproblb2}
    &\Zfactors\geq
    \exp\big\{-4k\frac{Z_{k_1+1}^{\smallsup{1,N}}
    Z_{k-k_1}^{\smallsup{2,N}}}{L_{\sN}}\big\}.
    \enalign
We define ${\cal B}_{\sN}(\vep,k)={\cal
B}_{\sN}^{\smallsup{1}}(\vep,k) \cup {\cal
B}_{\sN}^{\smallsup{2}}(\vep,k)$, which is equal to
    \eq
    \label{BNdef}
    {\cal B}_{\sN}(\vep,k)
    =\Big\{m\leq l\leq k-1-m:~l+1\in {\cal T}^{\smallsup{1,N}}_m(\vep),~k-l\in {\cal T}^{\smallsup{2,N}}_m(\vep)
    \Big\}.
    \en
We can summarize the obtained results by writing that with
probability $1-o_{\sN,m,\vep}(1)$, and when ${\cal
B}_{\sN}(\vep,k) \neq \varnothing$, we have
    \eqalign
    \label{rel451}
    \Zfactors&=
    \exp\big\{-\lambda_{\sN}\frac{Z_{k_1+1}^{\smallsup{1,N}}
    Z_{k-k_1}^{\smallsup{2,N}}}{L_{\sN}}\big\},
    \enalign
for all $k_1\in {\cal B}_{\sN}(\vep,k)$, where
$\lambda_{\sN}=\lambda_{\sN}(k)$ satisfies
    \eq
    \label{lambdabd}
    \frac 12 \leq \lambda_{\sN}(k)\leq 4k.
    \en
Relation (\ref{rel451}) is true for any $k_1\in {\cal
B}_{\sN}(\vep,k)$. However, our coupling fails when
$Z_{k_1+1}^{\smallsup{1,N}}$ or $Z_{k-k_1}^{\smallsup{2,N}}$ grows
too large, since we can only couple $Z_{j}^{\smallsup{i,N}}$ with
$\hat Z_{j}^{\smallsup{i,N}}$ up to the point where
$Z_{j}^{\smallsup{i,N}}\leq N^{\frac{1-\vep^2}{\tau-1}}$.
Therefore, we next take the maximal value over $k_1\in {\cal
B}_{\sN}(\vep,k)$ to arrive at the fact that, with probability
$1-o_{\sN,m,\vep}(1)$, on the event that ${\cal
B}_{\sN}(\vep,k)\neq \varnothing$,
    \eqalign
    \label{mainrewb}
    \Zfactors&=
    \max_{k_1\in {\cal B}_{\sN}(\vep,k)}
    \exp\big\{-\lambda_{\sN}\frac{Z_{k_1+1}^{\smallsup{1,N}}
    Z_{k-k_1}^{\smallsup{2,N}}}{L_{\sN}}\big\}
    =
    \exp\Big\{-\lambda_{\sN}\min_{k_1\in {\cal B}_{\sN}(\vep,k)}
    \frac{Z_{k_1+1}^{\smallsup{1,N}}
    Z_{k-k_1}^{\smallsup{2,N}}}{L_{\sN}}\Big\}.
    \enalign
We conclude that
    \eqalign
    \label{mainrew}
    &\prob(\{H_{\sN}>k\} \cap E_{m,\sN}\cap \{{\cal B}_{\sN}(\vep,k)\neq \varnothing\})\\
    &\qquad =\expec\Big[{\bf 1}_{E_{m,\sN}}
    \exp\Big\{-\lambda_{\sN}\min_{k_1\in {\cal B}_{\sN}(\vep,k)}
    \frac{Z_{k_1+1}^{\smallsup{1,N}}
    Z_{k-k_1}^{\smallsup{2,N}}}{L_{\sN}}\Big\}\Big]
    +o_{\sN,m,\vep}(1).\nn
    \enalign
From here on we take $k=k_{\sN}$ as in (\ref{kNldef}) with $l$ a
fixed integer.

In Section \ref{sec-pflemmas}, we prove the following lemma that
shows that, apart from an event of probability
$1-o_{\sN,m,\vep}(1)$, we may assume that ${\cal
B}_{\sN}(\vep,k_{\sN})\neq \varnothing$:
    \begin{lemma} For all $l$, with $k_{\sN}$ as in (\ref{kNldef}),
    \label{lem-BNempty}
    \eq
    \limsup_{\vep\downarrow 0} \limsup_{m\rightarrow \infty}
    \limsup_{N\rightarrow\infty}
    \prob(\{H_{\sN}>k_{\sN}\} \cap E_{m,\sN}\cap \{{\cal B}_{\sN}(\vep,k_{\sN})=\varnothing\})
    =0.\nn
    \en
    \end{lemma}

From now on, we will abbreviate ${\cal B}_{\sN}={\cal
B}_{\sN}(\vep,k_{\sN})$. Using (\ref{mainrew}) and Lemma
\ref{lem-BNempty}, we conclude

    \begin{corr}
    \label{cor-mainrew2}
    For all $l$, with $k_{\sN}$ as in (\ref{kNldef}),
    \eqalign
    &\prob\big(\{H_{\sN}>k_{\sN}\} \cap E_{m,\sN}\big)
    =\expec\Big[{\bf 1}_{E_{m,\sN}}
    \exp\Big\{-\lambda_{\sN}\min_{k_1\in {\cal B}_{\sN}}
    \frac{Z_{k_1+1}^{\smallsup{1,N}}
    Z_{k-k_1}^{\smallsup{2,N}}}{L_{\sN}}\Big\}\Big]+o_{\sN,m,\vep}(1).\nn
    \enalign
    \end{corr}

\subsection{Application of the coupling results}
\label{sec-step4} In this section, we use the coupling results in
Section \ref{sec-couplingresults}. Before doing so, we investigate
the minimum of the function $t\mapsto
\kappa^{t}y_1+\kappa^{n-t}y_2 $, where the minimum is taken over
the discrete set $\{0,1,\ldots,n\}$, and we recall that
$\kappa=(\tau-2)^{-1}$.

\begin{lemma}
\label{Lemma_mintZ} Suppose that $y_1>y_2>0$, and
$\kappa=(\tau-2)^{-1}>1$. For the integer $n>
\frac{-\log(y_2/y_1)}{\log \kappa}$,
\begin{equation*}
 t^*=\argmin_{t\in \{0,1,\ldots,n\}}
 \left(
\kappa^{t}y_1+\kappa^{n-t}y_2 \right)= \rm{round}\left(
\frac{n}{2}+\frac{\log(y_2/y_1)}{2\log \kappa}\right),
\end{equation*}
where {\it round}$(x)$ is $x$ rounded off to the nearest integer.
In particular,
    \begin{eqnarray*}
    \max
    \left\{
    \frac{\kappa^{t^*}y_1}{\kappa^{n-t^*}y_2},
    \frac{\kappa^{n-t^*}y_2}{\kappa^{t^*}y_1}
    \right\}
    \leq
    \kappa.
    \end{eqnarray*}
\end{lemma}

\proof Consider, for real valued $t\in[0,n]$, the function
    $$
    \psi(t)=\kappa^{t}y_1+\kappa^{n-t}y_2.
    $$
Then,
    \[
    \psi'(t)=(\kappa^{t}y_1-\kappa^{n-t}y_2)\log \kappa,\qquad
    \psi''(t)=(\kappa^{t}y_1+\kappa^{n-t}y_2)\log^2 \kappa.
    \]
In particular, $\psi''(t)>0$, so that the function $\psi$ is
strictly convex. The unique minimum of $\psi$ is attained at
${\hat t}$, satisfying $\psi'({\hat t})=0$, i.e.,
    $$
    {\hat t}=\frac{n}{2}+\frac{\log (y_2/y_1)}{2\log \kappa} \in (0,n),
    $$
because $n>-\log (y_2/y_1)/\log \kappa$. By convexity $t^*=\lfloor
{\hat t}\rfloor$ or $t^*=\lceil {\hat t}\rceil$. We will show that
$|t^*-{\hat t}|\leq \frac12$. Put $t_1^*=\lfloor {\hat t}\rfloor$
and $t_2^*=\lceil {\hat t}\rceil$. We have
    \eq
    \label{hattform}
    \kappa^{{\hat t}}y_1=\kappa^{n-{\hat t}}y_2=\kappa^{\frac{n}{2}}\sqrt{y_1y_2}.
    \en
Writing $t_i^*={\hat t}+t_i^*-{\hat t}$, we obtain for $i=1,2$,
    $$
    \psi(t_i^*)=\kappa^{\frac{n}{2}}\sqrt{y_1y_2}
    \{\kappa^{t_i^*-{\hat t}}+\kappa^{{\hat t}-t_i^*}
    \}.
    $$
For $0<x<1$, the function $x\mapsto \kappa^x+\kappa^{-x}$ is
increasing so $\psi(t_1^*)\le \psi(t_2^*)$ if and only if ${\hat
t}-t_1^*\le t_2^*-{\hat t} $, or ${\hat t}-t_1^*\le \frac12$,
i.e., if $\psi(t_1^*)\le \psi(t_2^*)$ and hence the minimum over
the discrete set $\{0,1,\ldots,n\}$ is attained at $t_1^*$, then
${\hat t}-t_1^*\le \frac12$. On the other hand, if $\psi(t_2^*)\le
\psi(t_1^*)$, then by the `only if' statement we find $t_2^*-{\hat
t}\le \frac12$. In both cases we have $|t^*-{\hat t}|\le \frac12$.
Finally, if $t^*=t^*_1$, then we obtain, using (\ref{hattform}),
    $$
    1\leq \frac{\kappa^{n-t^*}y_2}{\kappa^{t^*}y_1}
    =\frac{\kappa^{{\hat t}-t_1^*}}{\kappa^{t_1^*-{\hat t}}}
    =
    \kappa^{2({\hat t}-t_1^*)}\leq \kappa,
    $$
while for $t^*=t^*_2$, we obtain $1\leq
\frac{\kappa^{t^*}y_1}{\kappa^{n-t^*}y_2} \leq \kappa$.
 \qed
\vskip0.5cm

\noindent We continue with our investigation of
$\prob\big(\{H_{\sN}>k_{\sN}\} \cap E_{m,\sN}\big)$. We start from
Corollary \ref{cor-mainrew2}, substituting (\ref{Yidef2}),
    \eqalign
    &\prob\big(\{H_{\sN}>k_{\sN}\} \cap E_{m,\sN}\big)\\
    &=
    \expec\Big[{\bf 1}_{E_{m,\sN}}
    \exp\Big\{-\lambda_{\sN}\exp \Big[\min_{k_1\in {\cal B}_{\sN}}
    \big(\kappa^{k_1+1}Y_{k_1+1}^{\smallsup{1,N}}
    +\kappa^{k_{\sN}-k_1}Y_{k_{\sN}-k_1}^{\smallsup{2,N}}\big)
    -\log L_{\sN}\Big]
    \Big\}
    \Big]+o_{\sN,m,\eps}(1),\nn
    \enalign
where we rewrite, using (\ref{BNdef}) and (\ref{def-TmN}),
    \eq
    \label{BN-rew}
    {\cal B}_{\sN}=
    \{m\le k_1\leq k_{\sN}-1-m: \kappa^{k_1+1}Y_{m}^{\smallsup{1,N}}\leq  \frac{1-\vep^2}{\tau-1}
    \log N,
    \kappa^{k_{\sN}-k_1}Y_{m}^{\smallsup{2,N}}\leq \frac{1-\vep^2}{\tau-1}\log N\}.
    \en
Moreover, according to Proposition \ref{prop-weakconv2bsec}, and
with probability at least $1-o_{\sN,m,\vep}(1)$, we have that
$\min_{k_1\in {\cal
B}_{\sN}}(\kappa^{k_1+1}Y_{k_1+1}^{\smallsup{1,N}}
+\kappa^{k_{\sN}-k_1}Y_{k_{\sN}-k_1}^{\smallsup{2,N}})$ is between
    $$
    \min_{k_1\in {\cal B}_{\sN}}\big(\kappa^{k_1+1}(Y_{m}^{\smallsup{1,N}}-\vep^3)
    +\kappa^{k_{\sN}-k_1}(Y_{m}^{\smallsup{2,N}}-\vep^3)\big)
    $$
and
    $$
    \min_{k_1\in {\cal B}_{\sN}}\big(\kappa^{k_1+1}(Y_{m}^{\smallsup{1,N}}+\vep^3)
    +\kappa^{k_{\sN}-k_1}(Y_{m}^{\smallsup{2,N}}+\vep^3)\big).
    $$
To abbreviate the notation, we will write, for $i=1,2$,
    \eq
    \Ymplus{i}=Y_{m}^{\smallsup{i,N}}+\vep^3, \qquad \Ymmin{i}=Y_{m}^{\smallsup{i,N}}-\vep^3.
    \en

Define for $\vep > 0$,
    $$
    H_{m,\sN}=H_{m,\sN}(\vep)
    =\left\{
    \min_{0\le k_1 \le k_{\sN}-1}
    \big(\kappa^{k_1+1}\Ymmin{1}
    +\kappa^{k_{\sN}-k_1}\Ymmin{2}\big)\leq (1+\vep^2) \log N
    \right\}.
    $$
On the complement $H^c_{m,\sN}$, the minimum over $0\le k_1 \le
k_{\sN}-1$ of
$\kappa^{k_1+1}\Ymmin{1}+\kappa^{k_{\sN}-k_1}\Ymmin{2}$ exceeds
$(1+\vep^2)\log N$. Therefore, also the minimum over the set
${\cal B}_{\sN}$ of
$\kappa^{k_1+1}\Ymmin{1}+\kappa^{k_{\sN}-k_1}\Ymmin{2}$ exceeds
$(1+\vep^2)\log N$, so that, using Lemma \ref{lem-match}, and with
error at most $o_{\sN,m,\eps}(1)$,
    \eqalign
    &\prob\Big(\{H_{\sN}>k_{\sN}\} \cap E_{m,\sN}\cap H^c_{m,\sN}\Big)\nn\\
    &\qquad \leq\expec\Big[{\bf 1}_{H^c_{m,\sN}}\exp\Big\{- \frac{1}{2}\exp
    \Big[\min_{k_1\in {\cal B}_{\sN}}
    \big(\kappa^{k_1+1}Y_{k_1+1}^{\smallsup{1,N}}
    +\kappa^{k_{\sN}-k_1}Y_{k_{\sN}-k_1}^{\smallsup{2,N}}\big)-\log L_{\sN}\Big]
    \Big\}\Big]\nn\\
    &\qquad \leq
    \expec\Big[{\bf 1}_{H^c_{m,\sN}}\exp\Big\{-\frac{1}{2}\exp
    \Big[\min_{ k_1\in {\cal B}_{\sN}}
    \big(\kappa^{k_1+1}\Ymmin{1}
    +\kappa^{k_{\sN}-k_1}\Ymmin{2}\big)-\log L_{\sN}
    \Big]\Big\}\Big]\nn\\
    &\qquad \le \expec\Big[\exp\Big\{- \frac{1}{2}\exp \big((1+\vep^2)\log N-\log{L_{\sN}}
    \big)\Big\}\Big]
    \leq e^{-\frac{1}{2c}N^{\vep^2}},
    \enalign
because $L_{\sN}\leq cN$, {\bf whp}, as $N\to\infty$. Therefore,
in the remainder of the proof, we assume that $H_{m,\sN}$ holds.

We next show that \whpl,
    \eq
    \min_{k_1\in {\cal B}_{\sN}}
     \big(\kappa^{k_1+1}\Ymplus{1}+\kappa^{k_{\sN}-k_1}\Ymplus{2}\big)
    =
    \min_{0\leq k_1< k_{\sN}}
     \big(\kappa^{k_1+1}\Ymplus{1}+\kappa^{k_{\sN}-k_1}\Ymplus{2}\big),
     \label{fullmin1}
    \en
and
    \eq
    \min_{k_1\in {\cal B}_{\sN}}
    \big(\kappa^{k_1+1}\Ymmin{1}+\kappa^{k_{\sN}-k_1}\Ymmin{2}\big)
    =
    \min_{0\leq k_1< k_{\sN}}
    \big(\kappa^{k_1+1}\Ymmin{1}+\kappa^{k_{\sN}-k_1}\Ymmin{2}\big).
    \label{fullmin2}
    \en
We start with (\ref{fullmin1}), the proof of (\ref{fullmin2}) is
similar, and, in fact, slightly simpler, and is therefore omitted.
To prove (\ref{fullmin1}), we use Lemma \ref{Lemma_mintZ}, with
$n=k_{\sN}+1$, $t=k_1+1$, $y_1=\Ymplus{1}$ and $y_2=\Ymplus{2}$.
Let
    $$
    t^*=\argmin_{t\in \{0,1,\ldots,n\}}
    \left(
    \kappa^{t}y_1+\kappa^{n-t}y_2
    \right),
    $$
and assume (without restriction) that $\kappa^{t^*}y_1\geq
\kappa^{n-t^*}y_2$. We have to show that $t^*-1\in {\cal
B}_{\sN}$. According to Lemma \ref{Lemma_mintZ},
    \eq
    \label{x/ybd}
    1\leq \frac{\kappa^{t^*} \Ymplus{1}}{\kappa^{n-t^*}\Ymplus{2}}=
    \frac{\kappa^{t^*}y_1}{\kappa^{n-t^*}y_2}
    \leq  \kappa.
    \en
We define $x=\kappa^{t^*} \Ymplus{1}$ and
$y=\kappa^{n-t^*}\Ymplus{2}$, so that $x\geq y$. By definition, on
$H_{m,\sN}(\vep)$,
    $$
    \kappa^{t^*} \Ymmin{1}+
    \kappa^{n-t^*}\Ymmin{2}\leq (1+\vep^2)\log N.
    $$
Since, on $E_{m,\sN}$, we have that $Y^{\smallsup{1,N}}_m\geq
\vep$,
    \eq
    \Ymplus{1}
    \leq \frac{\vep+\vep^3}{\vep-\vep^3}\Ymmin{1}
    =\frac{1+\vep^2}{1-\vep^2}\Ymmin{1},
    \en
and likewise for $\Ymplus{2}$. Therefore, we obtain that on
$E_{m,\sN}\cap F_{m,\sN}\cap H_{m,\sN},$ and with $\eps$
sufficiently small,
    \eq
    \label{x+ybd}
    x+y \leq \frac{1+\vep^2}{1-\vep^2}\big[\kappa^{t^*} \Ymmin{1}
    +\kappa^{n-t^*}\Ymmin{2}
    \big]\leq \frac{(1+\eps^2)^2}{1-\eps^2}\log N\leq (1+\vep)\log N.
    \en
Moreover, by (\ref{x/ybd}), we have that
    \eq
    \label{x/ybd2}
    1\leq \frac{x}{y}\leq \kappa.
    \en
Hence, on $E_{m,\sN}\cap F_{m,\sN}\cap H_{m,\sN}$, we have
    \eq
    \label{xub}
    x=\frac{x+y}{1+\frac{y}{x}} \leq (1+\vep) \frac1{1+\kappa^{-1}}\log N
    =\frac{1+\vep}{\tau-1}\log{N},
    \en
when $\vep>0$ is sufficiently small. We claim that if
    \eq
    \label{xaim}
    x=\kappa^{t^*} \Ymplus{1}\leq \frac{1-\vep}{\tau-1}
    \log{N},
    \en
then $k^*=t^*-1\in {\cal B}_{\sN}(\vep,k_{\sN})$, so that
(\ref{fullmin1}) follows. Indeed, we use (\ref{xaim}) to see that
    \eq
    \label{boundYm1}
    \kappa^{k^*+1}Y^{\smallsup{1,N}}_m=\kappa^{t^*}Y^{\smallsup{1,N}}_m
    \leq \kappa^{t^*}\Ymplus{1}
    \leq \frac{1-\vep}{\tau-1} \log{N},
    \en
so that the first bound in (\ref{BN-rew}) is satisfied. The second
bound is satisfied, since
    \eq
    \kappa^{k_{\sN}-k^*}Y^{\smallsup{2,N}}_m=\kappa^{n-t^*}Y^{\smallsup{2,N}}_m\leq
    \kappa^{n-t^*} \Ymplus{2}=y\leq x\leq \frac{1-\vep}{\tau-1}
    \log{N},
    \en
where we have used $n=k_{\sN}+1$, and (\ref{xaim}).

Thus, in order to show that (\ref{fullmin1}) holds with
probability close to $1$, we have to show that the probability of
the intersection of the events $\{H_{\sN}>k_{\sN}\}$ and
    \eqalign
    \label{voorwaarde}
    {\cal E}_{m,\sN}={\cal E}_{m,\sN}(\vep)&=\Big\{\exists t:
    \frac{1-\vep}{\tau-1} \log N <\kappa^{t} \Ymplus{1}
    \le\frac{1+\vep}{\tau-1}\log N,\\
    &\qquad \qquad\kappa^{t} \Ymplus{1}+
    \kappa^{n-t}\Ymplus{2}\leq (1+\vep)\log N \Big\},\nn
    \enalign
can be made arbitrarily small by choosing $\vep$ close to $0$,
when first $N \to \infty$ and then $m \to \infty$. That is the
content of the following lemma, whose proof is deferred to Section
\ref{sec-pflemmas}:

    \begin{lemma}
    \label{lem-calE} For $k_{\sN}$ as in (\ref{kNldef}),
    \eq
    \limsup_{\vep\downarrow 0} \limsup_{m\rightarrow \infty}
    \limsup_{N\rightarrow\infty}
    \prob(E_{m,\sN}(\vep)\cap {\cal E}_{m,\sN}(\vep)\cap \{H_{\sN}>k_{\sN}\})
    =0.\nn
    \en
    \end{lemma}

\noindent Therefore, we finally arrive at
    \eqalign
    &\prob\big(\{H_{\sN}>k_{\sN}\} \cap E_{m,\sN}\big)\\
    &\quad\leq
    \expec\Big[{\bf 1}_{E_{m,\sN}}
    \exp\Big\{-\lambda_{\sN}\exp \Big[\min_{0\leq k_1< k_{\sN}}
    \big(\kappa^{k_1+1}\Ymplus{1}
    +\kappa^{k_{\sN}-k_1}\Ymplus{2}\big)-\log L_{\sN}\Big]
    \Big\}
    \Big]+o_{\sN,m,\vep}(1),\nn
    \enalign
and at a similar lower bound where $\Ymplus{i}$ is replaced by
$\Ymmin{i}$.

\subsection{Evaluating the limit}
\label{sec-step5} The final argument consists of letting $N\to
\infty$ and then $m \to \infty$. The argument has to be performed
with $\Ymplus{i}$ and $\Ymmin{1}$ separately, after which we let
$\vep \downarrow 0$. Since the precise value of $\vep$ plays no
role in the derivation, we only give the derivation for $\vep=0$.
Observe that
    \begin{eqnarray}
    \label{basis asymp}
    &&\min_{0\leq k_1< k_{\sN}}
     (\kappa^{k_1+1}Y^{\smallsup{1,N}}_m+\kappa^{k_{\sN}-k_1}Y^{\smallsup{2,N}}_m)-\log
    L_{\sN}\nonumber\\
    &&\qquad =\kappa^{\lceil k_{\sN}/2 \rceil}\min_{0\leq k_1< k_{\sN}}\left(
     \kappa^{k_1+1-\lceil k_{\sN}/2 \rceil}Y^{\smallsup{1,N}}_m
     +\kappa^{\lfloor k_{\sN}/2 \rfloor-k_1}Y^{\smallsup{2,N}}_m
     -\kappa^{-\lceil k_{\sN}/2 \rceil}\log L_{\sN}\right)\nonumber\\
    &&\qquad = \kappa^{\lceil k_{\sN}/2\rceil} \min_{-\lceil
    k_{\sN}/2\rceil+1\le t < \lfloor k_{\sN}/2\rfloor+1} (\kappa^{t}
    Y_{m}^{\smallsup{1,N}} +\kappa^{c_{l}-t} Y_{m}^{\smallsup{2,N}}
    -\kappa^{-\lceil k_{\sN}/2 \rceil}\log L_{\sN}),
    \end{eqnarray}
where $t=k_1+1-\lceil k_{\sN}/2 \rceil$, $c_{k_{\sN}}=c_l=\lfloor
l/2 \rfloor-\lceil l/2\rceil+1=1_{\{\mbox{$l$ is even}\}}$. We
further rewrite, using (\ref{kNldef}) and the definition of
$a_{\sN}$ in Theorem \ref{thm-tau(2,3)},
    \eq
    \label{kappaaNident}
    \kappa^{-\lceil k_{\sN}/2 \rceil}\log L_{\sN}
    =\kappa^{\frac{\log\log{N}}{\log \kappa}-\lfloor \frac{\log\log{N}}{\log \kappa}\rfloor -\lceil l/2\rceil}
    \frac{\log L_{\sN}}{\log{N}}
    =\kappa^{-a_{\sN} -\lceil l/2\rceil}
    \frac{\log L_{\sN}}{\log{N}}.
    \en
From Lemma \ref{Lemma_mintZ}, for $N\to \infty$ and on the event
$E_{m,\sN}$,
    \begin{equation}
    \label{minoverZ}
    \min_{-\lceil k_{\sN}/2\rceil+1\le t \le \lfloor k_{\sN}/2\rfloor}
    (\kappa^{t} Y_{m}^{\smallsup{1,N}} +\kappa^{c_{l
    }-t}
    Y_{m}^{\smallsup{2,N}}) =\min_{t\in \Z} (\kappa^{t}
    Y_{m}^{\smallsup{1,N}} +\kappa^{c_{l}-t} Y_{m}^{\smallsup{2,N}}),
    \end{equation}
because $Y_{m}^{\smallsup{i,N}}\in [\vep, \vep^{-1}]$ on
$E_{m,\sN}$. We define
    \eq
    W_{m,\sN}(k_{\sN})=\min_{t\in \Z}
    (\kappa^{t} Y_{m}^{\smallsup{1,N}}
    +\kappa^{c_{l}-t} Y_{m}^{\smallsup{2,N}})
    -\kappa^{-a_{\sN} -\lceil l/2\rceil}
    \frac{\log L_{\sN}}{\log{N}},
    \en
and, for $\vep>0$,
    \eq
    \widetilde F_{\sN}=\widetilde F_{\sN}(l,\vep)=\big\{W_{m,\sN}(k_{\sN})>\vep\big\},
    \qquad\qquad
    \widetilde G_{\sN}=\widetilde G_{\sN}(l,\vep)=\big\{W_{m,\sN}(k_{\sN})<-\vep\big\}.
    \label{tildeFGdef}
    \en
Observe that
    \eq\kappa^{\lceil k_{\sN}/2\rceil}W_{m,\sN}(k_{\sN})\cdot {\bf 1}_{\widetilde F_{\sN}}
    \geq \kappa^{\lceil k_{\sN}/2\rceil} \vep,
    \qquad \qquad \kappa^{\lceil k_{\sN}/2\rceil}W_{m,\sN}(k_{\sN})\cdot {\bf 1}_{\widetilde G_{\sN}}
    \leq \kappa^{\lceil k_{\sN}/2\rceil}(-\vep).
    \label{bdFGN}
    \en
We split
    \begin{eqnarray}
    \prob\big(\{H_{\sN}>k_{\sN}\} \cap E_{m,\sN}\big)&=&\expec\left[{\bf 1}_{E_{m,\sN}}
    \exp\big[-\lambda_{\sN}e^{\kappa^{\lceil k_{\sN}/2\rceil}W_{m,\sN}(k_{\sN})}\big]\right]
    +o_{\sN,m,\vep}(1)
    \nn\\
    &=&\prob(\widetilde G_{\sN}\cap E_{m,\sN}) +
    I_{\sN}+J_{\sN}+K_{\sN}+o_{\sN,m,\vep}(1),
    \end{eqnarray}
where we define
    \begin{eqnarray}
    \shift I_{\sN}&\!\!=\!\!&
    \expec\left[\exp\big[-\lambda_{\sN}e^{\kappa^{\lceil k_{\sN}/2\rceil}W_{m,\sN}(k_{\sN})}\big]
    {\bf 1}_{\widetilde F_{\sN}\cap E_{m,\sN}}\right],\\
    \shift J_{\sN}&\!\!=\!\!&\expec
    \left[\Big(\exp\big[-\lambda_{\sN}e^{\kappa^{\lceil k_{\sN}/2\rceil}W_{m,\sN}(k_{\sN})}\big]-1\Big)
    {\bf 1}_{\widetilde G_{\sN}\cap E_{m,\sN}}\right],\\
    \shift K_{\sN}&\!\!=\!\!&\expec\left[\exp\big[- \lambda_{\sN}
    e^{\kappa^{\lceil k_{\sN}/2\rceil}W_{m,\sN}(k_{\sN})}\big]
    {\bf 1}_{\widetilde F_{\sN}^c\cap
    \widetilde G_{\sN}^c \cap E_{m,\sN}}\right].
    \end{eqnarray}
We first show that $I_{\sN},J_{\sN}$ and $K_{\sN}$ are error
terms, and then prove convergence of $\prob(\widetilde G_{\sN}\cap
E_{m,\sN})$.

We start by bounding $I_{\sN}$. By the first bound in
(\ref{bdFGN}), for every $\vep>0$, and since $\lambda_{\sN}\geq
\frac12$,
    \eq
    \limsup_{N\rightarrow \infty} I_{\sN} \leq \limsup_{N\rightarrow \infty}
    \exp\big\{-\frac{1}{2}\exp\{\kappa^{\lceil k_{\sN}/2\rceil} \vep\}\big\}
    =0.
    \en
Similarly, by the second bound in (\ref{bdFGN}), for every
$\vep>0$, and since $\lambda_{\sN}\leq 4k_{\sN}$, we can bound
$J_{\sN}$ as
    \eq
    \limsup_{N\rightarrow \infty} |J_{\sN}| \leq \limsup_{N\rightarrow \infty}
    \expec\Big[1-\exp\big\{-4k_{\sN}\exp\{-\kappa^{\lceil k_{\sN}/2\rceil} \vep\}\big\}\Big]=0.
    \en
Finally, we bound $K_{\sN}$ by
    \eq
    K_{\sN} \leq \prob (\widetilde F_{\sN}^c\cap \widetilde G_{\sN}^c\cap E_{m,\sN}).
    \en
We will next show that
    \eqalign
    &\prob(\widetilde G_{\sN}\cap E_{m,\sN})\label{aim1}=\prob\big(\min_{t\in \Z} (\kappa^{t} Y^{\smallsup{1}}
    +\kappa^{c_l-t} Y^{\smallsup{2}})-\kappa^{-a_{\sN}-\lceil
    l/2 \rceil}\leq 0, Y^{\smallsup{1}}Y^{\smallsup{2}}>0\big) +o_{\sN,m,\vep}(1),
    \enalign
and
    \eqalign
    \prob (\widetilde F_{\sN}^c\cap \widetilde G_{\sN}^c\cap E_{m,\sN})
    &=o_{\sN,m,\vep}(1).
    \label{aim2}
    \enalign
Equation (\ref{aim2}) is the content of the following lemma, whose
proof is deferred to Section \ref{sec-pflemmas}:
    \begin{lemma} For all $l$,
    \label{lem-FcGEbd}
    \eqalign
    \limsup_{\vep\downarrow 0} \limsup_{m\rightarrow \infty}
    \limsup_{N\rightarrow\infty}\prob (\widetilde F_{\sN}(l,\vep)^c\cap \widetilde G_{\sN}(l,\vep)^c\cap E_{m,\sN}(\eps))
    &=0.\nn
    \enalign
    \end{lemma}
\vskip0.3cm

\noindent We now prove (\ref{aim1}). From the definition of
$\widetilde G_{\sN}$,
    \eq
    \widetilde G_{\sN}
    \cap E_{m,\sN}
    =\Big\{\min_{t\in \Z} (\kappa^{t} Y_{m}^{\smallsup{1,N}}
    +\kappa^{c_l-t} Y_{m}^{\smallsup{2,N}})- \kappa^{-a_{\sN}-\lceil
    l/2\rceil}\frac{\log{L_{\sN}}}{\log{N}}<-\vep, Y_{m}^{\smallsup{i,N}}\in [\vep, \vep^{-1}]\Big\}.
    \en
By Proposition \ref{prop-caft} and the fact that $L_{\sN}=\mu
N(1+o(1))$ with probability $1-o_{\sN}(1)$,
    \eq
    \label{tildeGprobeq}
    \prob(\widetilde G_{\sN}\cap E_{m,\sN})
    =\prob\Big(\min_{t\in \Z} (\kappa^{t} Y_{m}^{\smallsup{1}}
    +\kappa^{c_l-t} Y_{m}^{\smallsup{2}})- \kappa^{-a_{\sN}-\lceil
    l/2\rceil}<-\vep, Y_{m}^{\smallsup{i}}\in [\vep, \vep^{-1}]\Big)
    +o_{\sN}(1).
    \en
Since $Y_{m}^{\smallsup{i}}$ converges to $Y^{\smallsup{i}}$
almost surely, $\sup_{s\geq
m}|Y_{s}^{\smallsup{i}}-Y^{\smallsup{i}}|$ converges to 0 a.s.\ as
$m\rightarrow \infty$. Therefore,
    \eq
    \label{YmbyYrepl}
    \prob(\widetilde G_{\sN}\cap E_{m,\sN})
    =\prob\Big(M_l- \kappa^{-a_{\sN}-\lceil
    l/2\rceil}<-\vep, Y^{\smallsup{i}}\in [\vep, \vep^{-1}]\Big)
    +o_{\sN, m, \vep}(1),
    \en
where we define
    \eq
    M_l=\min_{t\in \Z} (\kappa^{t} Y^{\smallsup{1}}
    +\kappa^{c_l-t} Y^{\smallsup{2}}).
    \en
Moreover, since $Y^{\smallsup{1}}$ has a density on $(0,\infty)$
and an atom at 0 (see \cite{davies}),
    \eq
    \prob(Y^{\smallsup{1}}\not\in [\vep, \vep^{-1}], Y^{\smallsup{1}}>0)
    =o_{\vep}(1),\nn
    \en
so that, in turn,
    \eqalign
    \prob(\widetilde G_{\sN}\cap E_{m,\sN})
    &=\prob\Big(M_l- \kappa^{-a_{\sN}-\lceil
    l/2\rceil}<-\vep, Y^{\smallsup{1}}Y^{\smallsup{2}}>0\Big)
    +o_{\sN, m, \vep}(1)\nn\\
    &=q^2 \prob\Big(M_l- \kappa^{-a_{\sN}-\lceil
    l/2\rceil}<-\vep\Big|Y^{\smallsup{1}}Y^{\smallsup{2}}>0\Big)
    +o_{\sN, m, \vep}(1).
    \label{tildeGprob}
    \enalign

Recall from Section \ref{sec-BP} that for any $l$ fixed,
conditionally on $Y^{\smallsup{1}}Y^{\smallsup{2}}>0$, the random
variable $M_l$ has a density. We denote this density by $f_2$ and
the distribution function by $F_2$. Also, $\kappa^{-a_{\sN}-\lceil
l/2 \rceil} \in I_l=[\kappa^{-\lceil l/2 \rceil}, \kappa^{-\lceil
l/2 \rceil+1}]$. Then,
    \eq
    \label{distribd}
    \prob\Big(-\vep \leq M_l- \kappa^{-a_{\sN}-\lceil
    l/2\rceil}<0\Big)
    \leq \sup_{a\in I_l} [F_2(a)-F_2(a-\vep)].
    \en
The function $F_2$ is continuous on $I_l$, so that in fact $F_2$
is uniformly continuous on $I_l$, and we conclude that
    \eq
    \label{densest1}
    \limsup_{\vep\downarrow 0}\sup_{a\in I_l} [F_2(a)-F_2(a-\vep)]=0.
    \en
We conclude the results of Sections
\ref{sec-step1}--\ref{sec-step5} in the following corollary:
    \begin{corr}
    \label{cor-Y'spos}
    For all $l$, with $k_{\sN}$ as in (\ref{kNldef}),
    \eq
    \prob\big(\{H_{\sN}>k_{\sN}\} \cap \{Y^{\smallsup{1}}Y^{\smallsup{2}}>0\}\big)\nn\\
    =q^2 \prob\Big(M_l- \kappa^{-a_{\sN}-\lceil
    l/2\rceil}\leq 0\Big|Y^{\smallsup{1}}Y^{\smallsup{2}}>0\Big)
    +o_{\sN, m, \vep}(1).
    \en
    \end{corr}
\vskip0.3cm

\noindent We now come to the conclusion of the proofs of Theorems
\ref{thm-tau(2,3)} and \ref{thm-ll}. We combine the results in
Corollaries \ref{cor-Y's0} and \ref{cor-Y'spos}, together with the
fact that $q_m=q+o(1)$ as $m\rightarrow \infty$, to obtain that,
with $\ktauN=2\left\lfloor \frac{\log \log
N}{|\log(\tau-2)|}\right\rfloor$,
    \eq
    \prob\left(H_{\sN}> \ktauN+l\right)
    = 1-q^2 + q^2\prob\left(M_l- \kappa^{-a_{\sN}-\lceil
    l/2\rceil}\leq 0\Big|Y^{\smallsup{1}}Y^{\smallsup{2}}>0\right)
    +o_{\sN, m, \vep}(1).
    \en
Therefore,
    \eqalign
    \label{probunconn}
    \prob\left(H_{\sN}\leq \ktauN+l\right)
    &=q^2\prob\left(\min_{t\in \Z} (M_l>\kappa^{-a_{\sN}-\lceil
    l/2 \rceil}\big|Y^{\smallsup{1}}Y^{\smallsup{2}}>0\right)
    +o_{\sN, m, \vep}(1).
    \enalign
When $l\rightarrow \infty$, we claim that (\ref{probunconn})
implies that, when $N\rightarrow \infty$,
    \eq
    \label{probconn}
    \prob\left(H_{\sN}<\infty\right)=q^2+o(1).
    \en
Indeed, to see (\ref{probconn}), we prove upper and lower bounds.
For the lower bound, we use that for any $l\in \mathbb{Z}$
    \[
    \prob\left(H_{\sN}<\infty\right)
    \geq \prob\left(H_{\sN}\leq \ktauN+l\right),
    \]
and let $l\rightarrow \infty$ in (\ref{probunconn}), together with
the fact that $\kappa^{-a_{\sN}-\lceil l/2 \rceil}\rightarrow 0$
as $l\rightarrow \infty$. For the upper bound, we split
    \[
    \prob\left(H_{\sN}<\infty\right)
    =\prob\left(\{H_{\sN}<\infty\}\cap \{Y_m^{\smallsup{1,N}}Y_m^{\smallsup{2,N}}=0\}\right)
    +\prob\left(\{H_{\sN}<\infty\}\cap \{Y_m^{\smallsup{1,N}}Y_m^{\smallsup{2,N}}>0\}\right).
    \]
The first term is bounded by $\prob(H_{\sN}\leq m-1)=o_{\sN}(1)$,
by Lemma \ref{lem-Hsmall}. The second term is bounded from above
by, using Proposition \ref{prop-caft},
    \eq
    \prob\left(\{H_{\sN}<\infty\}\cap \{Y_m^{\smallsup{1,N}}Y_m^{\smallsup{2,N}}>0\}\right)
    \leq \prob\left(Y_m^{\smallsup{1,N}}Y_m^{\smallsup{2,N}}>0\right)
    =q_m^2+o_{\sN}(1),
    \en
which converges to $q^2$ as $m\rightarrow \infty$. This proves
(\ref{probconn}). We conclude from (\ref{probunconn}) and
(\ref{probconn}) that
    \eq
    \prob\left(H_{\sN}\leq \ktauN+l\Big|H_{\sN}<\infty\right)
    =\prob\left(M_l\geq (\tau-2)^{a_{\sN}+\lceil
    l/2 \rceil}\Big|Y^{\smallsup{1}}Y^{\smallsup{2}}>0\right)
    +o(1).
    \en
This completes the proofs of Theorems \ref{thm-tau(2,3)} and
\ref{thm-ll}. \qed


\section{Proofs of Lemmas \ref{lem-BNempty}, \ref{lem-calE} and \ref{lem-FcGEbd}}
\label{sec-pflemmas} In this section, we prove the three lemmas
used in Section \ref{sec-pftau(2,3)}. The proofs are similar in
nature. Denote
    \eq
    \{k\in \partial {\cal T}_m^{\smallsup{N}}\}
    =\{k\in {\cal T}_m^{\smallsup{N}}\}\cap \{k+1\not\in {\cal T}_m^{\smallsup{N}}\}.
    \en
We will make essential use of the following consequence of
Propositions \ref{prop-caft} and \ref{prop-weakconv2bsec}:

\begin{lemma}
\label{lem-bdsexp} Assume that Propositions \ref{prop-caft} and
\ref{prop-weakconv2bsec} hold. Then, for any $u>0$, and $i=1,2$,
    \eq
    \prob\big(\{k\in {\cal T}_m^{\smallsup{i,N}}\}
    \cap \{\vep\leq Y_m^{\smallsup{i,N}}\leq \vep^{-1}\}
    \cap \{Z_k^{\smallsup{i,N}}\in [N^{u(1-\vep)}, N^{u(1+\vep)}]\}\big)=
    o_{\sN, m, \vep}(1).
    \en
Consequently,
    \eq
    \label{ingred}
    \prob\big(\{k\in \partial{\cal T}_m^{\smallsup{i,N}}\}
    \cap \{\vep\leq Y_m^{\smallsup{i,N}}\leq \vep^{-1}\}
    \cap \{Z_k^{\smallsup{i,N}}\leq N^{\frac{1}{\kappa(\tau-1)}+\vep}\}\big)=
    o_{\sN, m, \vep}(1).
    \en
\end{lemma}

\proof By Proposition \ref{prop-weakconv2bsec}, \whps
    \eq
    Y_k^{\smallsup{i,N}}\leq Y_m^{\smallsup{i,N}}+\vep^3\leq Y_m^{\smallsup{i,N}}(1+\vep^2),
    \en
where the last inequality follows from $Y_m^{\smallsup{i,N}}\leq
\vep^{-1}$. Therefore, also
    \eq
    Y_m^{\smallsup{i,N}}\geq Y_k^{\smallsup{i,N}}(1-2\vep^2),
    \en
when $\vep$ is so small that $(1+\vep^2)^{-1}\geq 1-2\vep^2$. In a
similar way, we conclude that
    \eq
    \label{Ymub}
    Y_m^{\smallsup{i,N}}\leq Y_k^{\smallsup{i,N}}(1+2\vep^2).
    \en
Furthermore, the event $Z_k^{\smallsup{i,N}}\in [N^{u(1-\vep)},
N^{u(1+\vep)}]$ is equivalent to
    \eq
    (1-\vep)u\kappa^{-k}\log{N}\leq Y_k^{\smallsup{i,N}} \leq (1+\vep)u\kappa^{-k} \log{N}.
    \en
Therefore, we obtain that, with $u_{k,\sN}=u\kappa^{-k} \log{N}$,
    \eq
    Y_m^{\smallsup{i,N}}\leq (1+2\vep^2)(1+\vep) u\kappa^{-k} \log{N}\leq (1+2\vep) u_{k,\sN}.
    \en
Similarly, we obtain
    \eq
    Y_m^{\smallsup{i,N}}\geq (1-2\vep^2)(1-\vep)u\kappa^{-k} \log{N}\geq (1-2\vep) u_{k,\sN}.
    \en
We conclude that the events $k\in {\cal T}_m^{\smallsup{i,N}},
\vep\leq Y_m^{\smallsup{i,N}}\leq \vep^{-1} $ and
$Z_k^{\smallsup{i,N}}\in [N^{u(1-\vep)}, N^{u(1+\vep)}]$ imply
    \eq
    Y_m^{\smallsup{i,N}}\in u_{k,\sN}[1-2\vep, 1+2\vep]\equiv
    [u_{k,\sN}(1-2\vep), u_{k,\sN}(1+2\vep)].
    \en
Since $\vep\leq Y_m^{\smallsup{N}}\leq \vep^{-1}$, we therefore
must also have (when $\vep$ is so small that $1-2\vep\geq \frac
12$),
    \eq
    u_{k,\sN}\in [\frac{\vep}{2},\frac{2}{\vep}].
    \en
Therefore,
    \eqalign
    &\limsup_{\vep\downarrow 0}\limsup_{m\rightarrow \infty}
    \limsup_{N\rightarrow\infty}\prob\big(\{k\in {\cal T}_m^{\smallsup{i,N}}\}
    \cap \{\vep\leq Y_m^{\smallsup{i,N}}\leq \vep^{-1}\}
    \cap \{Z_k^{\smallsup{i,N}}\in [N^{u(1-\vep)}, N^{u(1+\vep)}]\}\big)\\
    &\qquad \leq
    \limsup_{\vep\downarrow 0}\limsup_{m\rightarrow \infty}
    \limsup_{N\rightarrow\infty}\sup_{x\in [\frac{\vep}{2},\frac{2}{\vep}]}
    \prob\big(Y_m^{\smallsup{i,N}}\in x[1-2\vep, 1+2\vep]\big).\nn
    \enalign

Since $Y_m^{\smallsup{i,N}}=Y_m^{\smallsup{i}}$ \whpl by
Proposition \ref{prop-caft}, we arrive at
    \eqalign
    &\limsup_{\vep\downarrow 0}\limsup_{m\rightarrow \infty}
    \limsup_{N\rightarrow\infty}\prob\big(\{k\in {\cal T}_m^{\smallsup{i,N}}\}
    \cap \{\vep\leq Y_m^{\smallsup{i,N}}\leq \vep^{-1}\}
    \cap \{Z_k^{\smallsup{i,N}}\in [N^{u(1-\vep)}, N^{u(1+\vep)}]\}\big)\\
    &\qquad \leq
    \limsup_{\vep\downarrow 0}\limsup_{m\rightarrow \infty}
    \sup_{x\in [\frac{\vep}{2},\frac{2}{\vep}]}
    \prob\big(Y_m^{\smallsup{i}}\in x[1-2\vep, 1+2\vep]\big).\nn
    \enalign
We next use that $Y_m^{\smallsup{i}}$ converges to
$Y^{\smallsup{i}}$ almost surely as $m\rightarrow \infty$ to
arrive at
    \eqalign
    &\limsup_{\vep\downarrow 0}\limsup_{m\rightarrow \infty}
    \limsup_{N\rightarrow\infty}\prob\big(\{k\in {\cal T}_m^{\smallsup{i,N}}\}
    \cap \{\vep\leq Y_m^{\smallsup{i,N}}\leq \vep^{-1}\}
    \cap \{Z_k^{\smallsup{i,N}}\in [N^{u(1-\vep)}, N^{u(1+\vep)}]\}\big)\\
    &\qquad \leq
    \limsup_{\vep\downarrow 0} \sup_{x\in [\frac{\vep}{2},\frac{2}{\vep}]}
    \prob\big(Y^{\smallsup{i}}\in x[1-2\vep, 1+2\vep]\big)
    \leq \limsup_{\vep\downarrow 0} \sup_{x>0} [F_1(x(1+2\vep))-F_1(x(1-2\vep))],
    \nn
    \enalign
where $F_1$ denotes the distribution function of
$Y^{\smallsup{i}}$, which is continuous for $x>0$. Moreover,
    \eq
    \lim_{x\rightarrow \infty} 1-F_1(x)=0.
    \en
Therefore, uniformly in $\vep<1/4$,
    \eq
    \sup_{x>K} [F_1(x(1+2\vep))-F_1(x(1-2\vep))]
    \leq 2 \sup_{x>K} [1-F_1(x(1-2\vep))]\rightarrow 0, \qquad (K\rightarrow \infty),
    \en
so that
    \eqalign
    &\sup_{x>0} [F_1(x(1+2\vep))-F_1(x(1-2\vep))]\\
    &\qquad\leq \sup_{0<x\leq K} [F_1(x(1+2\vep))-F_1(x(1-2\vep))]
    +\sup_{x>K} [F_1(x(1+2\vep))-F_1(x(1-2\vep))]\nn\\
    &\qquad =
    \sup_{0<x\leq K} [F_1(x(1+2\vep))-F_1(x(1-2\vep))]+o(1).\nn
    \enalign
Therefore,
    \eq
    \label{unifcont}
    \limsup_{\vep\downarrow 0} \sup_{x>0} [F_1(x(1+2\vep))-F_1(x(1-2\vep))]
    =\limsup_{K\uparrow \infty} \limsup_{\vep\downarrow 0}\sup_{0<x\leq K}[F_1(x(1+2\vep))-F_1(x(1-2\vep))]=0,
    \en
since $F_1$ is uniformly continuous on $(0,K]$. This completes the
proof of the first statement in Lemma \ref{lem-bdsexp}.

We turn to the second statement. The event that $k\in \partial
{\cal T}_m^{\smallsup{N}}$ implies that
    \eq
    Y_m^{\smallsup{i,N}}\geq \frac{1-\vep^2}{\tau-1}\kappa^{-(k+1)} \log{N}.
    \en
By (\ref{Ymub}), we therefore conclude that when $\vep$ is
sufficiently small
    \eq
    Y_k^{\smallsup{i,N}}\geq \frac{1-\vep}{\tau-1}\kappa^{-(k+1)} \log{N},
    \en
which is equivalent to
    \eq
    Z_k^{\smallsup{i,N}}\geq N^{\frac{1-\vep}{\kappa(\tau-1)}}
    \geq N^{\frac{1}{\kappa(\tau-1)}-\vep}.
    \en
Therefore,
    \eqalign
    &\limsup_{\vep\downarrow 0}\limsup_{m\rightarrow \infty}
    \limsup_{N\rightarrow\infty}\prob\big(\{k\in \partial{\cal T}_m^{\smallsup{i,N}}\}
    \cap \{\vep\leq Y_m^{\smallsup{i,N}}\leq \vep^{-1}\}
    \cap \{Z_k^{\smallsup{i,N}}\leq N^{\frac{1}{\kappa(\tau-1)}+\vep}\}\big)\\
    &\leq \limsup_{\vep\downarrow 0}\limsup_{m\rightarrow \infty}
    \limsup_{N\rightarrow\infty}\prob\big(\{k\in {\cal T}_m^{\smallsup{i,N}}\}
    \cap \{\vep\leq Y_m^{\smallsup{i,N}}\leq \vep^{-1}\}
    \cap \{Z_k^{\smallsup{i,N}}\in [N^{\frac{1}{\kappa(\tau-1)}-\vep},
    N^{\frac{1}{\kappa(\tau-1)}+\vep}]\}\big)\nn\\
    &=0,\nn
    \enalign
which follows from the first statement in Lemma \ref{lem-bdsexp}
with $u=\frac{1}{\kappa(\tau-1)}$. \qed \vskip0.5cm

\noindent {\bf Proof of Lemma \ref{lem-BNempty}.} By
(\ref{Fbdpf}), it suffices to prove that
    \eq
    \label{aimB}
    \limsup_{\vep\downarrow 0}\limsup_{m\rightarrow \infty}
    \limsup_{N\to\infty}\prob(\{H_{\sN}>k_{\sN}\} \cap E_{m,\sN}\cap F_{m,\sN}
    \cap\{{\cal B}_{\sN}(\vep,k_{\sN})= \varnothing\} )=0,
    \en
which shows that in considering the event $\{H_{\sN}>k_{\sN}\}\cap
E_{m,\sN}\cap F_{m,\sN}$, we may assume that ${\cal
B}_{\sN}(\vep,k_{\sN})\neq \varnothing$.

We define the random variable $l^*$ by
    \eq
    l^*=\sup\{k: {\cal B}_{\sN}(\vep,k)\neq \varnothing\}.
    \en
Observe that if ${\cal B}_{\sN}(\vep,k)=\varnothing$, then ${\cal
B}_{\sN}(\vep,k+1)=\varnothing$, so that $l^*$ is well defined.
Indeed, if $l\in {\cal B}_{\sN}(\vep,k+1)$ and $l\neq m$, then
$l-1\in {\cal B}_{\sN}(\vep,k)$. If, on the other hand, ${\cal
B}_{\sN}(\vep,k+1)=\{m\}$, then also $m\in {\cal
B}_{\sN}(\vep,k)$. Since $\{{\cal B}_{\sN}(\vep,k_{\sN})=
\varnothing\} =\{k_{\sN}\geq l^*+1\}$, we therefore have
    \eq
    \{{\cal B}_{\sN}(\vep,k_{\sN})= \varnothing\}
    =\{l^*<k_{\sN}\}
    =\{l^*\leq k_{\sN}-2\}\DDcup \{l^*=k_{\sN}-1\}.
    \en
We deal with each of the two events separately. We start with the
first.

Since the sets ${\cal B}_{\sN}(\vep,k)$ are $Z_m$-measurable, we
obtain, as in (\ref{rewprob3}),
    \eqalign
    \label{bdBNempt1}
    \prob(\{H_{\sN}>k_{\sN}\} \cap E_{m,\sN}\cap F_{m,\sN}
    \cap\{l^*\leq k_{\sN}-2\})
    &\leq \prob(\{H_{\sN}>l^*+2\} \cap E_{m,\sN}\cap F_{m,\sN})\\
    &=\expec\Big[{\bf 1}_{E_{m,\sN}\cap F_{m,\sN}}
    P_m(l^*+2,k_1)\Big] +o_{\sN,m,\vep}(1).\nn
    \enalign
We then use (\ref{rewprobub}) to bound
    \eqalign
    \expec\Big[{\bf 1}_{E_{m,\sN}\cap F_{m,\sN}}
    P_m(l^*+2,k_1)\Big]
    &\leq\expec\Big[{\bf 1}_{E_{m,\sN}\cap F_{m,\sN}}\min_{k_1\in {\cal B}_{\sN}(\vep,l^*)}
    \exp\{-\frac{Z_{k_1+1}^{\smallsup{1,N}}
    Z_{l^*+2-k_1}^{\smallsup{2,N}}}{2L_{\sN}}\}\Big].
    \label{l*small}
    \enalign
Now, since ${\cal B}_{\sN}(\vep,l^*)\neq \varnothing$, we can pick
$k_1$ such that $k_1-1\in {\cal B}_{\sN}(\vep,l^*)$. Since ${\cal
B}_{\sN}(\vep,l^*+1)=\varnothing$, we have $k_1-1\notin {\cal
B}_{\sN}(\vep,l^*+1)$, implying $l^*+1-k_1\in {\cal
T}_m^{\smallsup{2,N}}$ and $l^*+2-k_1\notin {\cal
T}_m^{\smallsup{2,N}}$ so that, by (\ref{Zlb}),
$Z_{l^*+2-k_1}^{\smallsup{2,N}} \geq N^{\frac{1-\vep}{\tau-1}}$.

Similarly, since $k_1\not\in {\cal B}_{\sN}(\vep,l^*+1)$ we have that
$k_1\in {\cal T}_m^{\smallsup{1,N}}$ and $k_1+1\notin {\cal
T}_m^{\smallsup{1,N}},$ so that, again by (\ref{Zlb}),
$Z_{k_1+1}^{\smallsup{1,N}}\geq N^{\frac{1-\vep}{\tau-1}}$.
Therefore, since $L_{\sN}\geq N$, \whps,
    \eq
    \label{bdBNempt5}
    \frac{Z_{k_1+1}^{\smallsup{1,N}}
    Z_{l^*+2-k_1}^{\smallsup{2,N}}}{L_{\sN}}
    \geq N^{\frac{2(1-\vep)}{\tau-1}-1},
    \en
and the exponent of $N$ is strictly positive for $\tau \in (2,3)$
and $\vep>0$ small enough. This bounds the contribution in
(\ref{l*small}) due to $\{l^*\leq k_{\sN}-2\}$.

We proceed with the contribution due to $\{l^*=k_{\sN}-1\}$. In
this case, there exists a $k_1$ with $k_1-1\in {\cal
B}_{\sN}(\vep,k_{\sN}-1)$ so that $k_1\in {\cal
T}_m^{\smallsup{1,N}} $ and $k_{\sN}-k_1 \in {\cal
T}_m^{\smallsup{2,N}}$. On the other hand, ${\cal
B}_{\sN}(\vep,k_{\sN})=\varnothing$, which together with $k_1-1\in
{\cal B}_{\sN}(\vep,k_{\sN}-1)$ implies that $k_{\sN}-k_1 \in
{\cal T}_m^{\smallsup{2,N}}$, and $k_{\sN}-k_1+1 \notin {\cal
T}_m^{\smallsup{2,N}}$. Similarly, we obtain that $k_1\in {\cal
T}_m^{\smallsup{1,N}}$ and $k_1+1 \notin {\cal
T}_m^{\smallsup{1,N}}$. Using Proposition \ref{prop-Tub}, we
conclude that $Z_{k_1+1}^{\smallsup{1,N}}\geq
N^{\frac{1-\vep}{\tau-1}}$.

There are two possible cases that we will treat separately: (a)
$Z_{k_{\sN}-k_1}^{\smallsup{2,N}} \leq
N^{\frac{\tau-2}{\tau-1}+\vep}$; and (b)
$Z_{k_{\sN}-k_1}^{\smallsup{2,N}}
>N^{\frac{\tau-2}{\tau-1}+\vep}$. By (\ref{ingred}) and the fact
that $k_{\sN}-k_1 \in \partial{\cal T}_m^{\smallsup{2,N}}$, case
(a) has small probability, so we need to investigate case (b)
only.

In case (b), we can bound
    \eqalign
    \label{bdBNempt4}
    &\shift\prob(\{H_{\sN}>k_{\sN}\} \cap E_{m,\sN} \cap F_{m,\sN}
    \cap\{l^*= k_{\sN}-1\}\cap \{Z_{k_{\sN}-k_1}^{\smallsup{2,N}}>
    N^{\frac{\tau-2}{\tau-1}+\vep}\})\nn\\
    &\quad= \expec\Big[{\bf 1}_{E_{m,\sN}\cap F_{m,\sN}\cap\{l^*= k_{\sN}-1\}}
    {\bf 1}_{\{Z_{k_{\sN}-k_1}^{\smallsup{2,N}}> N^{\frac{\tau-2}{\tau-1}+\vep}\}\cap
    \{k_1-1\in {\cal B}_{\sN}(\vep,k_{\sN}-1)\}} P_m(k_{\sN},k_1)\Big]+o_{\sN,m,\vep}(1),
    \enalign
and again use (\ref{rewprobub}) to obtain
    \eqalign
    \label{bdBNempt3}
    \shift&\prob(\{H_{\sN}>k_{\sN}\} \cap E_{m,\sN} \cap F_{m,\sN}
    \cap\{l^*= k_{\sN}-1\}\cap \{Z_{k_{\sN}-k_1}^{\smallsup{2,N}}> N^{\frac{\tau-2}{\tau-1}+\vep}\})\\
    &\leq \expec\Big[{\bf 1}_{E_{m,\sN}\cap F_{m,\sN}(\vep)\cap\{l^*= k_{\sN}-1\}}
    {\bf 1}_{\{Z_{k_{\sN}-k_1}^{\smallsup{2,N}}> N^{\frac{\tau-2}{\tau-1}+\vep}\}\cap
    \{k_1-1\in {\cal B}_{\sN}(\vep,k_{\sN}-1)\}} \exp\{-\frac{Z_{k_1+1}^{\smallsup{1,N}}
    Z_{k_{\sN}-k_1}^{\smallsup{2,N}}}{2L_{\sN}}\}\Big]+o_{\sN,m,\vep}(1).\nn
    \enalign
We note that by Proposition \ref{prop-Tub} and similarly to (\ref{bdBNempt5}),
    \eq
    Z_{k_1+1}^{\smallsup{1,N}}Z_{k_{\sN}-k_1}^{\smallsup{2,N}}
    \geq N^{\frac{1-\vep}{\tau-1}}N^{\frac{\tau-2}{\tau-1}+\vep}
    =N^{1+\big(1-\frac{1}{\tau-1}\big)\vep},
    \en
and again the exponent is strictly larger than 1, so that,
following the arguments in (\ref{bdBNempt1}--\ref{bdBNempt3}), we
obtain that also the contribution due to case (b) is small. \qed

\vskip0.5cm

\noindent {\bf Proof of Lemma \ref{lem-calE}.} Recall that
$x=\kappa^{t} \Ymplus{1}$ and $y=\kappa^{n-t}\Ymplus{2}$, and that
$x\geq y$. The event ${\cal E}_{m,\sN}$ in (\ref{voorwaarde}) is
equal to
    \eq
    \label{xbds}
    \frac{1-\vep}{\tau-1} \log N \leq x\leq \frac{1+\vep}{\tau-1}\log N,
    \qquad \text{ and }\qquad x+y\leq (1+\vep)\log N.
    \en
Therefore, by (\ref{x/ybd2}),
    \eq
    y\geq \frac{x}{\kappa} \geq (1-\vep)\frac{\tau-2}{\tau-1}\log{N}.
    \en
Also, by the bound on $x+y$ in (\ref{xbds}) and the lower bound on
$x$ in (\ref{xbds}),
    \eq
    y\leq (1+\vep)\log N-x \leq (1+\vep)\log N-\frac{1-\vep}{\tau-1}\log{N}
    = \big(1+\vep\frac{\tau}{\tau-2}\big)\frac{\tau-2}{\tau-1} \log N.
    \en
Therefore, by multiplying the bounds on $x$ and $y$, we obtain
    \eq
    (1-\vep)^2\frac{\tau-2}{(\tau-1)^2}\log^2{N} \leq
    \kappa^{k_{\sN}+1}\Ymplus{1}\Ymplus{2}\leq
    \big(1+\vep\frac{\tau}{\tau-2}\big)(1+\vep)\frac{\tau-2}{(\tau-1)^2} \log^2 N,
    \en
and thus
    \eqalign
    &\prob(E_{m,\sN}\cap{\cal E}_{m,\sN}\cap \{H_{\sN}>k_{\sN}\})
    \leq \prob\Big((1-\vep)^2 \leq
    \frac{\kappa^{k_{\sN}+1}}{c \log^2{N}}\Ymplus{1}\Ymplus{2}\leq
    \big(1+\vep\frac{\tau}{\tau-2}\big)(1+\vep)
    \Big),
    \enalign
where we abbreviate $c=\frac{\tau-2}{(\tau-1)^2}$. We conclude
that
    \eq
    \limsup_{\vep\downarrow 0}\limsup_{m\rightarrow \infty}
    \limsup_{N\to \infty}
    \prob(E_{m,\sN}\cap{\cal E}_{m,\sN}\cap \{H_{\sN}>k_{\sN}\}) =0,
    \en
analogously to the final part of the proof of Lemma
\ref{lem-bdsexp}. \qed \vskip0.5cm

\noindent {\bf Proof of Lemma \ref{lem-FcGEbd}.} We recall that
$M_l=\min_{t\in \Z} (\kappa^{t} Y^{\smallsup{1}} +\kappa^{c_l-t}
Y^{\smallsup{2}})$. We repeat the arguments leading to
(\ref{tildeGprobeq}--\ref{tildeGprob}) to see that, as first $N\to
\infty$ and then $m\to\infty$,
    \eqalign
    \label{PGNeta}
    \prob(\widetilde F^c_{\sN}\cap \widetilde G^c_{\sN}\cap E_{m,\sN})
    &\leq
    \prob\left(-\vep\leq M_l-\kappa^{-a_{\sN}-\lceil
    l/2 \rceil} \leq \vep, Y^{\smallsup{1}}Y^{\smallsup{2}}>0
    \right)+o_{\sN,m}(1)\\
    &=q^2\prob\left(-\vep\leq M_l-\kappa^{-a_{\sN}-\lceil
    l/2 \rceil} \leq \vep\Big|Y^{\smallsup{1}}Y^{\smallsup{2}}>0
    \right)+o_{\sN,m}(1)\nonumber.
    \enalign
Recall from Section \ref{sec-BP} that, conditionally on
$Y^{\smallsup{1}}Y^{\smallsup{2}}>0$, the random variable $M_l$
has a density. Recall that we denoted the distribution function of
$M_l$ given $Y^{\smallsup{1}}Y^{\smallsup{2}}>0$ by $F_2$.
Furthermore, $\kappa^{-a_{\sN}-\lceil l/2 \rceil}\in I_l=
[\kappa^{-\lceil l/2 \rceil},\kappa^{-\lceil l/2 \rceil+1}]$, so
that, uniformly in $N$,
    \eq
    \prob\left(-\vep\leq M_l-\kappa^{-a_{\sN}-\lceil
    l/2 \rceil} \leq \vep\Big|Y^{\smallsup{1}}Y^{\smallsup{2}}>0
    \right)\leq \sup_{u\in I_l} [F_2(u+\vep)-F_2(u-\vep)]=0,
    \nn
    \en
where the conclusion follows by repeating the argument leading to
(\ref{unifcont}). This completes the proof of Lemma
\ref{lem-FcGEbd}.\qed

\renewcommand{\thesection}{\Alph{section}}
\setcounter{section}{0}

\numberwithin{equation}{subsection}
\numberwithin{theorem}{subsection}



\section{Appendix: Proof of Propositions~\ref{prop-caft}--\ref{prop-Tub}}
The appendix is organised as follows. In Section \ref{preplemmas}
we prove three lemmas that are used in Section \ref{sec-caft} to
prove Proposition \ref{prop-caft}. In Section \ref{sec-prep} we
continue with preparations for the proofs of Proposition
\ref{prop-weakconv2bsec} and \ref{prop-Tub}. In this section we
formulate key Proposition \ref{dnziva_P1}, which will be proved in
Section \ref{sec-mainproofapp}. In Section \ref{sec-main idea} we
end the appendix with the proofs of Proposition
\ref{prop-weakconv2bsec} and \ref{prop-Tub}.

\subsection{Some preparatory lemmas}
\label{preplemmas}

In order to prove Proposition \ref{prop-caft}, we make essential
use of three lemmas, that also play a key role in Section
\ref{sec-mainproofapp} below.
The first of these three lemmas investigates
the tail behaviour of $1-G(x)$ under Assumption \ref{ass-gamma}.
Recall that $G$ is the distribution function of the probability mass function
$\{g_j\}$, defined in (\ref{outgoing degree}).
\begin{lemma}
\label{lem-G-x} If $F$ satisfies Assumption \ref{ass-gamma} then
there exists $K_{\tau}>0$ such that for $x$ large enough
    \eq
    \label{gerard1}
    x^{2-\tau-K_{\tau}\gamma(x)}\leq 1-G(x) \leq x^{2-\tau+K_{\tau}\gamma(x)},
    \en
where $\gamma(x)=(\log x)^{\gamma-1}$, $\gamma \in[0,1)$.
\end{lemma}
\proof Using definition~(\ref{outgoing degree}) we rewrite
$1-G(x)$ as
$$
1-G(x)=\sum\limits_{j=x+1}^\infty\frac{(j+1)f_{j+1}}{\mu}
=\frac{1}{\mu}\left[(x+2)\left[1-F(x+1)\right]+\sumu\limits_{j=x+2}^{\infty}
\left[1-F(j)\right]\right].
$$
Then we use \cite[Theorem 1, p.~281]{fellerb}, together with the
fact that $1-F(x)$ is regularly varying with exponent $1-\tau\ne1$
to deduce that there exists a constant $c=c_{\tau}>0$ such that
$$
\sumu\limits_{j=x+2}^{\infty}\left[1-F(j)\right] \le
c_{\tau}(x+2)\left[1-F(x+2)\right].
$$
Hence, if $F$ satisfies Assumption \ref{ass-gamma}, then
$$
\begin{array}{rl}
1-G(x)&\ge\frac{1}{\mu}(x+2)\left[1-F(x+1)\right]
\ge x^{2-\tau-K_{\tau}\gamma(x)},\\[5pt]
1-G(x)&\le\frac{1}{\mu}(c+1)(x+2)\left[1-F(x+1)\right] \le
x^{2-\tau+K_{\tau}\gamma(x)},
\end{array}
$$
for some $K_{\tau}>0$ and large enough $x$. \qed \vskip.5truecm

\begin{remark}
It follows from Assumption~\ref{ass-gamma} and Lemma~\ref{lem-G-x}
that for each $\vep>0$ and sufficiently large $x$,
    \eq
    \begin{array}{rlll}
    x^{1-\tau-\vep}&\leq 1-F(x)&\leq x^{1-\tau+\vep},&\qquad(a)\\[5pt]
    x^{2-\tau-\vep}&\leq 1-G(x)&\leq x^{2-\tau+\vep}.&\qquad(b)
    \end{array}
    \label{1-F,1-G bound}
    \en
We will often use (\ref{1-F,1-G bound}) with $\vep$ replaced by
$\vep^6$. \qed
\end{remark}
Let us define for $\vep>0$,
\begin{equation}
\label{gerard-alpha} \alpha=\frac{1-\vep^5}{\tau-1},\qquad
h=\vep^6,
\end{equation}
and  the auxiliary event $\FN$ by
    \eq
    \label{Fdef}
    \FN =\{\forall 1\le x\leq N^{\alpha}:
    |G(x)-G^{\smallsup{N}}(x)|\leq N^{-h} [1-G(x)]\},
    \en
    where $G^{\smallsup{N}}$ is the (random) distribution function of
    $\{g_j^{\smallsup{N}}\}$, defined in (\ref{convergenceoffspring}).

\begin{lemma}
\label{lem-Fbd}For $\vep$ small enough, and $N$ sufficiently
large,
    \eq
    \prob(\FN^c) \leq N^{-h}.
    \label{Fbd}
    \en
\end{lemma}
\proof First, we rewrite $1-G^{\smallsup{N}}(x)$, for $x\in
\mathbb{N}_0$, in the following way:
    \eqalign
    1-G^{\smallsup{N}}(x)&=\sum_{n=x+1}^{\infty} g^{\smallsup{N}}_n
    =\frac{1}{L_{\sN}}\sum_{j=1}^N \sum_{n=x+1}^{\infty} D_j{\bf 1}_{\{D_j=n+1\}}
    =\frac{1}{L_{\sN}}\sum_{j=1}^N D_j{\bf 1}_{\{D_j\geq x+2\}}
    \nn\\
    &=\frac{1}{L_{\sN}}\sum_{j=1}^N \sum_{l=1}^{D_j}{\bf 1}_{\{D_j \geq x+2\}}=\frac{1}{L_{\sN}}\sum_{l=1}^{\infty}\sum_{j=1}^N{\bf 1}_{\{D_j \geq (x+2)\vee l\}},
    \enalign
where $x\vee l$ is the maximum of $x$ and $l$. Writing
    \eq
    \label{Bydef}
    B_y^{\smallsup{N}}=\sum_{j=1}^{N}{\bf 1}_{\{D_j \geq y\}},
    \en
we thus end up with
    \eq
    1-G^{\smallsup{N}}(x)=\frac{1}{L_{\sN}}\sum_{l=1}^{\infty} B_{(x+2)\vee l}^{\smallsup{N}}.
    \en
We have a similar expression for $1-G(x)$ that reads
    \eq
    1-G(x)=\frac{1}{\mu} \sum_{l=1}^{\infty} \prob(D_1\geq (x+2)\vee l).
    \en
Therefore, with
$$
\beta=\frac{1-h}{\tau-1},\qquad\mbox{ and }
\qquad\chi=\frac{1+2h}{\tau-1},
$$
we can write
    \eq
    \begin{array}{rl}
[G(x)-G^{\smallsup{N}}(x)]
    &=\left(\frac{N\mu}{L_{\sN}}-1\right)[1-G(x)]\\[5pt]
    &\quad+\frac{1}{L_{\sN}}\sum_{l=1}^{N^{\beta}}
    \left[B_{(x+2)\vee l}^{\smallsup{N}}-N\prob\big(D_1\geq (x+2)\vee l\big)\right]\\[5pt]
    &\quad+\frac{1}{L_{\sN}}\sum_{l=N^{\beta}+1}^{N^{\chi}}
    \big[B_{(x+2)\vee l}^{\smallsup{N}}-N\prob\big(D_1\geq (x+2)\vee l\big)\big]\\[5pt]
    &\quad+\frac{1}{L_{\sN}}\sum_{l=N^{\chi}+1}^{\infty}
    \big[B_{(x+2)\vee l}^{\smallsup{N}}-N\prob\big(D_1\geq (x+2)\vee l\big)\big].
    \end{array}
    \en
Hence, for large enough $N$ and $x\le
N^{\alpha}<N^\beta<N^{\chi}$, we can bound
    \eq
    \label{(1-G)diff}
    \begin{array}{rll}
R_{\sN}(x)\equiv\Big|G(x)-G^{\smallsup{N}}(x)\Big|
    &\le\left|\frac{N\mu}{L_{\sN}}-1\right|[1-G(x)]
    &\qquad(a)\\[5pt]
    &\quad+\frac{1}{L_{\sN}}\sum_{l=1}^{N^{\beta}}
    \left|B_{(x+2)\vee l}^{\smallsup{N}}-N\prob\big(D_1\geq (x+2)\vee l\big)\right|
    &\qquad(b)\\[5pt]
    &\quad+\frac{1}{L_{\sN}}\sum_{l=N^{\beta}+1}^{N^{\chi}}
    \left|B_{l}^{\smallsup{N}}-N\prob\big(D_1\geq l\big)\right|
    &\qquad(c)\\[5pt]
    &\quad+\frac{1}{L_{\sN}} \sum_{l=N^{\chi}+1}^{\infty}B_{l}^{\smallsup{N}}
    &\qquad(d)\\[5pt]
    &\quad+\frac{1}{L_{\sN}}\sum_{l=N^{\chi}+1}^{\infty}N \prob\big(D_1\geq l\big).
    &\qquad(e)
    \end{array}
    \en
We use (\ref{1-F,1-G bound}(b)) to conclude that, in
order to prove $\prob(\FN^c)\le N^{-h}$, it suffices to show
that
    \eq
    \prob\left(\bigcup\limits_{1\le x\le N^{\alpha}}
    \left\{|R_{\sN}(x)| > C_gN^{-h} x^{2-\tau-h}\right\}
    \right)\leq N^{-h},
    \en
for large enough $N$, and for some $C_g$, depending on distribution function $G$.
We will define an auxiliary
event $A_{\sN,\vep}$, such that $|R_{\sN}(x)|$ is more easy to
bound on $A_{\sN,\vep}$ and such that $\prob(A_{\sN,\vep}^c)$
is sufficiently small. Indeed, we define, with $A=3(\beta+2h)$,
    \eq
    \begin{array}{rll}
A_{\sN,\vep}(a)&=\left\{|\frac{N\mu}{L_N}-1|\leq
     N^{-3h}\right\},&\qquad(a)\\[10pt]
A_{\sN,\vep}(b)&=\left\{\max_{1\le j \le N} D_j \leq N^{\chi}\right\},&\qquad(b)\\[10pt]
A_{\sN,\vep}(c)&=\bigcap\limits_{1\le x\le N^{\beta}}
        \left\{|B_{x}^{\smallsup{N}}-N\prob(D_1\geq x)|
    \leq \sqrt{A(\log{N}) N\prob(D_1\geq x)}\right\},&\qquad(c)
    \label{F4def}
    \end{array}
    \en
and
$$
A_{\sN,\vep}=A_{\sN,\vep}(a)\cap A_{\sN,\vep}(b)\cap A_{\sN,\vep}(c).
$$

By intersecting with $A_{\sN,\vep}$ and its complement, we have
    \eq
    \begin{array}{l}
    \prob\big(\bigcup\limits_{1\le x\le N^{\alpha}}
     \{|R_{\sN}(x)| > C_gN^{-h}
    x^{2-\tau-h}\}\big)\\[10pt]
    \qquad\qquad\qquad\qquad
    \le \prob\big(A_{\sN,\vep}\cap \left\{\bigcup\limits_{1\le x\le N^{\alpha}} \{|R_{\sN}(x)| > C_gN^{-h}
    x^{2-\tau-h}\}\right\}\big)+\prob(A_{\sN,\vep}^c).
    \end{array}
    \en
We will prove that $\prob(A_{\sN,\vep}^c)\le N^{-h}$, and that on
the event $A_{\sN,\vep}$, and for each $1\le x\le N^{\alpha}$, the
right hand side of (\ref{(1-G)diff}) can be bounded by
$C_gN^{-h}x^{2-\tau-h}$. We start with the latter statement.

Consider the right hand side of (\ref{(1-G)diff}). Clearly, on
$A_{\sN,\vep}(a)$, the first term of $|R_{\sN}(x)|$ is bounded by
$N^{-3h}[1-G(x)]\le C_gN^{-3h} x^{2-\tau+h}\le C_gN^{-h}
x^{2-\tau-h} $, where the one but last inequality follows from
(\ref{1-F,1-G bound}(b)), and the  last since $x\le N^{\alpha}< N$
so that $x^{2h}< N^{2h} $. Since for $l>N^{\chi}$ and each
$j,\, 1\le j\le N$,   we have $D_j>l$ is the empty set on
$A_{\sN,\vep}(b)$, the one but last term of $|R_{\sN}(x)|$ vanish
on $A_{\sN,\vep}(b)$. The last term of $|R_{\sN}(x)|$ can for $N$ large be bounded, using the inequality $L_{\sN}\ge N$
and (\ref{1-F,1-G bound}(a)),
\begin{equation}
\frac1{L_{\sN}}\sum_{l=N^{\chi}+1}^\infty N\prob(D_1\ge l)
\le \sum_{l=N^{\chi}+1}^\infty l^{1-\tau+h}\nn
\quad\le
\frac{N^{\chi(2-\tau+h)}}{\tau-2}
<C_g N^{-h+\alpha(2-\tau+h)} \le C_g N^{-h} x^{2-\tau-h},
\end{equation}
for all $x\le N^{\alpha}$, and where we also used that for $\vep$ sufficiently small and $\tau>2$,
$$
\chi(2-\tau+h)<-h+\alpha(2-\tau+h).
$$

We bound the
third term of $|R_{\sN}(x)|$ as
    \eqalign
    \frac{1}{L_{\sN}}\sum_{l=N^{\beta}+1}^{N^{\chi}}
    \big|B_{l}^{\smallsup{N}}-N\prob(D_1\geq l)\big|
    &\leq \frac{1}{N}\sum_{l=N^{\beta}+1}^{N^{\chi}}
    [B_{l}^{\smallsup{N}}+N\prob(D_1\geq l)]\nn\\
    &\leq N^{\chi} [N^{-1}B_{N^{\beta}}^{\smallsup{N}}+\prob(D_1\geq N^{\beta})].
    \label{twoterms}
    \enalign
We note that due to (\ref{1-F,1-G bound}(a)),
    \eq
    \prob(D_1\geq N^{\beta}) \geq N^{\beta(1-\tau-h)},
    \en
for large enough $N$, so that
    \eq
    a_N=\sqrt{A(\log{N})N\prob(D_1\geq N^{\beta})}\leq
    N\prob(D_1\geq N^{\beta}).
    \en
Therefore, on $A_{\sN,\vep}(c)$, we obtain that
    \eq
    B_{N^{\beta}}^{\smallsup{N}} \leq 2N\prob(D_1\geq N^{\beta}),
    \en
for $\vep$ small enough and large enough $N$. Furthermore as $\vep
\downarrow 0$,
\[
N^{\chi+\beta(1-\tau+h)}<C_g N^{-h+\alpha(2-\tau-h)}\le
C_g N^{-h}x^{2-\tau-h},
\]
for $x\le N^{\alpha}$, $2-\tau-h<0$, because (after multiplying by $\tau-1$ and dividing by
$\vep^5$)
\[
\chi+\beta(1-\tau+h)<-h+\alpha(2-\tau-h), \qquad \mbox{or} \qquad
\vep(2+2\tau-h)<\tau-2+h,
\]
as $\vep$ is sufficiently small. Thus, the third term of $|R_{\sN}(x)|$ satisfies the required bound.

We bound the second term of $|R_{\sN}(x)|$ on
$A_{\sN,\vep}(a)\cap A_{\sN,\vep}(c)$ by
\eq
\frac{1}{N}\sum_{l=1}^{N^{\beta}}\sqrt{A(\log{N})
N\prob\big(D_1\geq (x+2)\vee l\big)}
=\frac{\sqrt{A\log{N}}}{\sqrt{N}} \sum_{l=1}^{N^{\beta}}
    \sqrt{\prob\big(D_1\geq (x+2)\vee l\big)}.
    \label{rempf1}
\en
 Let $c$ be a constant such that
$(\prob(D_1>x))^{\frac12}\le c x^{(1-\tau+h)/2}$,
then for all $1\le x\leq
N^{\alpha}$,
\begin{eqnarray}
&&\frac{1}{L_{\sN}}\sum_{l=1}^{N^{\beta}}
    \big|B_{(x+2)\vee l}^{\smallsup{N}}-N\prob\big(D_1\geq (x+2)\vee l\big)\big|
    \leq  \frac{c\sqrt{A\log{N}}}{\sqrt{N}} \sum_{l=1}^{N^{\beta}}
    \big((x+2)\vee l\big)^{(1-\tau+h)/2}\nn\\
&&\qquad\qquad\leq \frac{c\sqrt{A\log{N}}}{\sqrt{N}}
    \big[x^{(3-\tau+h)/2}+N^{\beta(3-\tau+h)/2}\big]
    \leq \frac{2c\sqrt{A\log{N}}}{\sqrt{N}} N^{\beta(3-\tau+h)/2}\nn\\
&&\qquad\qquad\leq N^{h-1/2} N^{\beta(3-\tau+h)/2}
    <C_gN^{-h}N^{\alpha(2-\tau-h)}\le C_gN^{-h}x^{2-\tau-h},
\end{eqnarray}
because
\[
h-1/2+\beta(3-\tau+h)/2<-h+\alpha(2-\tau-h),\qquad \mbox{or} \qquad
h(5\tau-4-h) <2\eps^5(\tau-2+h),
\]
for $\vep$ small enough and $\tau \in (2,3)$. We have shown that
for $1\le x\leq N^{\alpha}$, $N$ sufficiently large, and on the
event $A_{\sN,\vep}$,
    \eq
    |R_{\sN}(x)| \leq C_gN^{-h} x^{2-\tau-h}.
    \en

It remains to prove that $\prob(A_{\sN,\vep}^c)\le N^{-h}$. We use
that
    \eq
    \prob(A_{\sN,\vep}^c)\leq \prob(A_{\sN,\vep}(a)^c)
    +\prob(A_{\sN,\vep}(b)^c)+\prob(A_{\sN,\vep}(c)^c),
    \en
and we bound each of the three terms separately.

The bound
    \eq
    \prob(A_{\sN,\vep}(a)^c)=
    \prob\left(\left|\frac1N \sum_{j=1}^N (D_j-\mu)\right|>
    N^{-3h}\cdot L_{\sN}/N\right)
    \leq \frac 13 N^{-h},
    \en
follows, since $N^{-\frac1{\tau-1}}\sum_{j=1}^N (D_j-\mu) $,
converges to a stable law, for $2<\tau<3$.

The bound on~$\prob(A_{\sN,\vep}(b)^c)$ is a trivial estimate
using (\ref{1-F,1-G bound}(a)). Indeed, for $N$ large,
    \eq
    \prob(A_{\sN,\vep}(b)^c)=\prob\Big(\max_{1\le j \le N} D_j > N^{\chi}\Big)
    \le N\prob(D_1\ge N^{\chi})
    \le N^{\chi(1-\tau+h)+1}
    \le \frac 13 N^{-h},
    \en
for small enough $\vep$, because $\tau>2+h$. For the third term
$\prob(A_{\sN,\vep}(c)^c)$, we will use a bound given by Janson
\cite{Jans02}, which states that for a binomial random variable
$X$ with parameters $N$ and $p$, and all $t>0$,
 \eq
 \label{binbd}
\prob(|X-Np|\geq t) \leq
2\exp\left\{-\frac{t^2}{2(Np+t/3)}\right\}. \en We will apply
(\ref{binbd}) with $t=a_{\sN}(x)=\sqrt{A(\log{N}) N\prob(D_1\geq x)}$,
and obtain that uniformly in $x\leq N^{\alpha}$,
    \begin{eqnarray}
    \label{Bbd}
    && \prob\left(|B_{x}^{\smallsup{N}}-N\prob(D_1\geq x)|>  a_{\sN}(x)\right)
    \leq 2\exp\left\{-\frac{a_{\sN}(x)^2}{2(N\prob(D_1\geq x)+a_{\sN}(x)/3)}\right\}\nn\\
    &&\leq 2\exp\left\{-\frac{A\log N}{2(1+\frac13\sqrt{A\log N/(N\prob(D_1\geq N^{\alpha}))})}\right\}
    \leq 2N^{-A/3},
    \end{eqnarray}
because
$$
\frac{\log N}{N\prob(D_1\geq N^{\alpha})}\le
\frac{\log N}{N^{1+\alpha(\tau-1-h)}} \to 0,
$$
as $N\to\infty$.
Thus, (\ref{Bbd})gives us, using $A=3(\beta+2h)$,
    \eqalign
    \label{F4cbd}
    \prob(A_{\sN,\vep}(c)^c)
    &\leq \sum_{x=1}^{N^{\beta}} \prob\big(|B_{x}^{\smallsup{N}}-N\prob(D_1\geq x)|
    > a_{\sN}(x)\big)
    \leq 2N^{\beta-A/3}=2N^{-2h}\leq \frac13 N^{-h}.
    \enalign
This completes the proof of the lemma. \qed \vskip0.5cm

For the third lemma we introduce some further notation. For any
$x\in \NN$, define
$$
\hat S_x^{\smallsup{N}}=\sum_{i=1}^x \hat
X_i^{\smallsup{N}},\qquad \hat
V_x^{\smallsup{N}}=\max\limits_{1\le i\le x} \hat
X_i^{\smallsup{N}},
$$
where $\{\hat X_i^{\smallsup{N}}\}_{i=1}^x$ have the same law, say
$\hat H^{\smallsup{N}}$, but are {\it not necessarily} independent.

\begin{lemma}[Sums with law $\hat H^{\smallsup{N}}$ on the good event]~\\
\label{lem-sums}

$({\rm i})$ If $\hat H^{\smallsup{N}}$ satisfies
\eq
\label{lem-sums-cond1} [1-\hat H^{\smallsup{N}}(z)]\leq
[1+2N^{-h}] [1-G(z)],\quad \forall \, z\le y,
\en
then for all $x\ge 1$, there exists  a constant $b'$, such that:
    \eq
    \prob\Big(\hat S_x^{\smallsup{N}}\geq y\Big)
    \leq b' x [1+2N^{-h}]\big[1-G\big(y)\big];
    \en

${\rm (ii)}$ If $\hat H^{\smallsup{N}}$ satisfies
\eq
\label{lem-sums-cond2} [1-\hat H^{\smallsup{N}}(y)]\geq
[1-2N^{-h}] [1-G(y)] ,
\en
and $\{\hat
X_i^{\smallsup{N}}\}_{i=1}^x$ are independent, then for all
$x\ge 1$,
\eq
\label{gerard21}
    \prob\left(\hat V_x^{\smallsup{N}}\leq  y\right)
    \leq \left(1-[1-2N^{-h}][1-G\left(y\right)]\right)^x.
\en
\end{lemma}

\proof We first bound $\prob\Big(\hat S_x^{\smallsup{N}} \geq
y\Big)$. We write
    \eq
    \label{probsplit}
    \prob\Big(\hat S_x^{\smallsup{N}}\geq y\Big)
    \leq
    \prob\Big(\hat S_x^{\smallsup{N}} \geq
    y,\hat V_x^{\smallsup{N}}
    \leq y\Big)
    +\prob\big(\hat V_x^{\smallsup{N}}> y\big).
    \en
The second term is bounded due to~(\ref{lem-sums-cond1}) by
    \eq
    \label{afschat-sup}
    x\prob\big(\hat X_1^{\smallsup{N}}> y\big)
    = x\Big[1-{\hat H}^\smallsup{N}\big(y\big)\Big]\\[10pt]
    \leq x[1+2N^{-h}]\Big[1-G\big(y\big)\Big].
    \en
We use the Markov inequality and~(\ref{lem-sums-cond1}) to bound
the first term on the right-hand side of (\ref{probsplit}) by
    \begin{eqnarray}
    &&\prob\Big(\hat S_x^{\smallsup{N}} \geq
    y,\hat V_x^{\smallsup{N}}\leq y\Big)\leq \frac1{y}\expec\left(\hat S_x^{\smallsup{N}}
    {\bf 1}_{\{\hat V_x^{\smallsup{N}}\leq y\}}\right)
    \leq \frac{x}{y} \expec\left(\hat X_1^{\smallsup{N}}
    {\bf 1}_{\{\hat X_1^{\smallsup{N}}\leq y\}}\right)\nn\\
    &&\qquad \leq  \frac{x}{y} \sum_{i=1}^{y}
    [1-\hat H^{\smallsup{N}}(i)] \leq  \frac{x}{y}[1+2N^{-h}]
    \sum_{i=1}^{y} [1-G(i)].
    \label{upperbd}
    \end{eqnarray}
For the latter sum, we use \cite[Theorem 1(b), p.~281]{fellerb},
together with the fact that $1-G(y)$ is regularly varying with
exponent $2-\tau\ne 1,$ to deduce that there exists a constant
$c_1$ such that
    \eq
    \sum_{i=1}^{y} [1-G(i)]\leq c_1 y[1-G(y)].
    \label{sum(1-G)}
    \en

Combining (\ref{probsplit}), (\ref{afschat-sup}), (\ref{upperbd}) and (\ref{sum(1-G)}), we conclude that
    \eq
    \prob\Big(\hat S_x^{\smallsup{N}} \geq
    y\Big) \leq
  b' x [1+2N^{-h}]
    \big[1-G\big(y)\big)\big],
    \en
where $b'=c_1+1$.
This completes the proof of Lemma \ref{lem-sums}(i).

For the proof of (ii), we use independence of $\{\hat
X_i^{\smallsup{N}}\}_{i=1}^x$, and condition~(\ref{lem-sums-cond2}), to conclude that
\[
\prob\Big(\hat V_x^{\smallsup{N}}\leq y\Big)
    =\Big(\hat H^\smallsup{N}\left(y\right)\Big)^x
       =\left(1-\Big[1-\hat H^\smallsup{N}(y)\Big]\right)^x
 \leq \left(1-[1-2N^{-h}][1-G\left(y\right)]\right)^x.
\]
Hence,  (\ref{gerard21}) holds. \qed

\begin{remark}
\label{rem-conf}In the proofs in the appendix, we will only use
that
\begin{itemize}
\item[{\rm (i)}] the event $\FN$ holds {\bf whp}; \item[{\rm
(ii)}] that $L_{\sN}$ is concentrated around its mean; \item[{\rm
(iii)}] that the maximal degree is bounded by $N^{\chi}$ for any
$\chi>1/(\tau-1)$, with {\bf whp}.
\end{itemize}
Moreover, the proof of Proposition \ref{prop-caft} relies on
\cite[Proposition A.3.1]{HHV03}, and in its proof it was further
used that
\begin{itemize}
\item[{\rm (iv)}] $p_{\sN}\leq N^{\alpha_2}$, {\bf whp}, for any
$\alpha_2>0$, where $p_{\sN}$ is the total variation distance
between $g$ and $g^{\smallsup{N}}$, i.e.,
    \eq
    p_{\sN}= \frac 12 \sum_n |g_n-g^{\smallsup{N}}_n|.
    \en
\end{itemize}
Therefore, if instead of taking the degrees i.i.d.\ with
distribution $F$, we would take the degrees in an exchangeable way
such that the above restrictions hold, then the proof carries on
verbatim. In particular, this implies that our results also hold
for the usual configuration model, where the degrees are fixed, as
long as the above restrictions are satisfied.\qed
\end{remark}


\subsection{Proof of Proposition \ref{prop-caft}}
\label{sec-caft} The proof makes use of \cite[Proposition
A.3.1]{HHV03}, which proves the statement in Proposition
\ref{prop-caft} under an additional condition.

In order to state
this condition, let, for $i=1,2$,
$\{\hat{Z}_j^{\smallsup{i,N}}\}_{j\geq 1}$ be two independent
copies of the delayed BP, where $\hat{Z}_1^{\smallsup{i,N}}$ has
law $\{f_n\}$ given in (\ref{kansen}), and where the offspring of any
individual in generation $j$ with $j > 1$ has law
$\{g_n^{\smallsup{N}}\}$, where $g_n^{\smallsup{N}}$ is defined in
(\ref{gnN}). Then, the conclusion of Proposition \ref{prop-caft}
follows from \cite[Proposition A.3.1]{HHV03}, for any $m$ such
that, for any $\eta>0$, and $i=1,2$,
    \begin{equation}
    \label{asstau}
    \prob(\sum_{j=1}^m \hat Z_j^{\smallsup{i,N}} \geq N^{\eta})=o(1).
    \end{equation}
By exchangeability it suffices to prove (\ref{asstau}) for $i=1$
only, we can therefore simplify notation and write further
${\hat Z}_{k}^{\smallsup{N}}$ instead of ${\hat Z}_{k}^{\smallsup{i,N}}$.
We turn to the proof of (\ref{asstau}).

By Lemma \ref{lem-Fbd} and (\ref{1-F,1-G bound}(b)), respectively, for every
$\eta>0$, there exists a $c_{\eta}>0$, such that {\bf whp} for all
$x\leq N^{\alpha}$,
    \eq
    \label{1-GNbd}
    1-G^{\smallsup{N}}(x)\leq [1+2N^{-h}][1-G(x)]\le c_{\eta} x^{2-\tau+\eta}.
    \en

We call a generation $j\geq 1$ {\bf good}, when
    \eq
    \hat Z_j^{\smallsup{N}} \leq \Big(\hat Z_{j-1}^{\smallsup{N}} \log{N}\Big)^{\frac{1}{\tau-2-\eta}},
    \en
and {\bf bad} otherwise, where as always $\hat
Z_{0}^{\smallsup{N}}=1$. We further write
    \eq
    H_m=\{\text{generations }1, \ldots, m \text{ are good}\}.
    \en
We will prove that when $H_m$ holds, then $\sum_{j=1}^m \hat
Z_j^{\smallsup{N}} \leq N^{\eta}$. Indeed, when generations $1,
\ldots, m$ are all good, then, for all $j\leq m$,
    \eq
    \label{ZjHbd}
    \hat Z_j^{\smallsup{N}}\leq (\log{N})^{\sum_{i=1}^{j} (\tau-2-\eta)^{-i}}.
    \en
Therefore,
    \eq
    \label{gerard25}
    \sum_{j=1}^m \hat Z_j^{\smallsup{N}} \leq
    m(\log N)^{\sum_{i=1}^{m} (\tau-2-\eta)^{-i}}
    \leq m (\log{N})^{\frac{(\tau-2-\eta)^{-{m-2}}}{(\tau-2-\eta)^{-1}-1}}
    \leq N^{\eta},
    \en
for any $\eta>0$, when $N$ is sufficiently large. We conclude that
    \begin{equation}
    \label{asstaured}
    \prob(\sum_{j=1}^m \hat Z_j^{\smallsup{N}} > N^{\eta})\leq \prob(H_m^c),
    \end{equation}
and Proposition \ref{prop-caft} follows if we show that
$\prob(H_m^c)=o(1)$. In order to do so, we write
    \eq
    \prob(H_m^c) = \prob(H_1^c)+\sum_{j=1}^{m-1} \prob(H_{j+1}^c\cap H_{j}).
    \en
For the first term, we use (\ref{1-F,1-G bound}(a)) to deduce that
    \eq
    \prob(H_1^c)= \prob(D_1>(\log{N})^{\frac{1}{\tau-2-\eta}})
    \leq (\log{N})^{-\frac{\tau-1-\eta}{\tau-2-\eta}}
    \leq (\log{N})^{-1}.
    \en
For $1\le j \le m$, we have
$\hat Z_{j}^{\smallsup{N}}\le \sum_{k=1}^m \hat Z_{k}^{\smallsup{N}}$, and
using (\ref{gerard25}),
\eq
    \sum_{j=1}^m \hat Z_j^{\smallsup{N}} \leq
    m(\log{N})^{\frac{(\tau-2-\eta)^{-{m-2}}}{(\tau-2-\eta)^{-1}-1}}
    =K_{\sN}.
\en
Using Lemma \ref{lem-sums}(i) with ${\hat H}^{\smallsup{N}}=G^{\smallsup{N}}$, $x=l$ and
$y=v_{\sN}=(l \log N)^{\frac1{\tau-2-\eta}}$, where (\ref{lem-sums-cond1})  follows from (\ref{1-GNbd}),
we obtain that
    \begin{eqnarray}
    &&\prob(H_{j+1}^c\cap H_{j})
    \leq \sum_{l=1}^{K_{\sN}}
    \prob\Big(\hat Z_{j+1}^{\smallsup{N}}
    \geq v_{\sN}\Big|\hat Z_{j}^{\smallsup{N}}=l\Big)
    \prob(\hat Z_{j}^{\smallsup{N}}=l)\leq \max_{1\le l \le K_{\sN}}
    \prob\Big(\hat S_l^{\smallsup{N}}
    \geq v_{\sN}\Big)\nn\\
    &&\qquad\leq  \max_{1\le l \le x_o}
    \prob\Big(\hat S_l^{\smallsup{N}}
    \geq v_{\sN}\Big)
    +b' \max_{x_o\le l \le K_{\sN}}
    l [1+2N^{-h}]\big[1-G( v_{\sN})\big].
    \end{eqnarray}
The case where $l\leq x_o$ can, by (\ref{1-GNbd}),  be bounded
as
 \begin{eqnarray}
 &&x_o \prob\Big(\sum_1^l \hat X_{j}^{\smallsup{N}}
 \geq v_{\sN})
\leq
x_o \prob\Big(\bigcup_{j=1}^l \{\hat X_{j}^{\smallsup{N}}
 \geq v_{\sN}/l\}\Big)
 \leq
x_o \sum_{j=1}^l
\prob\Big(\hat X_{j}^{\smallsup{N}}
 \geq v_{\sN}/l\Big)\nn\\
 &&\qquad
 \le
 c_{\eta} l x_o(v_{\sN}/l)^{2-\tau+\eta}
 =\frac{c_{\eta}x_o}{l^{2-\tau+\eta}}(\log N)^{-1}\le C (\log N)^{-1}.
  \end{eqnarray}
Furthermore by (\ref{1-F,1-G bound}(b)),
\eq
\max_{x_o\le l \le K_{\sN}}
    l [1-G(v_{\sN})]
    \le \max_{x_o\le l \le K_{\sN}} l v_{\sN}^{1-\tau+\eta}< (\log N)^{-1}.
\en
This completes the proof of Proposition  \ref{prop-caft}. \qed

\subsection{Some further preparations}

\label{sec-prep}

Before we can prove Propositions~\ref{prop-weakconv2bsec}
and~\ref{prop-Tub}    , we state a lemma that was proved in
\cite{HHV03}.

We introduce some notation. Suppose we have $L$ objects divided
into $N$ groups of sizes $d_1, \ldots, d_{\sN}$, so that
$L=\sum_{i=1}^N d_i$. Suppose we draw an object at random. This
gives a distribution $g^{\smallsup{\vec{d}}}$, i.e.,
    \eq
    g^{\smallsup{\vec{d}}}_n = \frac{1}{L} \sum_{i=1}^N d_i{\bf 1}_{\{d_i=n+1\}},
    \quad n=0,1,\ldots
    \en
Clearly, $g^{\smallsup{N}}=g^{\smallsup{\vec{D}}}$, where
$\vec{D}=(D_1, \ldots, D_{\sN})$. We further write
    \eq
    G^{\smallsup{\vec{d}}}(x) =\sum_{n=0}^x g^{\smallsup{\vec{d}}}_n.
    \en

We next label $M$ of the $L$ objects, and suppose that the
distribution $G^{\smallsup{\vec{d}}}_{\sM}(x)$ is obtained in a
similar way from drawing conditionally on drawing an unlabelled
object. More precisely, we remove the labelled objects from all
objects thus creating new $d_1', \ldots, d_{\sN}'$, and we let
$G^{\smallsup{\vec{d}}}_{\sM}(x)=G^{\smallsup{\vec{d}'}}(x).$ Even
though this is not indicated, the law
$G^{\smallsup{\vec{d}}}_{\sM}$ depends on what objects have been
labelled.

Lemma \ref{lem-stochbds} below shows that the law
$G^{\smallsup{\vec{d}}}_{\sM}$ can be stochastically bounded above
and below by two specific ways of labeling objects. Before we can
state the lemma, we need to describe those specific labellings.

For a vector $\vec{d}$, we denote by  $d_{\smallsup{1}}\leq
d_{\smallsup{2}}\le \ldots\le d_{\smallsup{N}}$ the ordered
coordinates. Then the laws $\overline{G}^{\smallsup{\vec{d}}}_{\sM}$
and $\underline{G}^{\smallsup{\vec{d}}}_{\sM}$, respectively, are
defined by successively decreasing $d_{\smallsup{N}}$ and
$d_{\smallsup{1}}$, respectively, by one. Thus,
    \begin{eqnarray}
    \overline{G}^{\smallsup{\vec{d}}}_{\sss 1}(x)
    &=&\frac{1}{L-1} \sum_{i=1}^{N-1}
    d_{\smallsup{i}}{\bf 1}_{\{d_{\smallsup{i}}\leq x+1\}} +
    \frac{d_{\smallsup{N}}-1}{L-1}{\bf 1}_{\{d_{\smallsup{N}}-1\leq x+1\}}\label{f1},\\
    \underline{G}^{\smallsup{\vec{d}}}_{\sss 1}(x)&=&\frac{1}{L-1} \sum_{i=2}^{N}
    d_{\smallsup{i}}{\bf 1}_{\{d_{\smallsup{i}}\leq x+1\}} +
    \frac{d_{\smallsup{1}}-1}{L-1}{\bf 1}_{\{d_{\smallsup{1}}-1\leq x+1\}}.\label{h1}
    \end{eqnarray}
For $\overline{G}^{\smallsup{\vec{d}}}_{\sM}$ and
$\underline{G}^{\smallsup{\vec{d}}}_{\sM}$, respectively, we
perform the above change $M$ times, and after each repetition we
reorder the groups. Here we note that when $d_{\smallsup{N}}=1$
(in which case $d_i=1$, for all $i$), and for
$\overline{G}^{\smallsup{\vec{d}}}_{\sss 1}$ we decrease
$d_{\smallsup{N}}$ by one, that we only keep the groups with $d_i=
1$. The same rule applies when $d_{\smallsup{1}}=1$ and for
$\underline{G}^{\smallsup{\vec{d}}}_{\sss 1}$ we decrease
$d_{\smallsup{1}}$ by one. Thus, in these cases, the number of
groups of objects, indicated by $N$, is decreased by 1. Applying
the above procedure to $\vec{d}=(D_1,\dots,D_{\sN})$ we obtain
that, for all $x\geq 1$,
    \eqalign
    \overline{G}^{\smallsup{N}}_{\sM}(x)
    \equiv \overline{G}^{\smallsup{\vec{D}}}_{\sM}(x)
    &\le\frac{1}{L_{\sN}-M}\sum_{i=1}^N D_i {\bf 1}_{\{D_i\leq x+1\}}\label{barG}
    =\frac{L_{\sN}}{L_{\sN}-M}G^{\smallsup{N}}(x),\\
    \underline{G}^{\smallsup{N}}_{\sM}(x)\equiv
    \underline{G}^{\smallsup{\vec{D}}}_{\sM}(x)
    &\ge\frac{1}{L_{\sN}-M}\Big[\sum_{i=1}^N D_i {\bf 1}_{\{D_i\leq x+1\}}-M\Big]
    =\frac{1}{L_{\sN}-M}\Big[L_{\sN} G^{\smallsup{N}}(x)-M\Big],
    \label{underbarG}
    \enalign
where equality is achieved precisely when $D_{\smallsup{N}}\geq
x+M$, and $\#\{i:D_i=1\}\geq M$, respectively.

Finally, for two distribution functions $F,G$, we write that $F
\preceq G$ when $F(x)\geq G(x)$ for all $x$. Similarly, we write
that $X \preceq Y$ when for the distribution functions $F_X, F_Y$
we have that $F_X \preceq F_Y$.

We next prove stochastic bounds on the distribution
$G^{\smallsup{\vec{d}}}_{\sM}(x)$ that are uniform in the choice
of the $M$ labelled objects. The proof of Lemma
\ref{lem-stochbds} can be found in \cite{HHV03}.
\begin{lemma}
\label{lem-stochbds} For all choices of $M$ labelled objects
    \begin{equation}
    \underline{G}^{\smallsup{\vec{d}}}_{\sM}\preceq G^{\smallsup{\vec{d}}}_{\sM}
    \preceq \overline{G}^{\smallsup{\vec{d}}}_{\sM}.
    \label{bdsdistr}
    \end{equation}

Moreover, when $X_1, \ldots, X_j$ are draws from
$G^{\smallsup{\vec{d}}}_{\sss M_1}, \ldots,
G^{\smallsup{\vec{d}}}_{\sss M_l}$, where the only dependence
between the $X_i$ resides in the labelled objects, then
    \eq
    \sum_{i=1}^j \underline{X}_i \preceq \sum_{i=1}^j X_i \preceq \sum_{i=1}^j \overline{X}_i,
    \label{bdsums}
    \en
where $\{\underline{X}_i\}_{i=1}^j$ and
$\{\overline{X}_i\}_{i=1}^j$, respectively, are i.i.d.\ copies of
$\underline{X}$ and $\overline{X}$ with laws
$\underline{G}^{\smallsup{N}}_{\sM}$ and
$\overline{G}^{\smallsup{N}}_{\sM}$ for $M=\max_{1\le i \le l}
M_i$, respectively.
\end{lemma}

We will apply Lemma \ref{lem-stochbds} to
$G^{\smallsup{\vec{D}}}=G^{\smallsup{N}}.$


\subsubsection{The inductive step}
Our key result, which will yield the proofs of Proposition
\ref{prop-weakconv2bsec} and \ref{prop-Tub}, is Proposition
\ref{dnziva_P1} below. This proposition will be proved in Section
\ref{sec-mainproofapp}. For its formulation we need some more
notation.

As before we simplify notation and write further
on $Z_{k}^{\smallsup{N}}$ instead of $Z_{k}^{\smallsup{i,N}}$.
Similarly, we write ${\cal Z}_{k}$ instead of ${\cal
Z}_{k}^{\smallsup{i}}$ and ${\cal T}^{\smallsup{N}}_m(\vep)$
instead of ${\cal T}^{\smallsup{i,N}}_m(\vep)$.
Recall that we have defined previously
    $$
    \kappa=\frac1{\tau-2}>1\qquad \mbox{and} \qquad  \alpha=\frac{1-\vep^5}{\tau-1}.
    $$
In the sequel we work with $Y_k^{\smallsup{N}}>\varepsilon$, for
$k$ large enough, i.e., we work with
$Z_{k}^{\smallsup{N}}>e^{\varepsilon \kappa^k}>1$, due to
definition~(\ref{Yidef2}). Hence, we can treat these definitions
as
    \eq
    Y_k^{\smallsup{N}}=\kappa^{-k} \log (Z_{k}^{\smallsup{N}}) \qquad \text{and} \qquad
    Y_k=\kappa^{-k}\log ({\cal Z}_{k}).
    \label{Yidef2_}
  \en
With $\gamma$ defined in the Assumption \ref{ass-gamma}, and
$0<\varepsilon<3-\tau$, we take $m_\vep$ sufficiently large to
have
\begin{equation}
\label{dnziva_2_} \sumu\limits_{k=m_\varepsilon}^{\infty}
(\tau-2+\varepsilon)^{k(1-\gamma)}\le \varepsilon^3 \qquad
\mbox{and}\qquad \sumu\limits_{k=m_\varepsilon}^{\infty} k^{-2}\le
\vep/2.
\end{equation}
For any $m_\vep\le m<k$, we denote
\begin{equation}
\label{dnziva_16} M_k^{\smallsup{N}}=\sum\limits_{j=1}^k
Z_{j}^{\smallsup{N}}, \quad \mbox{and} \quad  M_k=\sum\limits_{j=1}^k
{\cal Z}_{j}.
\end{equation}

As defined in Section 3 of \cite{HHV03} we speak of free stubs at
level $l$, as the free stubs connected to nodes at distance $l-1$
from the root; the total number of free stubs, obtained
immediately after pairing of all stubs of level $l-1$ equals
$Z_{l}^{\smallsup{N}}$ . For any $l\ge 1$ and $1\le x\le
Z_{l-1}^{\smallsup{N}}$, let $Z_{x,l}^{\smallsup{N}}$ denote the
number of constructed free stubs at level $l$ after pairing of the
first $x$ stubs of $Z_{l-1}^{\smallsup{N}}$. Note that for
$x=Z_{l-1}^{\smallsup{N}}$, we obtain
$Z_{x,l}^{\smallsup{N}}=Z_{l}^{\smallsup{N}}$. For general $x$,
the quantity $Z_{x,l}^{\smallsup{N}}$ is loosely speaking the sum
of the number of children of the first $x$ stubs at level $l-1$,
and according to the coupling at fixed times (Proposition
\ref{prop-caft})
 this number is for fixed $l$, {\bf whp} equal to the number of children of
 the first $x$ individuals in generation $l-1$ of the BP $\{{\cal Z}_k\}_{k\ge 1}$.

We introduce the event $ \hat F_{m,k}(\varepsilon)$,
\begin{equation}
\label{dnziva_2} \hat F_{m,k}(\varepsilon)=
\begin{array}{ll}
\{k\in{\cal T}^{\smallsup{N}}_m(\vep)\}&\qquad (a)\\
\cap\{\forall m<l\le k-1:~ |Y_{l}^{\smallsup{N}}-Y_m^{\smallsup{N}}|\leq \vep^3\}&\qquad (b)\\
\cap\{\vep\leq Y_m^{\smallsup{N}}\leq\vep^{-1}\}&\qquad (c)\\
\cap\{M_{m}^{\smallsup{N}}\leq 2Z_m^{\smallsup{N}}\}.&\qquad (d)
\end{array}
\end{equation}
In the proof of Proposition \ref{dnziva_P1} we compare the quantity $Z_{x,l}^{\smallsup{N}}$ to the sum
$\sum_{i=1}^x X_{i,l-1}^{\smallsup{N}}$ for part (a) and to $\max_{1\le i\le x} X_{i,l-1}^{\smallsup{N}}$
for part (b). We then couple  $X_{i,l-1}^{\smallsup{N}}$ to
$ {\overline X}_{i,l-1}^{\smallsup{N}}$ for part (a) and to $ {\underline X}_{i,l-1}^{\smallsup{N}}$ for part (b).
Among other things, the event $ \hat F_{m,k}(\varepsilon)$ ensures that these couplings hold.

\begin{prop}[Inductive step]
\label{dnziva_P1} Let $F$ satisfy Assumption~\ref{ass-gamma}. For
$\vep>0$ sufficiently small and $c_{\gamma}$ sufficiently large, there exist a constant
$b=b(\tau,\vep)>0$ such that, for $x=Z_{l-1}^{\smallsup{N}}\wedge
N^{\frac{(1-\vep/2)}{\kappa(\tau-1)}}$,
$$
\begin{array}{ll}
    \prob\Big(\hat F_{m,l}(\varepsilon)\cap \big\{
      Z_{x,l}^{\smallsup{N}}
      \geq (l^3x)^{\kappa+c_{\gamma}\gamma(x)}\big\}
    \Big)\leq bl^{-3},&\qquad(a)\\[15pt]
\prob\Big(\hat F_{m,l}(\varepsilon)\cap \big\{
      Z_{x,l}^{\smallsup{N}}
      \leq \left(\frac x {l^3}\right)^{\kappa-c_{\gamma}\gamma(x)}\big\}
    \Big)\leq bl^{-3}.&\qquad(b)
\end{array}
$$
\end{prop}

\label{proof(i,ii) idea} The proof of Proposition \ref{dnziva_P1}
is quite technical and is given in Section~\ref{sec-mainproofapp}.
In this section we give a short overview of the proof. For $l\ge
1$, let $SPG_{l}$ denote the shortest path graph containing all
nodes on distance $l-1$, and including all stubs at level $l$, i.e., the moment we have $Z_{l}^{\smallsup{N}}$
free stubs at level $l$. For
$i\in\{1,\dots ,x\}$, let $X_{i,l-1}^{\smallsup{N}}$ denote the
number of brother stubs of a stub attached to $i^{\rm th}$ stub of
$SPG_{l-1}$ (see Figure~\ref{dnziva_f1}).

\fig{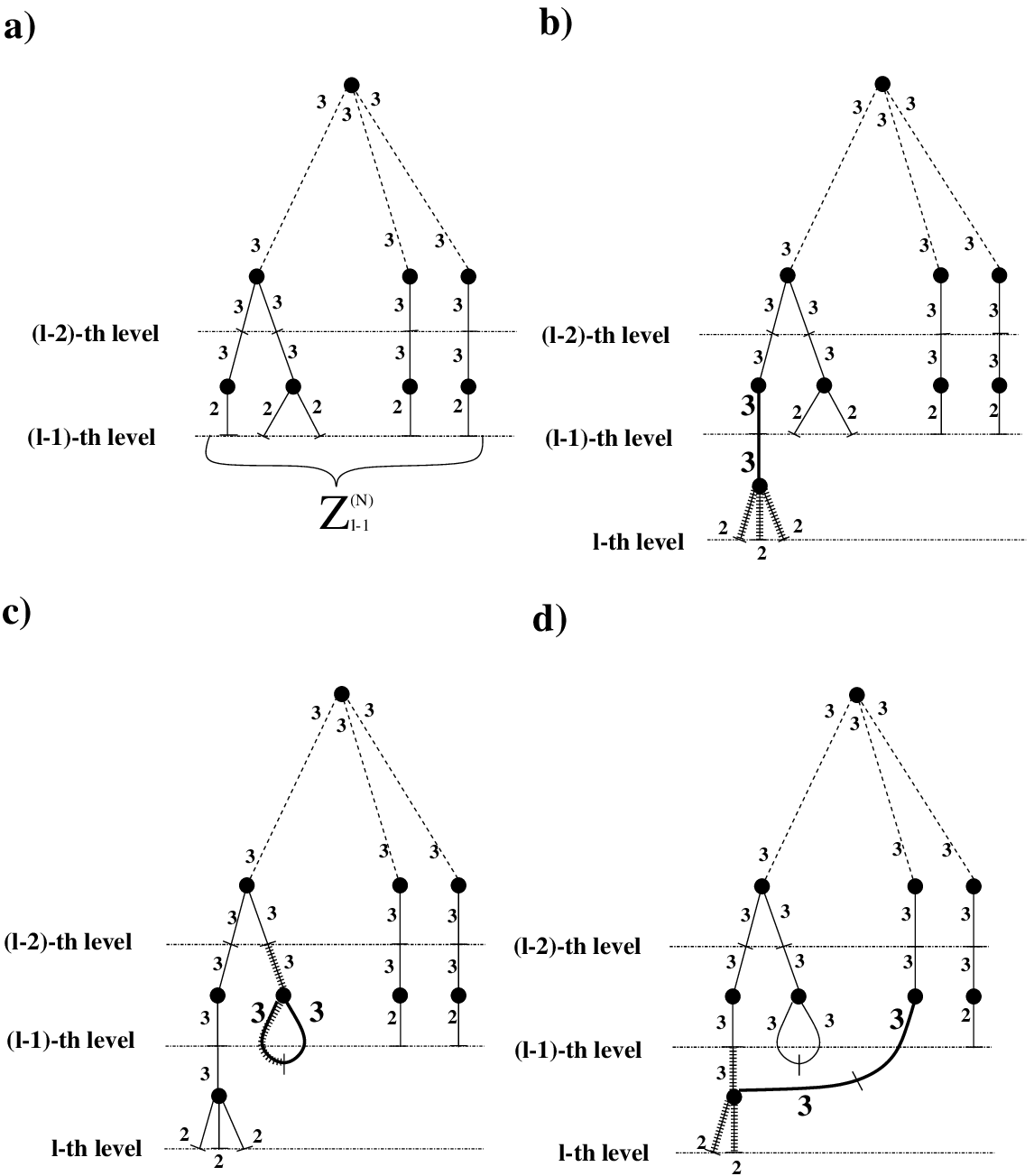}{The building of the $l^{\rm th}$ level of $SPG$.
The last paired stubs are marked by thick lines, the brother stubs
by dashed lines. In $a)$ the $(l-1)^{\rm st}$ level is completed,
in $b)$ the pairing with a new node is described, in $c)$ the
pairing within the $(l-1)^{\rm st}$ level is described, and in
$d)$ the pairing with already existing node at $l^{\rm th}$ level
is described.}{dnziva_f1}

Because $ Z_{x,l}^{\smallsup{N}}$ is the number of free stubs at
level $l$ after the pairing of the first $x$ stubs, one would
expect that
    \eq
    \label{recZ2_}
    Z_{x,l}^{\smallsup{N}}\sim
    \sum_{i=1}^{x} X_{i,l-1}^{\smallsup{N}},
    \en
where $\sim$ denotes that we have an uncontrolled error term.
Indeed, the intuition behind (\ref{recZ2_}) is that loops or
cykels should be rare for small $l$.
Furthermore, when $M^{\smallsup{N}}_{l-1}$ is much smaller than
$N$, then the law of $X_{i,l-1}^{\smallsup{N}}$ should be quite
close to the law $G^{\smallsup{N}}$, which, in turn, by Lemma
\ref{lem-Fbd} is close to $G$. If $X_{i,l-1}^{\smallsup{N}}$ would
have distribution $G(x)$, then we could use the theory of sums of
random variables with infinite expectation, as well as extreme
value theory, to obtain the inequalities of
Proposition~\ref{dnziva_P1}.

In order to make the above estimates rigorous, we use upper and
lower bounds. We note that the right-hand side of (\ref{recZ2_})
is a valid upper bound for  $Z_{x,l}^{\smallsup{N}}$. We show
below that $X_{i,l-1}^{\smallsup{N}}$ have the same law, and we wish
to apply Lemma \ref{lem-sums}(i). For this, we need to control the
law $X_{i,l-1}^{\smallsup{N}}$, for which we use Lemma
\ref{lem-stochbds} to bound each $X_{i,l-1}^{\smallsup{N}}$ from
above by a random variable with law
$\overline{G}^{\smallsup{N}}_{\sM}$. This coupling makes sense
only on the {\em good event} where
$\overline{G}^{\smallsup{N}}_{\sM}$ is sufficiently close to
$G$.

For the lower bound, we have to do more work. The basic idea from
the theory of sums of random variables with infinite mean is that
the sum has the same order as the maximal summand. Therefore, we
bound from below
    \eq
    \label{Zlb2}
    Z_{x,l}^{\smallsup{N}}\geq
    \underline Z_{x,l}^{\smallsup{N}}-x.
    \en
where
    \eq
    \label{underZdef}
    \underline Z_{x,l}^{\smallsup{N}}= \max_{1\le i\le x} X_{i,l-1}^{\smallsup{N}}.
    \en
However, this lower bound is only valid when the chosen stub
is not part of the shortest path graph up to that point.
We show in Lemma \ref{dnziva_L3} below that the chosen stub has
label 1 when $\underline Z_{x,l}^{\smallsup{N}}
>2M_{l-1}^{\smallsup{N}}.$
In this case, (\ref{Zlb2}) follows since the $x-1$ remaining stubs
can `eat up' at most $x-1\leq x$ stubs. To proceed with the lower
bound, we bound $(X_{1,l-1}^{\smallsup{N}}, \ldots,
X_{x,l-1}^{\smallsup{N}})$ stochastically from below, using Lemma
\ref{lem-stochbds}, by an i.i.d.\ sequence of random variables
with laws $\underline{G}^{\smallsup{N}}_{\sM}$, where $M$ is
chosen appropriately and serves as an upper bound on the number of
stubs with label 3. Again on the {\it good event},
$\underline{G}^{\smallsup{N}}_{\sM}$ is sufficiently close to
$G$. Therefore, we are now faced with the problem of studying
the maximum of a number of random variables with a law close to
$G$. Here we can use Lemma \ref{lem-sums}(ii), and we conclude
in the proof of Proposition~\ref{dnziva_P1}(a) that $\underline Z_{x,l}^{\smallsup{N}}$ is to leading order equal
to $x^{\kappa}$, when $x=Z_{l-1}^{\smallsup{N}}\wedge
N^{\frac{1-\vep/2}{\kappa(\tau-1)}}.$ For this choice of $x$, we
also see that $\underline Z_{x,l}^{\smallsup{N}}$ is of bigger
order than $M_{l-2}^{\smallsup{N}},$ so that the basic assumption
in the above heuristic is satisfied. This completes the overview
of the proof.

We now state and prove the Lemmas~\ref{lem-barFbd(i)},
\ref{dnziva_L3}. The proof of Proposition~\ref{dnziva_P1} then follows
in Section \ref{sec-mainproofapp}.
We define the good event
mentioned above by
\eq
\label{gerard9}
  F_{\vep, \sM}
 = \bigcap_{x=1}^{N^{\alpha}}
 \left\{[1-2N^{-h}][1-G(x)]
    \leq 1-\overline{G}^{\smallsup{N}}_{\sM}(x)
    \leq 1-\underline{G}^{\smallsup{N}}_{\sM}(x)
    \leq [1+2N^{-h}] [1-G(x)]\right\}.
\en
The following lemma says that for $M\le N^{\alpha}$ the probability
of the good event is close to one.
\begin{lemma}
\label{lem-barFbd(i)} Let $F$ satisfy Assumption~\ref{ass-gamma}.
Then, for $\vep>0$ sufficiently small,
$$
\prob(F_{\vep,\sN^{\alpha}}^c) \leq N^{-h},\qquad
\mbox{ for large }N.
$$
\end{lemma}

\noindent {\bf Proof.}
Due to Lemma~\ref{lem-Fbd} it suffices to show that for  $\vep$ small enough, and $N$ sufficiently we have
    \eq
    F_{\vep, \sN^{\alpha}}^c \subset \FN^c.
    \en
We will prove the equivalent statement that
\eq
\label{dnziva_32}
\FN\subset F_{\vep, \sN^{\alpha}}.
\en
It follows from~(\ref{barG})
and~(\ref{underbarG}) that for every $M$ and $x$
    \eq
    1-\underline{G}^{\smallsup{N}}_{\sM}(x)\leq
    1-G^{\smallsup{N}}(x)\leq 1-\overline{G}^{\smallsup{N}}_{\sM}(x),
    \en
and, in particular, that for $M\leq N^{\alpha}$,
    \eq
    [1-\overline{G}^{\smallsup{N}}_{\sM}(x)]
    -[1-\underline{G}^{\smallsup{N}}_{\sM}(x)]
    \leq \frac{M}{L_{\sss N}-M}\leq C N^{\alpha-1}.
    \en
Then we use (\ref{1-F,1-G bound}(b)) to obtain that for all $x\le
N^{\alpha}$, $\vep$ small enough, and $N$ sufficiently large,
    \eq
    \begin{array}{l}
    C N^{\alpha-1}\leq N^{\alpha-1+h}
    =N^{\frac{1-\vep^5}{\tau-1}-1+\vep^6}
    <N^{-2\vep^6}N^{\frac{1-\vep^5}{\tau-1}(2-\tau-\vep^6)}\\[10pt]
    =N^{-2h}N^{\alpha(2-\tau-h)}
    \le N^{-2h}x^{2-\tau-h}
    \le N^{-h}[1-G(x)].
    \end{array}
    \en
Therefore, for $M\leq N^{\alpha}$ and with the above choices of
$\vep$, $\alpha$ and $h$, we have on $\FN$,
$$
\begin{array}{rl}
[1-\underline{G}^{\smallsup{N}}_{\sM}(x)]
\le&1-G^{\smallsup{N}}(x)
+[1-\overline{G}^{\smallsup{N}}_{\sM}(x)]
    -[1-\underline{G}^{\smallsup{N}}_{\sM}(x)]
\le[1+2N^{-h}][1-G(x)],\\[15pt]
[1-\underline{G}^{\smallsup{N}}_{\sM}(x)]
\ge&1-G^{\smallsup{N}}(x)
-[1-\overline{G}^{\smallsup{N}}_{\sM}(x)]
    +[1-\underline{G}^{\smallsup{N}}_{\sM}(x)]
\ge[1-2N^{-h}][1-G(x)],
\end{array}
$$
i.e. we have~(\ref{gerard9}), so that indeed $\FN \subset
F_{\vep, \sN^{\alpha}}$ . \qed

For the coupling of $X_{i,l-1}^{\smallsup{N}}$ with the random
variables with laws $\underline{G}^{\smallsup{N}}_{\sM}(x)$ and
$\overline{G}^{\smallsup{N}}_{\sM}(x)$ we need the following
lemma.

\begin{lemma}
\label{dnziva_L3} For any $l\ge 1$ there are at most
$2M_l^{\smallsup{N}}$ stubs with label $3$ in $SPG_{l+1}$, while the
number of stubs with label 2 is equal to $Z_{l+1}^{\smallsup{N}}$.
\end{lemma}
\noindent{\bf Proof.} The proof is by induction on $l$. There are
$Z_1^{\smallsup{N}}$ free stubs in $SPG_1$. Some of these stubs
will be paired with stubs with label 2 or 3, others will be paired
to stubs with label 1 (see Figure~\ref{dnziva_f1}). This gives us
at most $2Z_1^{\smallsup{N}}$ stubs with label $3$ in $SPG_2$.
This initializes the induction. We next advance the induction.
Suppose that for some $l\ge1$ there are at most
$2M_l^{\smallsup{N}}$ stubs with label $3$ in $SPG_{l+1}$. There are
$Z_{l+1}^{\smallsup{N}}$ free stubs (with label $2$) in $SPG_{l+1}$.
Some of these stubs will be paired with stubs with label 2 or 3,
others will be linked with stubs with label 1 (again see
Figure~\ref{dnziva_f1}). This gives us at most
$2Z_{l+1}^{\smallsup{N}}$ new stubs with label $3$ in $SPG_{l+2}$.
Hence the total number of these stubs is at most
$2M_l^{\smallsup{N}}+2Z_{l+1}^{\smallsup{N}}=2M_{l+1}^{\smallsup{N}}$.
This advances the induction hypothesis, and proves the claim. \qed

\vskip0.5cm


\subsection{The proof of Proposition ~\ref{dnziva_P1}}
\label{sec-mainproofapp} We state and prove some
consequences of the event $\hat F_{m,k}(\varepsilon)$, defined in (\ref{dnziva_2}). We refer to the remark, following
definition (\ref{dnziva_2}), to explain where we use these consequences.
\begin{lemma}
\label{dnziva_L1} The event $\hat F_{m,k}(\varepsilon)$ implies, for sufficiently large $N$,
the following bounds:
\begin{equation}
\label{dnziva_3}
\begin{array}{rl}
(a)&\qquad M_{k-1}^{\smallsup{N}}<N^{\frac{1-3\eps^4/4}{\kappa(\tau-1)}},\\[5pt]
(b)&\qquad \mbox{for any }\delta>0, \, N^{-\delta}\le k^{-3},\\[5pt]
(c)&\qquad \kappa^{k-1}(\varepsilon-\varepsilon^3) \le
\log\left(Z_{k-1}^{\smallsup{N}}\right)
\le \kappa^{k-1}(\varepsilon^{-1}+\varepsilon^3), \quad \mbox{for}\quad  k-1\ge m,   \\[5pt]
(d)&\qquad M_{k-1}^{\smallsup{N}}\leq 2Z_{k-1}^{\smallsup{N}}\quad \mbox{for}\quad  k-1\ge m .
\end{array}
\end{equation}
\end{lemma}

\noindent {\bf Proof. } Assume that~(\ref{dnziva_2}(a)-(d)) holds.
We start by showing (\ref{dnziva_3}(b)), which is evident if
we show the following claim:
    \begin{equation}
    \label{dnziva_7}
    k\le\frac{\log\left(\frac{1-\eps^2}{\eps(\tau-1)}\log N\right)}{\log\kappa},
    \end{equation}
for $N$ large enough. In order to prove (\ref{dnziva_7}), we note
that if $k\in{\cal T}^{\smallsup{N}}_m(\vep) $ then, due to
definition~(\ref{def-TmN}),
    \begin{equation}
    \label{dnziv3}
    \kappa^{k-m}
    \le\frac{1-\eps^2}{\tau-1}
    \frac{\log N}{\log\left( Z_m^{\smallsup{N}}\right)}<\frac{1-\eps^2}{\eps(\tau-1)}\kappa^{-m}\log N,
    \end{equation}
where the latter inequality follows from $Y_m^{\smallsup{N}}>\eps$
and (\ref{Yidef2_}). Multiplying by $\kappa^{m}$ and taking
logarithms on both sides yields (\ref{dnziva_7}).

We now turn to (\ref{dnziva_3}(a)). Since
    $$
    M_{k-1}^{\smallsup{N}}=\sumu\limits_{l=1}^{k-1}Z_l^{\smallsup{N}}
    \le k\max\limits_{1\le l\le k-1} Z_{l}^{\smallsup{N}},
    $$
the inequality~(\ref{dnziva_3}(a)) follows when we show that for
any $l\le k-1$,
    \begin{equation}
    \label{dnziva_3_}
    Z_{l}^{\smallsup{N}}
    \le N^{\frac{1-\eps^4}{\kappa(\tau-1)}}.
    \end{equation}
Observe that for $l<m$ we have, due to~(\ref{dnziva_2}(c)) and
(\ref{dnziva_2}(d)), for any $\vep>0$ and $m$ fixed and by taking
$N$ sufficiently large,
    \eq
    Z_{l}^{\smallsup{N}}
    \leq M_{m}^{\smallsup{N}}\leq 2 Z_{m}^{\smallsup{N}}
    \leq 2 e^{\kappa^m \vep^{-1}}<N^{\frac{1-\eps^4}{\kappa(\tau-1)}}.
    \en
Consider $m\le l\le k-1$. Due to~(\ref{Yidef2_}),
inequality (\ref{dnziva_3_}) is equivalent to
    \begin{equation}
    \label{dnziva_4}
    \kappa^{l+1}
    Y_{l}^{\smallsup{N}}\le
    \frac{1-\eps^4}{\tau-1}\log N.
    \end{equation}
To obtain~(\ref{dnziva_4}) we will need two inequalities.
Firstly,~(\ref{dnziva_2}(a)) and $l+1\le k$ imply that
    \begin{equation}
    \label{dnziva_5}
    \kappa^{l+1}
    Y_{m}^{\smallsup{N}}\le
    \frac{1-\varepsilon^2}{\tau-1}\log N.
    \end{equation}
Secondly, (\ref{dnziva_5}) and~(\ref{dnziva_2}(c)) imply that
    \begin{equation}
    \label{dnziva_6}
    \kappa^{l+1}
    \le\frac{1-\eps^2}{\eps(\tau-1)}\log N.
    \end{equation}
Given~(\ref{dnziva_5}) and~(\ref{dnziva_2}(b)),
we obtain, when $Y_{m}^{\smallsup{N}}\geq \vep$, and for $m\le l
\leq k-1$,
    \begin{equation}
    \begin{array}{rl}
    \kappa^{l+1}
    Y_{l}^{\smallsup{N}}
    &\le
    \kappa^{l+1}
    (Y_{m}^{\smallsup{N}}+\varepsilon^3)
    \le
    \kappa^{l+1}
    Y_{m}^{\smallsup{N}}(1+\varepsilon^2)\\
    &\le
    \frac{(1-\varepsilon^2)(1+\varepsilon^2)}{\tau-1}\log N
    =
    \frac{1-\eps^4}{\tau-1}\log N.
    \end{array}
    \end{equation}
Hence we have~(\ref{dnziva_4}) or equivalently~(\ref{dnziva_3_})
for $m\le l\le k-1$.

The bound in~(\ref{dnziva_3}(c)) is an immediate consequence
of~(\ref{Yidef2_}) and~(\ref{dnziva_2}(b,c)) that imply for
$k-1>m$,
    $$
    \varepsilon-\varepsilon^3
    \le Y_{k-1}^{\smallsup{N}}\le\varepsilon^{-1}+\varepsilon^3.
    $$

We complete the proof by establishing (\ref{dnziva_3}(d)). We use
induction to prove that for all $l\geq m$, the bound
$M_{l}^{\smallsup{N}}\leq 2Z_{l}^{\smallsup{N}}$ holds. The
initialization of the induction hypothesis for $l=m$ follows from
(\ref{dnziva_2}(d)). So assume that for some $m\leq  l <k-1$ the
inequality $M_{l}^{\smallsup{N}}\leq 2Z_{l}^{\smallsup{N}}$ holds,
then
    \eq
    M_{l+1}^{\smallsup{N}}=Z_{l+1}^{\smallsup{N}}+M_{l}^{\smallsup{N}}
    \leq Z_{l+1}^{\smallsup{N}}+2Z_{l}^{\smallsup{N}},
    \en
so that it suffices to bound $2Z_{l}^{\smallsup{N}}$ by
$Z_{l+1}^{\smallsup{N}}$ . We note that $\hat
F_{m,k}(\varepsilon)$ implies that
    \eq
    |Y_{l+1}^{\smallsup{N}}-Y_{l}^{\smallsup{N}}|
    \leq |Y_{l+1}^{\smallsup{N}}-Y_{m}^{\smallsup{N}}|
    +|Y_{l}^{\smallsup{N}}-Y_{m}^{\smallsup{N}}|
    \leq 2\eps^3\leq 3\eps^2 Y_{l+1}^{\smallsup{N}}.
    \en
Therefore,
    \eq
    2Z_{l}^{\smallsup{N}}= 2e^{\kappa^{l}Y_{l}^{\smallsup{N}}}
    \leq 2e^{(1+3\eps^2)\kappa^{l} Y_{l+1}^{\smallsup{N}}}
    = 2\big(Z_{l+1}^{\smallsup{N}}\big)^{(1+3\eps^2) \kappa^{-1}}\leq Z_{l+1}^{\smallsup{N}},
    \en
when $\eps>0$ is so small that $\omega=(1+3\eps^2) \kappa^{-1}< 1$
and where we take $m$ large enough to ensure that for $l\geq m$,
the lower bound
$Z_{l+1}^{\smallsup{N}}=\exp\{\kappa^{l+1}Y_{l+1}^{\smallsup{N}}
\}>\exp\{\kappa^{l+1}\vep\}>2^{\frac1{1-\omega}}$ is satified.
\qed \vskip0.5cm

\noindent {\bf Proof of Proposition~\ref{dnziva_P1}(a).} Recall
that $\alpha=\frac{1-\vep^5}{\tau-1}$. We write \eq
\label{dnziva_34}
\begin{array}{ll}
    \shift\shift\prob\Big(
    \hat F_{m,l}(\varepsilon) \cap \big\{
      Z_{x,l}
      \geq \left(l^3x\right)^{\kappa+c_{\gamma}\gamma(x)}\big\}
    \Big)
&\le \prob_{{\sss N^{\alpha}}} \Big(\hat F_{m,l}(\varepsilon) \cap
\big\{
      Z_{x,l}
      \geq \left(l^3x\right)^{\kappa+c_{\gamma}\gamma(x)}\big\}
    \Big) +\prob(F_{\vep, \sN^{\alpha}}^c)\\[15pt]
&\le \prob_{{\sss N^{\alpha}}} \Big(\hat F_{m,l}(\varepsilon) \cap
\big\{
      Z_{x,l}
      \geq \left(l^3x\right)^{\kappa+c_{\gamma}\gamma(x)}
      \big\}\Big)
+l^{-3},
\end{array}
\en where $\prob_{\sM}$ is the conditional probability given that
$F_{\vep, \sM}$ holds, and where we have used
Lemma~\ref  {lem-barFbd(i)}  with $N^{-h}<l^{-3}$. It remains to bound the
first term on the right-hand side  of~(\ref{dnziva_34}). For this bound we
aim to use Lemma  \ref{lem-sums}. Clearly we have
    \eq
    Z_{x,l}^{\smallsup{N}}
    \le \sumu\limits_{i=1}^x X_{i,l-1}^{\smallsup{N}},
    \en
because loops and cycles can occur (in Figure~\ref{dnziva_f1} only
the case $b)$ contributes to $Z_{x,l}^{\smallsup{N}}$,  the cases
$c)$ and $d)$ do not contribute). Since the free stubs of
$SPG_{l-1}$ are exchangeable, each free stub will choose any stub
with label unequal to 3 with the same probability. Therefore, all
$X_{i,l-1}^{\smallsup{N}}$ have the same law which we denote by
$H^{\smallsup{N}}$. Then we observe that due to~(\ref{bdsums}),
$X_{i,l-1}^{\smallsup{N}}$ can be coupled with
$\overline{X}^{\smallsup{N}}_{i,l-1}$ having law
$\overline{G}^{\smallsup{N}}_{\sM}$, where $M$ is equal to the
number of stubs with label $3$ at the moment we generate
$X_{i,l-1}^{\smallsup{N}}$, which is at most the number of stubs
with label $3$ in $SPG_{l}$ plus $1$. The last number is due to
Lemma~\ref{dnziva_L3} at most $2M_{l-1}^{\smallsup{N}}+1$. By Lemma \ref{dnziva_L1}(a), we have
that
    \eq
    2M_{l-1}^{\smallsup{N}}+1
    \leq 2N^{\frac{1-3\eps^4/4}{\kappa(\tau-1)}}+1\leq N^{\frac{1-\vep^5}{\tau-1}}= N^{\alpha},
    \en
and hence, due to~(\ref{bdsums}), we can take as the largest
possible number $M=N^{\alpha}$.
We now verify whether we can apply Lemma \ref{lem-sums}(i).
Observe that $x\le N^{\frac{1-\vep/2}{\kappa(\tau-1)}}$
so that for N large and each $c_{\gamma}$, we have
\eq
y=(l^3 x)^{\kappa+c_{\gamma}\gamma(x)}<N^{\alpha},
\en
since by (\ref{dnziva_7}), we can bound $l$ by a double logarithm.
Hence (\ref{lem-sums-cond1}) holds, because we condition on $F_{\vep, \sN^{\alpha}}$.
We therefore can apply Lemma \ref{lem-sums}(i),
with ${\hat S}_x^{\smallsup{N}}=\sum_{i=1}^x\overline{X}^{\smallsup{N}}_{i,l-1}$,
${\hat H}^{\smallsup{N}}=\overline{G}^{\smallsup{N}}_{\sN^{\alpha}}$,
 and obtain also using the upper  bound in (\ref{gerard1}),
\begin{eqnarray}
&&\prob_{{\sss N^{\alpha}}} \Big(\hat F_{m,l}(\varepsilon) \cap
\big\{
Z_{x,l}
\geq \left(l^3x\right)^{\kappa+c_{\gamma}\gamma(x)}
\big\}\Big)\le  b' x [1+2N^{-h}]\big[1-G\big(y)\big]\nn\\
&&\qquad\leq  2b' x  y^{-\kappa^{-1}+K_{\tau}\gamma(y)}
= 2b' x  (l^3 x)^{(-\kappa^{-1}+K_{\tau}\gamma(y))(\kappa+c_{\gamma}\gamma(x))}
\le bl^{-3},
\label{teveel}
\end{eqnarray}
if we show that
\eq
\label{ger1}
c_{\gamma}\gamma(x)\left(-\kappa^{-1}+K_{\tau}\gamma(y)\right)+\kappa K_{\tau}\gamma(y)<0.
\en
Inequality (\ref{ger1}) holds, because
$\gamma(y)=(\log y)^{\gamma-1},\, \gamma\in [0,1)$, can be made arbitrarily small by taking $y$ large, which
follows from (\ref{dnziva_3}(c)): $l^3 x\ge l^3 \exp\{\kappa^m\vep/2\}$,
and because $m$ can be taken large.
 \qed \vskip0.5cm

\noindent{\bf Proof of Proposition~\ref{dnziva_P1}(b).} Similarly to
(\ref{dnziva_34}), we have \eq \label{dnziva_35}
\begin{array}{l}
\prob\Big(\hat F_{m,l}(\varepsilon)\cap
      \big\{
      Z_{x,l}^{\smallsup{N}}
      \leq \left(\frac x{l^3}\right)^{\kappa-c_{\gamma}\gamma(x)}\big\}
    \Big)
\leq \prob_{\sss N^{\alpha}} \Big(\hat F_{m,l}(\varepsilon)\cap
      \big\{
      Z_{x,l}^{\smallsup{N}}
      \leq \left(\frac x{l^3}\right)^{\kappa-c_{\gamma}\gamma(x)}\big\}
    \Big)+l^{-3},
\end{array}
\en and it remains to bound the first term on the right-hand side
of~(\ref{dnziva_35}). Recall that
$$
\underline Z_{x,l}^{\smallsup{N}}= \max_{1\le i\le x}
X_{i,l-1}^{\smallsup{N}},
$$
where, for $1\le i \le x$, $X_{i,l-1}^{\smallsup{N}}$ is the
number of brother stubs of a stub attached to the $i^{\rm th}$ free
stub of $SPG_{l-1}$. Suppose we can bound the first term on the
right-hand side of~(\ref{dnziva_35}) by $bl^{-3}$, when
$Z_{x,l}^{\smallsup{N}}$ is replaced by $\underline
Z_{x,l}^{\smallsup{N}}$ after adding an extra factor 2, e.g.,
suppose that
    \[
    \prob_{\sss N^{\alpha}}
    \Big(\hat F_{m,l}(\varepsilon)\cap
      \big\{
      \underline Z_{x,l}^{\smallsup{N}}
      \leq 2\left(\smfrac x{l^3}\right)^{\kappa-c_{\gamma}\gamma(x)}\big\}
    \Big)\le bl^{-3}.
    \]
Then we bound
    \eqalign
    \prob_{\sss N^{\alpha}}
    &\Big(\hat F_{m,l}(\varepsilon)\cap
      \big\{
      Z_{x,l}^{\smallsup{N}}
      \leq \left(\smfrac x{l^3}\right)^{\kappa-c_{\gamma}\gamma(x)}\big\}
    \Big)
    \label{splitPNs}\\
    &\qquad\le
    \prob_{\sss N^{\alpha}}
    \Big(\hat F_{m,l}(\varepsilon)\cap
      \big\{
      \underline Z_{x,l}^{\smallsup{N}}
      \leq 2\left(\smfrac x{l^3}\right)^{\kappa-c_{\gamma}\gamma(x)}\big\}
    \Big)\nn\\
    &\qquad\quad+\prob_{\sss N^{\alpha}}
    \Big(\hat F_{m,l}(\varepsilon)\cap
      \big\{
      Z_{x,l}^{\smallsup{N}}
      \leq \left(\smfrac x{l^3}\right)^{\kappa-c_{\gamma}\gamma(x)}\big\}
      \cap \big\{
      \underline Z_{x,l}^{\smallsup{N}}
      >2\left(\smfrac x{l^3}\right)^{\kappa-c_{\gamma}\gamma(x)}\big\}
    \Big).\nn
    \enalign
    By assumption, the first term is bounded by $bl^{-3}$, and we must bound the
second term. We will prove that the second term in
(\ref{splitPNs}) is equal to 0.

For $x$ sufficiently large we obtain from $l\leq
C\log{x}$, $\kappa>1$, and $\gamma(x)\to 0$,
    \eq
    \label{powerxlb}
    2\left(\frac x{l^3}\right)^{\kappa-c_{\gamma}\gamma(x)}
    >6 x.
    \en
Hence for $x=Z_{l-1}^{\smallsup{N}}>(\vep-\vep^3) \kappa^{l-1}$, it follows from Lemma
\ref{dnziva_L1}(d), that $ \underline
Z_{x,l}^{\smallsup{N}}>2\left(\frac
x{l^3}\right)^{\kappa-c_{\gamma}\gamma(x)}$ induces
    \eq
    \underline Z_{x,l}^{\smallsup{N}}>6Z_{l-1}^{\smallsup{N}}\geq 2M_{l-1}^{\smallsup{N}}+2Z_{l-1}^{\smallsup{N}}.
    \en
On the other hand, when $x=N^{\frac{(1-\vep/2)}{\kappa(\tau-1)}}<Z_{l-1}^{\smallsup{N}}$, then,
by Lemma \ref{dnziva_L1}(a), and where we use again  $l\leq
C\log{x}$, $\kappa>1$, and $\gamma(x)\to 0$,
\begin{eqnarray}
&&\underline Z_{x,l}^{\smallsup{N}}
\ge
2\left(\smfrac x{l^3}\right)^{\kappa-c_{\gamma}\gamma(x)}
= 2\left(
\frac{N^{\smfrac{1-\vep/2}{\kappa(\tau-1)}}}
{l^3}
\right)^{\kappa-c_{\gamma}\gamma(x)}\nn\\
&&\qquad
> 2 N^{\frac{1-3\eps^4/4}{\kappa(\tau-1)}}+2N^{\frac{1-\vep/2}{\kappa(\tau-1)}}
>2 M^{(N)}_{l-1}+2x.
\end{eqnarray}
We conclude that in both cases we have that $\underline
Z_{x,l}^{\smallsup{N}}
\ge2M_{l-1}^{\smallsup{N}}+2x\ge 2M_{l-2}^{\smallsup{N}}+2x $.
We claim that the event $\underline
Z_{x,l}^{\smallsup{N}}>2M_{l-2}^{\smallsup{N}}+2x$ implies that
    \eq
    \label{eqiib_L1}
    Z_{x,l}^{\smallsup{N}}\ge \underline Z_{x,l}^{\smallsup{N}}-x.
    \en
Indeed, let $i_0\in\{1,\dots,N\}$ be the node
such that
    $$
    D_{i_0}=\underline Z_{x,l}^{\smallsup{N}}+1,
    $$
and suppose that $i_0\in SPG_{l-1}$. Then
$D_{i_0}$ is at most the total number of stubs
with labels $2$ and $3$, i.e., at most
$2M_{l-2}^{\smallsup{N}}+2x$. Hence $\underline
Z_{x,l}^{\smallsup{N}}<D_{i_0}\le
2M_{l-2}^{\smallsup{N}}+2x$, and  this is a contradiction with the
assumption that $\underline
Z_{x,l}^{\smallsup{N}}>2M_{l-2}^{\smallsup{N}}+2x$. Since by
definition $i_0\in SPG_{l}$, we conclude that
$i_0\in SPG_{l}\setminus SPG_{l-1}$, which is
equivalent to saying that the chosen stub with $\underline
Z_{x,l}^{\smallsup{N}}$ brother stubs had label 1. Then, on
$\underline Z_{x,l}^{\smallsup{N}}>2M_{l-2}^{\smallsup{N}}+2x$, we
have~(\ref{eqiib_L1}). Indeed, the one stub from level $l-1$
connected to $i_0$ gives us $\underline
Z_{x,l}^{\smallsup{N}}$ free stubs at level $l$ and the other
$x-1$ stubs from level $l-1$ can `eat up' at most $x$ stubs.

We conclude from the above that
    \eqalign
    &\prob_{\sss N^{\alpha}}
    \Big(\hat F_{m,l}(\varepsilon)\cap
      \big\{
      Z_{x,l}^{\smallsup{N}}
      \leq \left(\smfrac x{l^3}\right)^{\kappa-c_{\gamma}\gamma(x)}\big\}
      \cap \big\{
      \underline Z_{x,l}^{\smallsup{N}}
      >2\left(\smfrac x{l^3}\right)^{\kappa-c_{\gamma}\gamma(x)}\big\}
    \Big)\\
    &\qquad \leq
    \prob_{\sss N^{\alpha}}
    \Big(\hat F_{m,l}(\varepsilon)\cap
      \big\{
      Z_{x,l}^{\smallsup{N}}
      \leq \left(\smfrac x{l^3}\right)^{\kappa-c_{\gamma}\gamma(x)}\big\}
      \cap \big\{
      Z_{x,l}^{\smallsup{N}}
      >2\left(\smfrac x{l^3}\right)^{\kappa-c_{\gamma}\gamma(x)}-x\big\}
    \Big)=0,\nn
    \enalign
since (\ref{powerxlb}) implies that
    \[2\left(\frac x{l^3}\right)^{\kappa-c_{\gamma}\gamma(x)}-x\geq
    \left(\frac x{l^3}\right)^{\kappa-c_{\gamma}\gamma(x)}.
    \]

Below we prove in two steps that there exists $b$ such that
    \begin{equation}
    \label{gerard6}
    \prob_{\sss N^{\alpha}}\Big( \hat F_{m,l}(\varepsilon) \cap\{\underline Z_{x,l}^{\smallsup{N}}
    \leq 2\left(\smfrac x{l^3}\right)^{\kappa-c_{\gamma}\gamma(x)}\}\Big)\le bl^{-3}.
    \end{equation}
First we couple $\{X^{\smallsup{N}}_{i,l-1}\}_{i=1}^x$ with a
sequence of i.i.d.\ random variables $\{\underline
X^{\smallsup{N}}_{i,l-1}\}_{i=1}^x$ with law
$\underline{G}^{\smallsup{N}}_{\sss N^{\alpha}}$, such that
    \eq
    \label{dnziva_33}
    X^{\smallsup{N}}_{i,l-1}\ge \underline X^{\smallsup{N}}_{i,l-1},\qquad\qquad i=1,2,\dots,x,
    \en
and hence
    \eq
    \label{gerard8}
    \underline Z_{x,l}^{\smallsup{N}}
    \ge V_{x}^{\smallsup{N}}\stackrel{def}=\max\limits_{1\le i \le x}
    \underline X^{\smallsup{N}}_{i,l-1}.
    \en
Then we apply Lemma~\ref{lem-sums}(ii) with $\hat
X_i^{\smallsup{N}}=\underline{X}^{\smallsup{N}}_{i,l-1}$ and
$y=2 (x/l^3)^{\kappa-c_{\gamma}\gamma(x)}$.

We prove the fact that we can couple
$\{X^{\smallsup{N}}_{i,l-1}\}_{i=1}^x$ with a sequence of i.i.d.\
random variables $\{\underline X^{\smallsup{N}}_{i,l-1}\}_{i=1}^x$
with law $\underline{G}^{\smallsup{N}}_{\sss N^{\alpha}}$ by
induction on $x$. For $x=1$, the claim follows from Lemma
\ref{lem-stochbds} with $l=1$. Observe that for the $x^{\rm th}$
stub at level $l-1$, we sample a uniform stub from the stubs with
labels $1$ and $2$. Hence, due to~(\ref{bdsums}), and
conditionally on $\{X^{\smallsup{N}}_{i,l-1}\}_{i=1}^{x-1}$, we
can bound $X^{\smallsup{N}}_{x,l-1}$ from below by
$\underline{X}^{\smallsup{N}}_{x,l-1}$, which has law
$\overline{G}^{\smallsup{N}}_{\sM}$, where $M$ is equal to the
number of stubs with label $3$ at the moment we generate
$X_{x,l-1}^{\smallsup{N}}$. This number is bounded from above by
$2M_{l-2}^{\smallsup{N}}+2x \leq N^{\alpha}$ by
(\ref{dnziva_3}(a)) and the fact that $x\leq
N^{\frac{(1-\vep/2)}{\kappa(\tau-1)}} \leq \frac 14 N^{\alpha}$.
Indeed, the maximal possible value for $M$ corresponds to the
moment we sample $X_{x,l-1}^{\smallsup{N}}$, i.e., $M$ is at most
the number of stubs with label $3$ in $SPG_{l-1}$ plus $2x$ for
the pairing at most $x$ free stubs of $SPG_{l-1}$. Hence, due to
Lemma~\ref{dnziva_L3}, $M$ is at most
$2M_{l-2}^{\smallsup{N}}+2x\leq N^{\alpha}$, and due
to~(\ref{bdsums}) we can take $M=N^{\alpha}$. Therefore,
conditionally on $\{X^{\smallsup{N}}_{i,l-1}\}_{i=1}^{x-1}$, we
can bound $X^{\smallsup{N}}_{x,l-1}$ from below by
$\underline{X}^{\smallsup{N}}_{x,l-1}$, which has law
$\overline{G}^{\smallsup{N}}_{\sM}$, and this conditional coupling
is equivalent to fact that each component of
$\{X^{\smallsup{N}}_{i,l-1}\}_{i=1}^x$ can be bounded from below
by the components of $\{\underline
X^{\smallsup{N}}_{i,l-1}\}_{i=1}^x$, where $\{\underline
X^{\smallsup{N}}_{i,l-1}\}_{i=1}^x$ are i.i.d.\ copies with law
$\overline{G}^{\smallsup{N}}_{\sss N^{\alpha}}$.

We finally restrict to $x=Z_{l-1}^{\smallsup{N}}\wedge
N^{\frac{1-\vep/2}{\kappa(\tau-1)}}$. Note that $y=2(x/l^3)^{\kappa-c_{\gamma}\gamma(x)}\le N^{\alpha}$, so that
$F_{\alpha,{\sss N^{\alpha}}}$ holds, which in turn implies condition~(\ref{lem-sums-cond2}). We can therefore apply
Lemma~\ref{lem-sums}(ii) with $\hat
X_i^{\smallsup{N}}=\underline{X}^{\smallsup{N}}_{i,l-1}$,
$i=1,2,\dots,x$, $\hat
H^{\smallsup{N}}=\underline{G}^{\smallsup{N}}_{\sss N^{\alpha}}$,
and $y=2(x/l^3  )^{\kappa-c_{\gamma}\gamma(x)}$ to obtain from
(\ref{gerard8}),
\begin{eqnarray}
&&\prob_{\sss N^{\alpha}} \Big(\hat F_{m,l}(\varepsilon)\cap
      \big\{
      \underline Z_{x,l}^{\smallsup{N}}
      \leq 2\left(\smfrac x{l^3}\right)^{\kappa-c_{\gamma}\gamma(x)}\big\}
    \Big)\nn\\
    &&\qquad
    \leq \prob\left(\max_{1\le i\le x} \underline{X}^{\smallsup{N}}_{i,l-1}\le y\right)
\leq \left(1-[1-2N^{-h}][1-G\left(y\right)]\right)^x .
\label{boundxmacht}
\end{eqnarray}
From the lower bound of (\ref{gerard1}),
\begin{eqnarray}
\label{boven1-Gy}
[1-G(y)]\ge y^{-\kappa^{-1}-K_{\tau}\gamma(y)}=2^{\kappa-c_{\gamma}\gamma(x)}(x/l^3)^{(-\kappa^{-1}-K_{\tau}\gamma(y))(\kappa-c_{\gamma}\gamma(x))}
\ge \frac{l^3}{x},
\end{eqnarray}
because $x/l^3>1$ and
\[
\kappa^{-1}c_{\gamma}\gamma(x)-\kappa K_{\tau}\gamma(y)+c_{\gamma}K_{\tau}\gamma(x)\gamma(y)
\ge c_{\gamma}\kappa^{-1}\gamma(x)-\kappa K_{\tau}\gamma(y)\ge0,
\]
by choosing $c_{\gamma}$ large and using $\gamma(x)\geq \gamma(y)$.
Combining (\ref{boundxmacht}) and (\ref{boven1-Gy}) and taking $1-2N^{-h}>\frac12$, we conclude that
\eq
\left(1-[1-2N^{-h}][1-G\left(y\right)]\right)^x \leq \left(1-\frac{l^3}{2x}\right)^{x}\le e^{-l^3/2}\le l^{-3},
\en
because $l>m$ and $m$ can be chosen large. This yields (\ref{gerard6}) with $b=1$. \qed

In the proof of Proposition \ref{prop-weakconv2bsec}, in Section
\ref{sec-main idea}, we often    use a corollary of Proposition
\ref{dnziva_P1} that we formulate and prove below.

\begin{corr}
\label{dnziva_C1} Let $F$ satisfy Assumption~\ref{ass-gamma}. For
any $\vep>0$ sufficiently small, there exists an integer $m$ such
that such that for any $k>m$,
    \eq
    \label{dnziva_C1eq}
    \prob\left(\hat F_{m,k}(\varepsilon)\cap
        \{|Y_k^{\smallsup{N}}-Y_{k-1}^{\smallsup{N}}|
    >(\tau-2+\varepsilon)^{k(1-\gamma)}\}\right)\leq k^{-2},
    \en
for sufficiently large $N$.
\end{corr}
{\bf Proof.} We use that part $(a)$ and part $(b)$ of Proposition
\ref{dnziva_P1} together imply:
    \begin{equation}
    \label{dnziva_C2}
    \begin{array}{l}
    \prob\Big(\hat F_{m,k}(\varepsilon)\cap
      \big\{\left|\log\left(Z_k^{\smallsup{N}}\right)-\kappa \log\left(Z_{k-1}^{\smallsup{N}}\right)\right|
      \geq
      \kappa\log(k^3)
            +c_{\gamma}\gamma\left(Z_{k-1}^{\smallsup{N}}\right)\log\left(k^3Z_{k-1}^{\smallsup{N}}\right)
            \big\}
    \Big)\leq 2bk^{-3}.
    \end{array}
    \end{equation}
Indeed applying Proposition \ref{dnziva_P1}, with $l=k$ and
$x=Z_{k-1}^{\smallsup{N}}$, and hence $Z_{x,k}= Z_k$, yields:
    \begin{eqnarray}
    &&    \prob\Big(\hat F_{m,k}(\varepsilon)\cap
        \big\{
    Z_k^{\smallsup{N}} \ge (k^3 x)^{\kappa+c_{\gamma}\gamma(x)}\big\}
        \Big)\leq bk^{-3},\\
    &&    \prob\Big(\hat F_{m,k}(\varepsilon)\cap
        \big\{
    Z_k^{\smallsup{N}} \le (x/k^3 )^{\kappa-c_{\gamma}\gamma(x)}\big\}
        \Big)\leq bk^{-3},
    \end{eqnarray}
and from the identities
    \begin{eqnarray*}
    &&\{Z_k^{\smallsup{N}} \ge (k^3 x)^{\kappa+c_{\gamma}\gamma(x)}\}
    = \{\log(Z_k^{\smallsup{N}})-\kappa\log (Z_{k-1}^{\smallsup{N}})
    \ge \log((k^3 x)^{\kappa+c_{\gamma}\gamma(x)})-\kappa\log x\},\\
    &&\{Z_k^{\smallsup{N}} \le (x/k^3)^{\kappa-c_{\gamma}\gamma(x)}\}
    = \{\log(Z_k^{\smallsup{N}})-\kappa\log (Z_{k-1}^{\smallsup{N}})
    \le \log((x/k^3)^{\kappa+c_{\gamma}\gamma(x)})-\kappa\log x\},
    \end{eqnarray*}
we obtain (\ref{dnziva_C2}).

Applying (\ref{dnziva_C2}) and (\ref{Yidef2_}), we arrive at
    \eqalign
    \label{dnziva_15}
    &\prob\left(\hat F_{m,k}(\varepsilon)\cap
        \{|Y_k^{\smallsup{N}}-Y_{k-1}^{\smallsup{N}}|
    >(\tau-2+\varepsilon)^{k(1-\gamma)}\}\right)\\
    &\qquad\le \prob\left(\hat F_{m,k}(\varepsilon)\cap
        \big\{
        \kappa^{-k}\left[\kappa\log(k^3)
    +c_{\gamma}\gamma\left(Z_{k-1}^{\smallsup{N}}\right)\log\left(k^3Z_{k-1}^{\smallsup{N}}\right)\right]
    >(\tau-2+\varepsilon)^{k(1-\gamma)}\big\}\right)+2bk^{-3}.\nn
    \enalign

Observe that, due to Lemma~\ref{dnziva_L1}(c), and since $\gamma(x)=(\log x)^{\gamma-1}$,
where $0\le \gamma<1$, we have on $\hat F_{m,k}(\vep)$,
    $$
    \begin{array}{l}
    \shift\shift\kappa^{-k}\left[\kappa\log(k^3)
    +c_{\gamma}\gamma\left(Z_{k-1}^{\smallsup{N}}\right)\log\left(k^3Z_{k-1}^{\smallsup{N}}\right)\right]\\[10pt]
    =\kappa^{-k}\left[\kappa\log(k^3)
    +c_{\gamma}\left(\log\left(Z_{k-1}^{\smallsup{N}}\right)\right)^{\gamma-1}
    \left(\log(k^3)+\log\left(Z_{k-1}^{\smallsup{N}}\right)\right)\right]\\[10pt]
    \le\kappa^{-k}\left[\kappa\log(k^3)
    +c_{\gamma}\log(k^3)+c_{\gamma}\left(\log\left(Z_{k-1}^{\smallsup{N}}\right)\right)^\gamma\right]\\[10pt]
    \le\kappa^{-k}\left[(c_{\gamma}+\kappa)\log(k^3)
    +c_{\gamma}\left(\kappa^{k-1}(\vep^{-1}+\vep^3)\right)^\gamma\right]\\[10pt]
    \le\kappa^{-k(1-\gamma)}\left[\kappa^{-k\gamma}(c_{\gamma}+\kappa)\log(k^3)
    +c_{\gamma}\left(\kappa^{-1}(\vep^{-1}+\vep^3)\right)^\gamma\right]\le (\tau-2+\vep)^{k(1-\gamma)} ,
    \end{array}
    $$
because, for $k$ large, and since $\kappa^{-1}=\tau-2$,
    $$
    \begin{array}{rl}
    \left(\frac{\tau-2}{\tau-2+\varepsilon}\right)^{k(1-\gamma)}
    [\kappa^{-k\gamma}(c_{\gamma}+\kappa)\log(k^3)]&\le\frac12,\qquad
    \left(\frac{\tau-2}{\tau-2+\varepsilon}\right)^{k(1-\gamma)}
    c_{\gamma}\left(\kappa^{-1}(\vep^{-1}+\vep^3)\right)^\gamma\le\frac12.
    \end{array}
    $$
We conclude that the first term on the right-hand side of~(\ref{dnziva_15}) is $0$,
for sufficiently large $k$, and the second term is bounded by
$2bk^{-3}\leq k^{-2}$, and hence the statement of the
corollary follows. \qed

\subsection{Proof of Proposition~\ref{prop-weakconv2bsec}
and Proposition~\ref{prop-Tub}} \label{sec-main idea}

\vskip0.5cm \noindent {\bf Proof of
Proposition~\ref{prop-weakconv2bsec}(a).}
We have to show that
    \[
 \prob(\vep\leq Y_m^{\smallsup{i,N}}\leq \vep^{-1}, \max_{k\in {\cal T}^{\smallsup{i,N}}_m(\vep)}
    |Y_k^{\smallsup{i,N}}-Y_m^{\smallsup{i,N}}|>\vep^3)=
    o_{\sN,m,\vep}(1).
 \]
Fix $\varepsilon>0$,
such that $\tau-2+\vep<1$. Then, take $m=m_\vep$, such
that~(\ref{dnziva_2_}) holds, and increase $m$, if necessary,
until~(\ref{dnziva_C1eq}) holds.

We use the inclusion that (recall the definition of ${\cal
T}^{\smallsup{N}}_m(\vep)$ given in (\ref{def-TmN})),
    \eq
    \label{incluYs1}
    \{\max_{k\in{\cal T}^{\smallsup{N}}_m(\vep)}
    |Y_k^{\smallsup{N}}-Y_m^{\smallsup{N}}|>\vep^3\}
    \subset \Big\{\sum\limits_{k\in{\cal T}^{\smallsup{N}}_m(\vep)}\kern-1em
    |Y_k^{\smallsup{N}}-Y_{k-1}^{\smallsup{N}}|
    >\sum\limits_{k\geq m}(\tau-2+\varepsilon)^{k(1-\gamma)}\Big\}.
    \en
If the event on the right-hand side of (\ref{incluYs1}) holds, then there must be a
$k\in {\cal T}^{\smallsup{N}}_m(\vep)$ such that
$|Y_k^{\smallsup{N}}-Y_{k-1}^{\smallsup{N}}|>(\tau-2+\varepsilon)^{k(1-\gamma)}$,
and therefore
    \eq
    \label{incluYs2}
    \{\max_{k\in{\cal T}^{\smallsup{N}}_m(\vep)}
    |Y_k^{\smallsup{N}}-Y_m^{\smallsup{N}}|>\vep^3\}
    \subset \bigcup_{k\in{\cal T}^{\smallsup{N}}_m(\vep)}
    F_{m,k-1}\cap F_{m,k}^c,
    \en
where we denote
    \eq
     F_{m,k}=F_{m,k}(\varepsilon)=\bigcap\limits_{j=m+1}^{k}\left\{
    |Y_j^{\smallsup{N}}-Y_{j-1}^{\smallsup{N}}|
    \le (\tau-2+\varepsilon)^{j(1-\gamma)}\right\}.
    \en
Since (\ref{dnziva_2_}) implies that on $F_{m,k-1}$ we have $
|Y_j^{\smallsup{N}}- Y_m^{\smallsup{N}}|\leq\vep^3$, $m<j\leq
k-1$, we find,
    \eq
    \label{incluYs3}
    F_{m,k-1}\cap F_{m,k}^c\subset
    \{|Y_{l}^{\smallsup{N}}-Y_m^{\smallsup{N}}|\leq\vep^3, \forall l\,: m<l\le k-1\}
    \cap\Big\{|Y_k^{\smallsup{N}}-Y_{k-1}^{\smallsup{N}}|
   >(\tau-2+\varepsilon)^{k(1-\gamma)}\Big\}.
    \en

Take $N$ sufficiently large such that, by
Proposition~\ref{prop-caft},
    \eqalign
    \label{prop-caftused}
    \prob\left(M_{m}^{\smallsup{N}}>2Z_{m}^{\smallsup{N}}, \eps\leq Y_{m}^{\smallsup{N}}\leq \eps^{-1}\right)
    &\leq \prob\left(\exists l\leq m: Y_l^{\smallsup{N}}\neq Y_l\right)
    +\prob\left(M_{m}>2{\cal Z}_{m}, \eps\leq Y_{m}\leq \eps^{-1}\right)\nn\\
    &\le \prob\left(M_{m}>2{\cal Z}_{m}, \eps\leq Y_{m}\leq \eps^{-1}\right)+\varepsilon/4,
    \enalign
Next, we use that
    \eq
    \label{prop-caftusedb}
    \lim_{m\rightarrow\infty} \prob\left(M_{m}>2{\cal Z}_{m}, \eps\leq Y_{m}\leq \eps^{-1}\right)=0,
    \en
since $Y_l=(\tau-2)^l \log{{\cal Z}_l}$ converges a.s., so that when
$Y_{m}\geq \eps$ and $m$ is large, $M_{m-1}$ is much smaller than
$Z_{m}$, so that $M_m=M_{m-1}+{\cal Z}_m>2{\cal Z}_m$ has small probability, as $m$ is large.

Then we use (\ref{incluYs1})--(\ref{prop-caftusedb}), together
with (\ref{dnziva_2}), to derive that
    \eqalign
    \label{dnziva_1}
    &\prob\left(\vep\leq Y_m^{\smallsup{N}}\leq \vep^{-1},
    \max_{k\in{\cal T}^{\smallsup{N}}_m(\vep)}
        |Y_k^{\smallsup{N}}-Y_m^{\smallsup{N}}|>\vep^3\right)\\
    &\quad \le
    \prob\left(\vep\leq Y_m^{\smallsup{N}}\leq \vep^{-1},
    \max\limits_{k\in{\cal T}^{\smallsup{N}}_m(\vep)}
        |Y_k^{\smallsup{N}}-Y_m^{\smallsup{N}}|>\vep^3,\,
    Y_l^{\smallsup{N}}=Y_l, \,\forall l\leq m\right)
    +\prob\left(\exists l\leq m: Y_l^{\smallsup{N}}\neq Y_l\right)\nn\\
    &\quad \le
    \sum\limits_{k>m} \prob\left(
    F_{m,k-1}\cap F_{m,k}^c\cap \{k \in {\cal T}^{\smallsup{N}}_m(\vep)\}
    \cap \{\vep\leq Y_m^{\smallsup{N}}\leq \vep^{-1}\} \cap \{
    Y_l^{\smallsup{N}}=Y_l, \, \forall l\leq m\}\right)
    +\frac{\varepsilon}{2}\nn\\
         &\quad  \le
    \sum\limits_{k>m}
    \prob\left(\hat F_{m,k}(\varepsilon)\cap \Big\{
        |Y_k^{\smallsup{N}}-Y_{k-1}^{\smallsup{N}}|
    >(\tau-2+\varepsilon)^{k(1-\gamma)}\Big\}\right)+\eps<3\varepsilon/2,
    \nn
    \enalign
by Corollary~\ref{dnziva_C1}. \qed

\vskip0.5cm \noindent {\bf Proof of
Proposition~\ref{prop-weakconv2bsec}(b).} We first
show~(\ref{monZT}), then~(\ref{bdZT}). Due to
Proposition~\ref{prop-weakconv2bsec}(a), and using that
$\{Y_m^{\smallsup{N}}\leq \vep^{-1}\}$, we find
$$
Y_k^{\smallsup{N}}\le Y_m^{\smallsup{N}}+\eps^3\leq
Y_m^{\smallsup{N}}(1+\eps^2),
$$
apart from an event with probability $o_{\sN,m,\vep}(1)$,
for all $k\in{\cal
T}^{\smallsup{N}}_m$. By~(\ref{Yidef2_})
and because $k\in {\cal T}^{\smallsup{N}}_m $, this is equivalent
to
$$
Z_k^{\smallsup{N}}
\le\left(Z_m^{\smallsup{N}}\right)^{\kappa^{k-m}(1+\eps^2)} \le
N^{\frac{1-\eps^2}{\tau-1}(1+\eps^2)}=N^{\frac{1-\eps^4}{\tau-1}}.
$$

We next show~(\ref{bdZT}). Observe that $k\in{\cal
T}^{\smallsup{N}}_m$ implies that either $k-1\in{\cal
T}^{\smallsup{N}}_m$, or $k-1=m$. Hence, from $k\in{\cal
T}^{\smallsup{N}}_m$ and Proposition~\ref{prop-weakconv2bsec}(a),
we obtain,  apart from an event with probability $o_{\sN,m,\vep}(1)$,
    \eq
    Y_{k-1}^{\smallsup{N}}\geq Y_m^{\smallsup{N}}-\vep^3\geq \vep-\vep^3\geq \frac{\vep}{2},
    \en
for $\vep>0$ sufficiently small, and
    \eq
    Y_k^{\smallsup{N}}=Y_k^{\smallsup{N}}-Y_m^{\smallsup{N}}
    +Y_m^{\smallsup{N}}-Y_{k-1}^{\smallsup{N}}+Y_{k-1}^{\smallsup{N}}\ge Y_{k-1}^{\smallsup{N}}-2\eps^3
    \geq Y_{k-1}^{\smallsup{N}} (1-4\vep^2),
    \en
By~(\ref{Yidef2_}) this is equivalent to
$$
Z_k^{\smallsup{N}}
\ge\left(Z_{k-1}^{\smallsup{N}}\right)^{\kappa(1-4\vep^2)} \geq
Z_{k-1}^{\smallsup{N}},
$$
when $\vep>0$ is so small that $\kappa(1-4\vep^2)\geq 1$, since
$\tau\in(2,3)$, and $\kappa=(\tau-2)^{-1}$. \qed

\vskip0.5cm

\noindent {\bf Proof of Proposition~\ref{prop-Tub}.} We must show
that
    \eq
    \label{aimprop3.3}
    \prob(k\in \partial{\cal T}^{\smallsup{N}}_m(\vep),
    \vep\leq Y_m^{\smallsup{N}}\leq \vep^{-1},
    Z_{k+1}^{\smallsup{N}}\leq N^{\frac{1-\vep}{\tau-1}})
    =o_{\sN,m,\vep}(1),
    \en
where
\[
\{k\in \partial {\cal T}_m^{\smallsup{N}}\}
    =\{k\in {\cal T}_m^{\smallsup{N}}\}\cap \{k+1\not\in {\cal T}_m^{\smallsup{N}}\}.
\]
In the proof, we
will make repeated use of Propositions \ref{prop-weakconv2bsec}
and~\ref{prop-caft}, whose proofs are now complete.

According to the definition of $\hat
F_{m,k}(\varepsilon)$ in (\ref{dnziva_2}),
    \begin{eqnarray}
    &&\prob\big(\{k\in \partial{\cal T}^{\smallsup{N}}_m(\vep)\}\cap
    \{\vep\leq Y_m^{\smallsup{N}}\leq \vep^{-1}\}\cap \hat F_{m,k}(\varepsilon)^c\big)\\
    &&\qquad\le
    \prob(\vep\leq Y_m^{\smallsup{N}}\leq \vep^{-1},\max_{l\in{\cal T}^{\smallsup{N}}_m(\vep)}
        |Y_l^{\smallsup{N}}-Y_m^{\smallsup{N}}|>\vep^3)+\prob\left(
        \vep\leq Y_m^{\smallsup{N}}\leq \vep^{-1},M_{m}^{\smallsup{N}}>2Z_{m}^{\smallsup{N}}\right)\nn.
\label{gerard10}
    \end{eqnarray}
In turn Propositions~\ref{prop-weakconv2bsec}(a)
and~\ref{prop-caft}, as well as
(\ref{prop-caftused}--\ref{prop-caftusedb}) imply that both
probabilities on the right-hand side of (\ref{gerard10}) are $o_{\vep}(1)$, as first $N\to\infty$, and then $m\to \infty$. Therefore, it
suffices to show
\begin{eqnarray}
    &&\prob\big(\{k\in \partial{\cal T}^{\smallsup{N}}_m(\vep),
    \vep\leq Y_m^{\smallsup{N}}\leq \vep^{-1},
    Z_{k+1}^{\smallsup{N}}\leq N^{\frac{1-\vep}{\tau-1}}\}\cap \hat F_{m,k}(\varepsilon)\big)\nn\\
    &&\qquad =\prob\big(\{k+1\not\in {\cal T}^{\smallsup{N}}_m(\vep),
    Z_{k+1}^{\smallsup{N}}\leq N^{\frac{1-\vep}{\tau-1}}\}\cap \hat F_{m,k}(\varepsilon)\big)
    =o_{\sN,m,\vep}(1).
    \label{furtheraimprop3.3}
\end{eqnarray}

Let $x=N^{\frac{1-\vep/2}{\kappa(\tau-1)}}$, and define the event
$I_{\sN,k}=I_{\sN,k}(a)\cap I_{\sN,k}(b)\cap I_{\sN,k}(c)\cap I_{\sN,k}(d)$, where
    \begin{eqnarray}
    I_{\sN,k}(a)&=&\{M_{k-1}^{\smallsup{N}}< N^{\frac{1-3\eps^4/4}{\kappa(\tau-1)}}\},
    \label{dnziva_28a}\\
    I_{\sN,k}(b)&=&\{x\le Z_{k}^{\smallsup{N}}\},
    \label{dnziva_28b}\\
    I_{\sN,k}(c)&=&\{Z_{k}^{\smallsup{N}}\le N^{\frac{1-\vep^4}{\tau-1}}\},
    \label{dnziva_28c}\\
    I_{\sN,k}(d)&=&\{Z_{k+1}^{\smallsup{N}}\ge Z_{x,k+1}^{\smallsup{N}}-
    Z_{k}^{\smallsup{N}}\}.
    \label{dnziva_28d}
    \end{eqnarray}
We split
    \eqalign
    \label{A518}
    &\prob\big(\{k+1\not\in {\cal T}^{\smallsup{N}}_m(\vep), Z_{k+1}^{\smallsup{N}}\leq N^{\frac{1-\vep}{\tau-1}}\}\cap \hat F_{m,k}(\vep)\big)\\
    &\qquad = \prob\big(\{k+1\not\in {\cal T}^{\smallsup{N}}_m(\vep), Z_{k+1}^{\smallsup{N}}\leq N^{\frac{1-\vep}{\tau-1}}\}
    \cap \hat F_{m,k}(\vep)\cap I_{\sN,k}\big)\nn\\
    &\qquad\quad
    +\prob\big(\{k+1\not\in {\cal T}^{\smallsup{N}}_m(\vep),Z_{k+1}^{\smallsup{N}}\leq N^{\frac{1-\vep}{\tau-1}}\}\cap \hat F_{m,k}(\vep)
    \cap I_{\sN,k}^c\big).\nn
    \enalign
We claim that both probabilities are small, which would complete
the proof. We start to show that
    \eq
    \prob\big(\{k+1\not\in {\cal T}^{\smallsup{N}}_m(\vep),
    Z_{k+1}^{\smallsup{N}}\leq N^{\frac{1-\vep}{\tau-1}}\}\cap \hat F_{m,k}(\vep)\cap I_{\sN,k}\big)=o_{\sN,m,\vep}(1).
    \label{gerard4}
    \en
Indeed, by Lemma \ref{lem-bdsexp} and (\ref{monZT}),
    \eqalign
    &\prob\big(\{k+1\notin {\cal T}^{\smallsup{N}}_m(\vep)\}\cap
    Z_{k}^{\smallsup{N}}\geq N^{\frac{1-\vep}{\tau-1}}\}\cap \hat F_{m,k}(\vep)\cap I_{\sN,k}\big)\\
    &\qquad \leq \prob\big(\{k\in {\cal T}^{\smallsup{N}}_m(\vep)\}\cap
    \{\vep\leq Y_m^{\smallsup{N}}\leq \vep^{-1}\}\cap\{
    Z_{k}^{\smallsup{N}}\in [N^{\frac{1-\vep}{\tau-1}}, N^{\frac{1-\vep^4}{\tau-1}}]\}\big)+o_{\sN,m,\vep}(1)=o_{\sN,m,\vep}(1),\nn
    \label{remco1}
    \enalign
where $u=(\tau-1)^{-1}$. Therefore, we are left to deal with the
case where $Z_{k}^{\smallsup{N}}\leq N^{\frac{1-\vep}{\tau-1}}$.
For this, and assuming $I_{\sN,k}$, we can use
Proposition~\ref{dnziva_P1}(b) with
$x=N^{\frac{1-\vep/2}{\kappa(\tau-1)}}\leq Z_{k}^{\smallsup{N}}$
by $I_{\sN,k}(b)$, and $l=k+1$ to obtain that, {\bf whp},
    \eq
    \label{ger4}
    Z_{k+1}^{\smallsup{N}}\ge Z_{x,k+1}^{\smallsup{N}}-Z_{k}^{\smallsup{N}}
    \ge x^{\kappa(1-\vep/2)}-N^{\frac{1-\vep}{\tau-1}}
    = N^{\frac{(1-\vep/2)^2}{\tau-1}}-N^{\frac{1-\vep}{\tau-1}}
    > N^{\frac{1-\vep}{\tau-1}},
    \en
where we have used that when $k\in{\cal T}^{\smallsup{N}}_m(\vep)$
and $Y_m^{\smallsup{N}}>\vep$, then we have $k\le c\log\log N$,
for some $c=c(\tau,\vep)$, and hence, for $N$ large enough,
    $$
    (k+1)^{3(\kappa-c_{\gamma}\gamma(x))}x^{c_{\gamma}\gamma(x)}
    \le (k+1)^{3\kappa} x^{c_{\gamma}\gamma(x)}
    \le x^{\vep\kappa/2}.
    $$
This proves (\ref{gerard4}).

For the second probability on the right-hand side of (\ref{A518}) it suffices to prove that
    \begin{eqnarray}
    \label{gerard5}
    \prob\big(\{k+1\not\in {\cal T}^{\smallsup{N}}_m(\vep)
    \}\cap \hat F_{m,k}(\vep) \cap I_{\sN,k}^c\big)
    =o_{\sN,m,\vep}(1).
    \end{eqnarray}
In order to prove (\ref{gerard5}), we prove that (\ref{gerard5})
holds with $I_{\sN,k}^c$ replaced by each one of the four events
$I^c_{\sN,k}(a),\ldots,I^c_{\sN,k}(d)$. For the intersection with the
event $I^c_{\sN,k}(a)$, we apply Lemma~\ref{dnziva_L1}(a), which
states that $\hat F_{m,k}(\vep) \cap I^c_{\sN,k}(a)$ is an empty set.

It follows from~(\ref{def-TmN}) that if $k+1\not\in {\cal
T}^{\smallsup{N}}_m(\vep)$, then
\eq \label{dnziva_22}
\kappa^{k+1}Y_m^{\smallsup{N}}>\frac{1-\varepsilon^2}{\tau-1}\log
N. \en
If  $\hat F_{m,k}(\varepsilon)$ holds
then by definition (\ref{dnziva_2}), and Corollary
\ref{dnziva_C1}, {\bf whp},
    \eq
    \label{dnziva_23}
     Y^{\smallsup{N}}_k \ge Y^{\smallsup{N}}_{k-1}\ge
    Y^{\smallsup{N}}_m-\varepsilon^3\geq
    Y^{\smallsup{N}}_m(1-\varepsilon^2).
    \en
Hence, if $\hat F_{m,k}(\varepsilon)$ holds and $k+1\not\in {\cal
T}^{\smallsup{N}}_m(\vep)$, then, by
(\ref{dnziva_22})--(\ref{dnziva_23}), {\bf whp},
    \eq
    \begin{array}{l}
    \kappa\log(Z_k^{\smallsup{N}})=\kappa^{k+1}Y_k^{\smallsup{N}} \ge
    (1-\vep^2)\kappa^{k+1}Y^{\smallsup{N}}_m \geq
    \frac{(1-\vep^2)^2}{\tau-1} \log{N},
    \end{array}
    \en
so that, {\bf whp},
    \eq
    Z_k^{\smallsup{N}}
    \ge x=N^{\frac{1-\vep/2}{\kappa(\tau-1)}},
    \en
for small enough $\vep>0$ and sufficiently large $N$, i.e., we
have
\[
\prob(\{k+1\not\in {\cal T}^{\smallsup{N}}_m(\vep)
    \}\cap \hat F_{m,k}(\vep) \cap I^c_{\sN,k}(b))=o_{\sN,m,\vep}(1).
\]
From Proposition~\ref{prop-weakconv2bsec}(b) it is immediate that
\[
\prob(\{k+1\not\in {\cal T}^{\smallsup{N}}_m(\vep)
    \}\cap \hat F_{m,k}(\vep) \cap I^c_{\sN,k}(c))=o_{\sN,m,\vep}(1).
\]
Finally, recall that $Z_{x,k+1}^{\smallsup{N}}$ is the number of
constructed free stubs at level $k+1$ after pairing of the first $x$
stubs at level $k$. After the pairing of the remaining
$Z_{k}^{\smallsup{N}}-x$ stubs at level $k$ they can `eat up' at
most $Z_{k}^{\smallsup{N}}-x\le Z_{k}^{\smallsup{N}}$ stubs, so
that $I_{N,k}(d)$ holds with probability 1, and hence the event in
(\ref{gerard5})
 intersected with $I^c_{\sN,k}(d)$ has probability equal to $0$.

This completes the proof of (\ref{gerard5}) and hence the
Proposition. \qed

\section*{Acknowledgement}
We thank Eric Cator for useful discussions on uniform continuity
used in the proof of Lemma \ref{lem-bdsexp}. The work of RvdH and
DZ was supported in part by Netherlands Organisation for
Scientific Research (NWO). The work of all authors was performed
in part at the Mittag-Leffler Institute. The work of RvdH was
performed in part also at Microsoft Research.


\smallskip

\parindent0pt
Corresponding author:\\
G. Hooghiemstra\\
EEMCS, Delft University of Technology\\
P.O.Box 5031\\
2600GA Delft\\
The Netherlands\\
E-mail: G.Hooghiemstra@ewi.tudelft.nl

\end{document}